\documentclass[12pt]{article}
\usepackage{subfigure}
\usepackage{tikz}
\usepackage{pstricks, amssymb, xr, multind, multicol,textcomp}
\def\hfl#1{\smash{\mathop{\hbox to 10mm{\rightarrowfill}}\limits^{\textstyle
#1}}}
\newtheorem{proposition}[equation]{Proposition} 
 
\newtheorem{theorem}[equation]{Theorem} 
\newtheorem{exa}[equation]{Example} 
\newtheorem{ex}[equation]{Exercise} 
\newtheorem{s-ex}[equation]{Side-exercise} 
\newtheorem{exas}[equation]{Examples} 
\newtheorem{lemma}[equation]{Lemma} 
\newtheorem{sublemma}[equation]{Sublemma} 
\newtheorem{remar}[equation]{Remark} 
\newtheorem{remars}[equation]{Remarks} 
\newtheorem{nota}[equation]{Notation} 
\newtheorem{sremar}[equation]{Side-remark} 
\newtheorem{definitio}[equation]{Definition}

\newenvironment{remark}{\begin{remar} \rm }{\end{remar}}

\newenvironment{definition}{\begin{definitio} \rm }{\end{definitio}}

\newcommand{\eqprop}{F} 
\newcommand{\deleqprop}{G} 
\newcommand{\deleqpropb}{H}
\newcommand{\CONS}{{\cal C}} 
 
\newcommand{\KK}{\mathbb{K}}

\newcommand{\CA}{{\cal A}}

\newcommand{\TCM}{\tilde{C}_2(M)}
\newcommand{\TCE}{\tilde{C}_2(\EXT)}
 
\newcommand{\CB}{{\cal B}} 
\newcommand{\CaC}{{\cal C}} 
\newcommand{\CD}{{\cal D}} 
\newcommand{\CE}{{\cal E}}

\newcommand{\CI}{{\cal I}} 
 
\newcommand{\CL}{{\cal L}}

\newcommand{\CP}{{\cal P}} 
\newcommand{\CT}{{\cal T}} 
 
\newcommand{\CQ}{{\cal Q}}
\newcommand{\CS}{{\cal S}}
\newcommand{\EXT}{E}
\newcommand{\EXTone}{E^-}
\newcommand{\EXTtwo}{E^+}
\newcommand{\tauone}{\tau^-}
\newcommand{\tautwo}{\tau^+}
\newcommand{\Kone}{K^-}
\newcommand{\Ktwo}{K^+}

\newcommand{\Cone}{C^-}
\newcommand{\Ctwo}{C^+}
\newcommand{\Rone}{R^-}
\newcommand{\Rtwo}{R^+}
\newcommand{\EXTonestwo}{E}
\newcommand{\EXTonetwo}{E^{-,+}}
\newcommand{\BOF}{D^{-,+}}
\newcommand{\Mbet}{M}
\newcommand{\Aman}{A}

\newcommand{\fvarM}{X}
\newcommand{\svarM}{Y}

\newcommand{\cvarM}{X}
\newcommand{\sK}{s}

\newcommand{\ZZ}{\mathbb{Z}}

\newcommand{\thetap}{\overline{\theta}} 
 
\newcommand{\funcE}{f_E} 
\newcommand{\func}{f} 
\newcommand{\functE}{\tilde{f}_E} 
\newcommand{\RR}{\mathbb{R}} 
\newcommand{\QQ}{\mathbb{Q}} 
\newcommand{\CC}{\mathbb{C}} 
\newcommand{\NN}{\mathbb{N}}

\newcommand{\bp}{\noindent {\sc Proof: }} 
\newcommand{\eop}{\nopagebreak \hspace*{\fill}{$\square$} \medskip}

\newcommand{\tetra}{ \begin{tikzpicture}\useasboundingbox (-.2,0) rectangle (.8,.3);
\begin{scope}[yshift=-.1cm]
\draw (.6,0) -- (0,0) -- (.3,.5) -- (.6,0);
\draw [very thick, dotted] (.6,0) -- (.3,.2) (0,0) -- (.3,.2) -- (.3,.5);
\fill (.3,.5) circle (1.5pt) (0,0) circle (1.5pt) (.3,.2) circle (1.5pt) (.6,0) circle (1.5pt);\end{scope}
\end{tikzpicture}}
\newcommand{\tadpole}{\begin{tikzpicture} \useasboundingbox (.05,-.1) rectangle (1.5,.2);
\draw (.3,0) circle (.15) (.45,0) -- (1.15,0);
\fill (.45,0) circle (1.5pt) (1.15,0) circle (1.5pt);
\draw [dashed] (1.3,0) circle (.15);
\end{tikzpicture}}
\newcommand{\tadpoleplain}{\begin{tikzpicture} \useasboundingbox (.05,-.1) rectangle (1.5,.2);
\draw (.3,0) circle (.15) (1.3,0) circle (.15) (.45,0) -- (1.15,0);
\fill (.45,0) circle (1.5pt) (1.15,0) circle (1.5pt);
\end{tikzpicture}}
\newcommand{\trivAS}{\begin{tikzpicture} \useasboundingbox (-.4,-.3) rectangle (.4,.4);
\draw [->] (70:.3) arc (70:110:.3);
\draw (0,0) -- (0,-.4) (0,0) -- (50:.4) (0,0) -- (130:.4);
\fill (0,0) circle (1.5pt);
\end{tikzpicture}}
\newcommand{\trivASop}{\begin{tikzpicture} \useasboundingbox (-.4,-.3) rectangle (.4,.4);
\draw (0,0) -- (0,-.4);
\draw (0,0) .. controls (50:.3) ..  (130:.5);
\draw [draw=white,double=black,very thick] (0,0) .. controls (130:.3) ..  (50:.5);
\fill (0,0) circle (1.5pt);
\end{tikzpicture}}
\newcommand{\ihxone}{\begin{tikzpicture} \useasboundingbox (-.4,-.3) rectangle (.4,.4);
\draw [thick] (-90:.4) -- (0,0) -- (30:.4) (0,0) -- (150:.4);
\draw [thick] (-.3,-.25) -- (0,-.25);
\fill (0,-.25) circle (1.5pt) (0,0) circle (1.5pt);
\end{tikzpicture}}
\newcommand{\ihxtwo}{\begin{tikzpicture} \useasboundingbox (-.4,-.3) rectangle (.4,.4);
\draw [thick] (-90:.4) -- (0,0) -- (30:.4) (0,0) -- (150:.4);
\draw [out=0,in=-60, draw=white,double=black,very thick] (-.3,-.25) to (30:.25);
\fill (30:.25) circle (1.5pt) (0,0) circle (1.5pt);
\end{tikzpicture}}
\newcommand{\ihxthree}{\begin{tikzpicture} \useasboundingbox (-.4,-.3) rectangle (.4,.4);
\draw [thick]  (-90:.4) -- (0,0) -- (30:.4) (0,0) -- (150:.4);
\draw [draw=white,double=black,very thick] (-.3,-.25) .. controls (-.2,-.25) and  (.25,-.2)  .. (.25,0)  .. controls (.25,.2) and (120:.35) .. (150:.25);
\fill (150:.25) circle (1.5pt) (0,0) circle (1.5pt);
\end{tikzpicture}}
\newcommand{\thetabibell}{\begin{tikzpicture}\useasboundingbox (-.1,-.6) rectangle (2.1,0);
\begin{scope}[yshift=-.5cm]
\draw [-] (0,0) .. controls (0,.6) and (.75,.6) .. (1,.6);
\draw [->] (0,0) .. controls (0,-.6) and (.75,-.6) .. (1,-.6);
\draw [->] (0,0) -- (1,0);
\draw (1,0) -- (2,0);
\draw (1,.6) .. controls (1.25,.6) and (2,.6) .. (2,0) (1,-.6) .. controls (1.25,-.6) and (2,-.6) .. (2,0) (1,-.1) node[above]{\scriptsize $B(i,\sigma)$} (1,-.4) node{\scriptsize $B(\ell,\tau)$};
\fill (0,0) circle (1.5pt) (2,0) circle (1.5pt);\end{scope}
\end{tikzpicture}}
\newcommand{\ansothree}{\tilde{\alpha}}
\begin{document} 
\title{A universal equivariant finite type knot invariant defined from configuration space integrals}
\author{Christine Lescop \thanks{Institut Fourier, UJF Grenoble, CNRS}}
\maketitle
\begin{abstract} 
In a previous article, we constructed an invariant $\tilde{Z}$ for null-homologous knots in rational homology spheres, from equivariant intersections in configuration spaces. Here we present an equivalent definition of $\tilde{Z}$ in terms of configuration space integrals, we prove that $\tilde{Z}$ is multiplicative under connected sum,
and we prove null Lagrangian-preserving surgery formulae for $\tilde{Z}$. 
Our formulae generalize similar formulae that are satisfied by the Kricker rational lift of the Kontsevich integral for null Borromean surgeries. They imply that $\tilde{Z}$ is universal with respect to a natural filtration.
According to results of Garoufalidis and Rozansky, they therefore imply that $\tilde{Z}$ is equivalent to the Kricker lift of the Kontsevich integral for null-homologous knots with trivial Alexander polynomial in integral homology spheres.

\vskip.5cm

\end{abstract}

\tableofcontents

\maketitle
\section{Introduction}
In this article, a \emph{rational homology sphere} is a compact oriented $3$-manifold $R$ with the same homology as the usual $3$-sphere $S^3$. A \emph{null-homologous knot} in such an $R$ is a knot
whose homology class is trivial in $H_1(R;\ZZ)$.
In \cite{lesbonn}, we constructed a series of invariants for null-homologous knots in rational homology spheres, from equivariant intersections in configuration spaces. For a knot $K$, our series lives in an algebra of Feynman-Jacobi diagrams beaded by rational functions whose denominators divide the Alexander polynomial of $K$. Our construction is similar to configuration space integral constructions of invariants for knots or rational homology spheres performed by Altsch\"uler, Freidel, Axelrod, Singer, Bar-Natan, Bott, Cattaneo, Taubes, Guadagnini, Martellini, Mintchev, Poirier, Polyak and Viro in \cite{af,axelsingI,axelsingII,barnatanper,bottcat,botttaubes,cat,gmm,ko,kt,lesconst,poirier,pv} and others ...  that emerged after the Witten work \cite{witten}, in an equivariant setting as in \cite{marcheeq}.
Configuration space integrals associated with knots in $\RR^3$ of \cite{af,poirier} yield a series of invariants that is known to be equivalent to the Kontsevich integral. This series is a universal Vassiliev (or finite type) knot invariant.
Similarly, the series of configuration space invariants for integral homology spheres of \cite{kt,lesconst} is equivalent to the LMO invariant of Le, Murakami and Ohtsuki, and it is a universal finite type invariant of these homology spheres. According to recent work by Massuyeau and Moussard \cite{moussardAGT,massuyeausplit}, similar results hold for rational homology spheres.
In order to prove universality for such an invariant, one proves that the invariant induces an isomorphism from the graded space associated with a filtration of the vector space generated by knots or integral homology spheres, to the target space of the invariant.

In this article, we prove null Lagrangian-preserving surgery formulae for the series of equivariant invariants of \cite{lesbonn}. 
Our formulae generalize formulae that are satisfied by the Kricker rational lift of the Kontsevich integral \cite[Theorem 4]{garroz}. According to results of Garoufalidis and Rozansky in \cite{garroz}, they imply that our series $\tilde{Z}=(\tilde{Z}_n)_{n\in \NN}$ is universal for null-homologous knots with trivial Alexander polynomial in integral homology spheres, among the finite type invariants with respect to the Garoufalidis-Rozansky filtration. Therefore, $\tilde{Z}$ is equivalent to the Kricker rational lift of the Kontsevich integral for these knots, where two knot invariants are {\em equivalent\/} when two knots are distinguished by one of them if and only if they are distinghished by the other one.

Results of Delphine Moussard \cite{moussardthese} imply that $\tilde{Z}$ is also equivalent to this Kricker invariant for null-homologous knots with trivial Alexander polynomial in rational homology spheres, and I conjecture that $\tilde{Z}$ is also equivalent to the Kricker invariant for knots with any given equivariant linking form.

We also prove that $\tilde{Z}$ is multiplicative under connected sum.

\subsection{Conventions and notation}
Unless otherwise mentioned, manifolds are smooth, compact and oriented.
Boundaries are always oriented with the outward normal first convention. 
The fiber $N_xV$ of the normal bundle $N(V)$ to an oriented submanifold $V$ in an oriented manifold $M$ is oriented so that the tangent bundle $T_xM$ to the ambient manifold $M$ at some $x \in V$ is oriented as  $T_xM = N_xV \oplus T_xV$.
If $V$ and $W$ are two oriented transverse submanifolds of an oriented manifold $M$, their intersection is oriented so that the normal bundle to $T_x(V \cap W)$ at $x$ is the sum 
$N_xV \oplus N_xW$.
If the two manifolds are of complementary dimensions, then the sign of an intersection point is $+1$ if the orientation of its normal bundle coincides with the orientation of the ambient space (this is equivalent to: $T_xM=T_xV \oplus T_xW$ as oriented vector spaces). Otherwise, the sign is $-1$. If $V$ and $W$ are compact and if $V$ and $W$ are of complementary dimensions in $M$, their algebraic intersection is the sum of the signs of the intersection points, it is denoted by
$\langle V, W \rangle_M$.

Let $K$ be a null-homologous knot in a rational homology sphere $R$. A \emph{meridian} of $K$ is the boundary of a disk that intersects $K$ once transversely, with a positive sign. Let $\EXT$ be the \emph{exterior} of $K$, that is the complement of an open tubular neighborhood of $K$ in $R$, and let $p_{\EXT} \colon \tilde{\EXT} \rightarrow \EXT$ be the infinite cyclic covering of $\EXT$. The action of the meridian $m(K)$ of $K$ on $\tilde{\EXT}$ is denoted by $\theta_{\EXT}$. The induced action on $H_1(\tilde{\EXT};\QQ)$ is denoted as the multiplication by $t$ so that $H_1(\tilde{\EXT};\QQ)$ is a $\QQ[t,t^{-1}]$--module. Let $\delta(K)$ denote the annihilator of $H_1(\tilde{\EXT};\QQ)$. It is a divisor of the Alexander polynomial  $\Delta(K)$ of $K$ and it coincides with $\Delta(K)$ when $\Delta(K)$ has only simple roots. Both $\delta(K)$ and $\Delta(K)$ are normalized in $\QQ[t,t^{-1}]$ so that $\Delta(K)(1)=1=\delta(K)(1)$, $\Delta(K)(t)=\Delta(K)(t^{-1})$ and
$\delta(K)(t)=\delta(K)(t^{-1})$ or $\delta(K)(t)=t\delta(K)(t^{-1})$ (see \cite{moussardJKTR}).

Let $(J,L)$ be a two-component link of $\tilde{\EXT}$ such that
$p_{\EXT}(J) \cap p_{\EXT}(L)=\emptyset$. If $J$ bounds an oriented surface $\Sigma$
in $\tilde{\EXT}$ transverse to $\cup_{n \in \ZZ}\theta_{\EXT}^{n}(L)$, define the {\em equivariant intersection\/} $\langle \Sigma, L \rangle_e$ as
$$\langle \Sigma, L \rangle_e = \sum_{n \in \ZZ} t^n \langle \Sigma, \theta_{\EXT}^{n}(L) \rangle_{\tilde{\EXT}}$$
where $\langle \Sigma, \theta_{\EXT}^{n}(L) \rangle_{\tilde{\EXT}}$ is the usual $\QQ$--bilinear algebraic intersection.
Then the {\em equivariant linking number\/} of $J$ and $L$ is 
$$lk_e(J,L)= \langle \Sigma, L \rangle_e.$$
In general, $\delta(K)(\theta_{\EXT})(J)$ bounds a rational chain $\delta(K)(\theta_{\EXT})\Sigma$ and
$$lk_e(J,L)= \frac{\langle \delta(K)(\theta_{\EXT})\Sigma, L \rangle_e}{\delta(K)(t)}.$$
For any two one-variable polynomials $P$ and $Q$ in $\QQ[t,t^{-1}]$,
$$lk_e(P(\theta_{\EXT})J,Q(\theta_{\EXT})L)=P(t)Q(t^{-1})lk_e(J,L).$$

\subsection{The space $\CA^{h}_n(\delta)$ of beaded diagrams}

The invariants we constructed in \cite{lesbonn} for null-homologous knots in rational homology spheres live in the following spaces of beaded Jacobi diagrams.

A {\em (finite) trivalent graph\/} is said to be {\em oriented\/} if each of its vertices is equipped with a {\em vertex orientation\/}, that is a cyclic order of the three half-edges that meet at this vertex.
When an oriented graph is represented by the image of one of its planar immersions, the vertex orientation is induced by the counterclockwise order of the half-edges meeting at this vertex.
$$
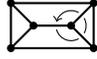
\begin{tikzpicture} \useasboundingbox (-.8,-.3) rectangle (.3,.4);
\draw (-165:.2) arc (-165:-60:.2) (-30:.2) arc (-30:30:.2);
\draw [->] (60:.2) arc (60:165:.2);
\draw [thick] (-.5,0) -- (-.8,.3) -- (-.8,-.3) -- (-.5,0) -- (0,0) --  (.3,-.3) -- (.3,.3) -- (0,0) (-.8,.3) -- (.3,.3) (-.8,-.3) -- (.3,-.3);
\fill (-.5,0) circle (1.5pt) (0,0) circle (1.5pt) (-.8,.3) circle (1.5pt) (.3,.3) circle (1.5pt) (-.8,-.3) circle (1.5pt) (.3,-.3) circle (1.5pt);
\end{tikzpicture}$$

When $\delta \in \QQ[t^{\pm 1}]$ and $\delta(1)=1$,
let $\CA^{h}_n(\delta)$ be the quotient of the rational vector space generated by oriented trivalent graphs with $2n$ vertices whose edges are oriented and equipped with some rational functions of $\frac1{\delta(t)}\QQ[t^{\pm 1}]$, by the following relations:
\begin{itemize}
\item \emph{Conjugation:} Reversing the orientation of an edge beaded by $P(t)$ and transforming this $P(t)$ into $P(t^{-1})$
gives the same element in the quotient.
$$\begin{tikzpicture} \useasboundingbox (-1.2,-.5) rectangle (.4,.2);
\begin{scope}[yshift=-.3cm]
\draw [thick] (-.25,-.55) -- (0,-.4) -- (0,.4) -- (.25,.55) (.25,-.55) -- (0,-.4) (0,.4) -- (-.25,.55);
\draw [thick,->] (0,-.4) -- (0,.05);
\draw (0,-.1) node[left]{\footnotesize $P(t)$};
\fill (0,-.4) circle (1.5pt) (0,.4) circle (1.5pt);
\end{scope}
\end{tikzpicture} = \begin{tikzpicture} \useasboundingbox (-.4,-.5) rectangle (1.2,.5);
\begin{scope}[yshift=-.3cm]
\draw [thick] (-.25,-.55) -- (0,-.4) -- (0,.4) -- (.25,.55) (.25,-.55) -- (0,-.4) (0,.4) -- (-.25,.55);
\draw [thick,->] (0,.4) -- (0,-.05);
\draw (0,-.1) node[right]{\footnotesize $P(t^{-1})$};
\fill (0,-.4) circle (1.5pt) (0,.4) circle (1.5pt);
\end{scope}
\end{tikzpicture}$$

\item \emph{Multilinearity:} If two graphs only differ by the label of one oriented edge, that is $P(t)$
for one of them and $Q(t)$ for the other one,
then the class of their sum is the class of the same graph with label $(P(t)+Q(t))$.

$$\begin{tikzpicture} \useasboundingbox (-2.6,-.5) rectangle (.4,.2);
\begin{scope}[yshift=-.3cm]
\draw [thick] (-.25,-.55) -- (0,-.4) -- (0,.4) -- (.25,.55) (.25,-.55) -- (0,-.4) (0,.4) -- (-.25,.55);
\draw [thick,->] (0,-.4) -- (0,.05);
\draw (0,-.1) node[left]{\footnotesize $P(t)+Q(t)$};
\fill (0,-.4) circle (1.5pt) (0,.4) circle (1.5pt);
\end{scope}
\end{tikzpicture} = \begin{tikzpicture} \useasboundingbox (-1,-.5) rectangle (.4,.2);
\begin{scope}[yshift=-.3cm]
\draw [thick] (-.25,-.55) -- (0,-.4) -- (0,.4) -- (.25,.55) (.25,-.55) -- (0,-.4) (0,.4) -- (-.25,.55);
\draw [thick,->] (0,-.4) -- (0,.05);
\draw (0,-.1) node[left]{\footnotesize $P(t)$};
\fill (0,-.4) circle (1.5pt) (0,.4) circle (1.5pt);
\end{scope}
\end{tikzpicture} + \begin{tikzpicture} \useasboundingbox (-1,-.5) rectangle (.4,.2);
\begin{scope}[yshift=-.3cm]
\draw [thick] (-.25,-.55) -- (0,-.4) -- (0,.4) -- (.25,.55) (.25,-.55) -- (0,-.4) (0,.4) -- (-.25,.55);
\draw [thick,->] (0,-.4) -- (0,.05);
\draw (0,-.1) node[left]{\footnotesize $Q(t)$};
\fill (0,-.4) circle (1.5pt) (0,.4) circle (1.5pt);
\end{scope}
\end{tikzpicture}$$

\item \emph{Holonomy:} Multiplying by $t$ the three rational functions of three edges adjacent to a vertex, oriented towards that vertex, does not change the element in the quotient.

$$
\begin{tikzpicture} \useasboundingbox (-1.5,-.6) rectangle (1.5,.3);
\begin{scope}[yshift=-.2cm]
\draw [thick] (0,-.7) -- (0,0) -- (30:.7) (150:.7) -- (0,0);
\draw [thick,-<] (0,0) -- (0,-.4);
\draw [thick,-<] (0,0) -- (30:.4);
\draw [thick,-<] (0,0) -- (150:.4);
\draw (0,-.4) node[left]{\footnotesize $P(t)$};
\draw (150:.3) node[above]{\footnotesize $R(t)$};
\draw (.25,0) node[right]{\footnotesize $Q(t)$};
\fill (0,0) circle (1.5pt);
\end{scope}
\end{tikzpicture}
= \begin{tikzpicture} \useasboundingbox (-1.5,-.6) rectangle (1.5,.3);
\begin{scope}[yshift=-.2cm]
\draw [thick] (0,-.7) -- (0,0) -- (30:.7) (150:.7) -- (0,0);
\draw [thick,-<] (0,0) -- (0,-.4);
\draw [thick,-<] (0,0) -- (30:.4);
\draw [thick,-<] (0,0) -- (150:.4);
\draw (0,-.4) node[left]{\footnotesize $tP(t)$};
\draw (150:.3) node[above]{\footnotesize $tR(t)$};
\draw (.25,0) node[right]{\footnotesize $tQ(t)$};
\fill (0,0) circle (1.5pt);
\end{scope}
\end{tikzpicture} \;\;\;\mbox{and}\;\;\; 
\begin{tikzpicture} \useasboundingbox (-1,-.6) rectangle (1.2,.3);
\begin{scope}[yshift=-.2cm]
\draw [thick] (0,-.7) -- (0,0);
\draw [thick,-<] (0,0) -- (0,-.4);
\draw [thick] (.25,.255) arc (0:270:.25);
\draw [thick,->] (0,0) arc (-90:0:.25);
\draw (0,-.4) node[left]{\footnotesize $P(t)$};
\draw (.2,.25) node[right]{\footnotesize $Q(t)$};
\fill (0,0) circle (1.5pt);
\end{scope}
\end{tikzpicture}
=
\begin{tikzpicture} \useasboundingbox (-1,-.6) rectangle (1.2,.3);
\begin{scope}[yshift=-.2cm]
\draw [thick] (0,-.7) -- (0,0);
\draw [thick,-<] (0,0) -- (0,-.4);
\draw [thick] (.25,.255) arc (0:270:.25);
\draw [thick,->] (0,0) arc (-90:0:.25);
\draw (0,-.4) node[left]{\footnotesize $tP(t)$};
\draw (.2,.25) node[right]{\footnotesize $Q(t)$};
\fill (0,0) circle (1.5pt);
\end{scope}
\end{tikzpicture}
$$

\item \emph{Antisymmetry or AS:} Changing the orientation of a vertex multiplies the element of the quotient by $(-1)$.

$$\trivAS + \trivASop =0$$

\item \emph{Jacobi or IHX:} The sum of three graphs that coincide outside a disk, where they look as in the picture below, vanishes in the quotient.
(The complete edges of the relation are equipped with the polynomial $1$ that is not written.)

$$\ihxone + \ihxtwo + \ihxthree =0
\;\;\;\;\;\;\;\;\mbox{where} \begin{tikzpicture} \useasboundingbox (-1.2,-.5) rectangle (.4,.2);
\begin{scope}[yshift=-.3cm]
\draw [thick] (-.25,-.55) -- (0,-.4) -- (0,.4) -- (.25,.55) (.25,-.55) -- (0,-.4) (0,.4) -- (-.25,.55);
\draw [thick,->] (0,-.4) -- (0,.05);
\draw (0,-.1) node[left]{\footnotesize $1$};
\fill (0,-.4) circle (1.5pt) (0,.4) circle (1.5pt);
\end{scope}
\end{tikzpicture} =\begin{tikzpicture} \useasboundingbox (-1.2,-.5) rectangle (.4,.2);
\begin{scope}[yshift=-.3cm]
\draw [thick] (-.25,-.55) -- (0,-.4) -- (0,.4) -- (.25,.55) (.25,-.55) -- (0,-.4) (0,.4) -- (-.25,.55);
\fill (0,-.4) circle (1.5pt) (0,.4) circle (1.5pt);
\end{scope}
\end{tikzpicture}.
$$
\end{itemize}

A more general version of this space of diagrams was introduced in
\cite[Definition 3.8]{garkri}.
Let $\CA_n$ denote the quotient of the rational vector space generated by oriented trivalent graphs with $2n$ vertices whose edges are neither oriented nor beaded, by the relations AS and IHX. There are two natural well-defined maps $i \colon \CA_n \rightarrow \CA_n^h(\delta)$ that sends a diagram to itself viewed as a diagram whose edges are beaded by $1$, and  $p \colon \CA_n^h(\delta) \rightarrow \CA_n$ that sets $t=1$ in the beads that are next seen as coefficients. Since $p \circ i$ is the identity, $i$ injects $\CA_n$ into $\CA_n^h(\delta)$ for any $\delta$.

The topological vector spaces $\CA=\prod_{n\in \NN}\CA_n$ and $\CA^{h}(\delta)=\prod_{n\in \NN}\CA^{h}_n(\delta)$ are graded algebras, where the graded product of diagrams is their disjoint union, and the degree $n$ is called the \emph{loop-degree}.

There is a hair map from $\CA_n^h(\delta)$ (where diagrams with circle components without trivalent vertices are allowed) to an algebra $\CB$ of diagrams with legs, that is isomorphic to the target space of the original Kontsevich integral of knots, by a Poincar\'e-Birkhoff-Witt symmetrization isomorphism. 
This hair map multilinearly replaces an oriented edge labeled by $$P(t=\exp(x))=\sum_{n\in \NN}p_nx^n$$ by a combination of edges with hair (or legs), that is with new edges with univalent vertices as in the following figure.

$$\begin{tikzpicture} \useasboundingbox (-1.2,-.5) rectangle (.4,.2);
\begin{scope}[yshift=-.3cm]
\draw [thick] (-.25,-.55) -- (0,-.4) -- (0,.4) -- (.25,.55) (.25,-.55) -- (0,-.4) (0,.4) -- (-.25,.55);
\draw [thick,->] (0,-.4) -- (0,.05);
\draw (0,-.1) node[left]{\footnotesize $P(t)$};
\fill (0,-.4) circle (1.5pt) (0,.4) circle (1.5pt);
\end{scope}
\end{tikzpicture} \mapsto \sum_{n=0}^{\infty} p_n \begin{tikzpicture} \useasboundingbox (-.4,-.5) rectangle (6,.2);
\begin{scope}[yshift=-.3cm]
\draw [thick, dotted] (0,-.2) -- (0,0);
\draw [thick] (-.25,-.55) -- (0,-.4) -- (0,-.2) (0,0) -- (0,.4) -- (.25,.55) (.25,-.55) -- (0,-.4) (0,.4) -- (-.25,.55) (0,.25) -- (.2,.25) (0,-.25) -- (.2,-.25) (0,.1) -- (.2,.1);
\fill (0,-.4) circle (1.5pt) (0,.4) circle (1.5pt) (0,.25) circle (1.5pt) (0,.1) circle (1.5pt) (0,-.25) circle (1.5pt) (.2,.25) circle (1.5pt) (.2,.1) circle (1.5pt) (.2,-.25) circle (1.5pt);
\draw (.2,-.1) node[right]{\footnotesize $n$ horizontal edges called \emph{legs}};
\end{scope}\end{tikzpicture}$$

Patureau-Mirand proved that this hair map is non-injective \cite{patureauhair}, using the work of Vogel on Jacobi diagrams \cite{vogel}.
The Kricker lift of the Kontsevich integral of a knot $K$ is valued in a space similar to $\CA^{h}(\Delta(K))$ and its composition by a modified hair map (and the Poincar\'e-Birkhoff-Witt isomorphism) is the Kontsevich integral. See \cite[Theorem 14 and Section 7.1]{garkri}.

Our equivariant invariant $\tilde{Z}(K)$ is valued in $\CA^{h}(\delta(K))$, its definition from equivariant intersections in configuration spaces is given in Section~\ref{secrecdef}, and an alternative equivalent definition based on integrals on equivariant configuration spaces is presented in Section~\ref{secaltdefzform}. I expect that the universal finite type invariant of \cite{botttaubes,af,poirier} for knots in $\RR^3$ is obtained from $\tilde{Z}(K)$ in the same way as the Kontsevich integral is obtained from its Kricker rational lift.

\subsection{Statement of the main theorem}
\setcounter{equation}{0}

Let $K$ be a knot in a rational homology sphere $R$ and let $\EXT$ be its exterior, this exterior is implicitly equipped with the meridian of $K$ so that it determines $K$.

The main result of this article describes part of the behaviour of the $\tilde{Z}_n$ under the null Lagrangian-preserving surgeries
defined below.

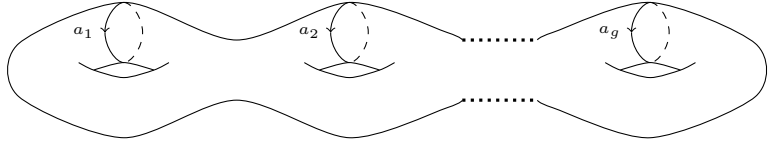
\begin{figure}[ht!]
\begin{center}
 \begin{tikzpicture} \useasboundingbox (0,-1) rectangle (11,1); 
\draw [->] (2,.9) .. controls (1.9,.9) and (1.75,.7) .. (1.75,.5) node[left]{\tiny $a_1$};  
\draw (1.75,.5) .. controls (1.75,.3) and (1.9,.1) .. (2,.1);
\draw [dashed] (2,.9) .. controls (2.1,.9) and (2.25,.7) .. (2.25,.5) .. controls (2.25,.3) and (2.1,.1) .. (2,.1);
\draw plot[smooth] coordinates{(1.4,.1) (1.6,0) (2,-.1) (2.4,0) (2.6,.1)};
\draw plot[smooth] coordinates{(1.6,0) (2,.1) (2.4,0)}; 
\draw plot[smooth] coordinates{(6.5,.4) (6.3,.5) (5,.9) (3.5,.4) (2,.9) (.6,.35) (.6,-.35) (2,-.9) (3.5,-.4) (5,-.9) (6.3,-.5) (6.5,-.4)};
\draw [very thick, dotted] (6.5,-.4) -- (7.5,-.4)  (6.5,.4) -- (7.5,.4);
\draw plot[smooth] coordinates{(7.5,.4) (7.7,.5) (9,.9) (10.4,.35) (10.4,-.35) (9,-.9) (7.7,-.5) (7.5,-.4)};
\begin{scope}[xshift=3cm]
 \draw [->] (2,.9) .. controls (1.9,.9) and (1.75,.7) .. (1.75,.5) node[left]{\tiny $a_2$};  
\draw (1.75,.5) .. controls (1.75,.3) and (1.9,.1) .. (2,.1);
\draw [dashed] (2,.9) .. controls (2.1,.9) and (2.25,.7) .. (2.25,.5) .. controls (2.25,.3) and (2.1,.1) .. (2,.1);
\draw plot[smooth] coordinates{(1.4,.1) (1.6,0) (2,-.1) (2.4,0) (2.6,.1)};
\draw plot[smooth] coordinates{(1.6,0) (2,.1) (2.4,0)};
\end{scope}
\begin{scope}[xshift=7cm]
 \draw [->] (2,.9) .. controls (1.9,.9) and (1.75,.7) .. (1.75,.5) node[left]{\tiny $a_g$};  
\draw (1.75,.5) .. controls (1.75,.3) and (1.9,.1) .. (2,.1);
\draw [dashed] (2,.9) .. controls (2.1,.9) and (2.25,.7) .. (2.25,.5) .. controls (2.25,.3) and (2.1,.1) .. (2,.1);
\draw plot[smooth] coordinates{(1.4,.1) (1.6,0) (2,-.1) (2.4,0) (2.6,.1)};
\draw plot[smooth] coordinates{(1.6,0) (2,.1) (2.4,0)};
\end{scope}
\end{tikzpicture}
\caption{The standard handlebody
$H_g$}\label{fig1}
\end{center}
\end{figure}

A {\em genus $g$ $\QQ$--handlebody}\/ is an (oriented, compact) 3--manifold $A$
with 
the same homology with rational coefficients as the standard (solid) handlebody
$H_g$ of Figure~\ref{fig1}. 
Note that the boundary of such a $\QQ$--handlebody $A$ is homeomorphic to the 
boundary $(\partial H_g =\Sigma_g)$ of $H_g$.

For a (compact, oriented) $3$-manifold $A$,
${\cal L}_A$ denotes the kernel of the map induced by the inclusion
$$ H_1(\partial A;\QQ) \longrightarrow H_1( A;\QQ).$$ 
It is a Lagrangian of $(H_1(\partial A;\QQ),\langle,\rangle_{\partial A})$, we
call it the {\em {Lagrangian}}\/ of $A$.
A $\QQ$--handlebody $A$ sitting in a 3--manifold $\EXT$ is said to be {\em rationally null-homologous\/} if the map induced by the inclusion maps $H_1(A;\QQ)$ to zero in $H_1(\EXT;\QQ)$.
A {\em null Lagrangian-preserving surgery\/} or {\em null LP--surgery} $(A,A^{\prime})$ is the replacement of a rationally null-homologous
$\QQ$--handlebody $A$ embedded in a 3--manifold $\EXT$ by
another such $A^{\prime}$ with identical (identified via a homeomorphism) boundary and Lagrangian.

There is a canonical isomorphism 
$$\partial_{MV} \colon H_2(A \cup_{\partial A}
-A^{\prime};\QQ) \rightarrow {\cal L}_A$$
that maps the class of a closed surface in the closed 3--manifold $(A \cup_{\partial A}
-A^{\prime})$ to the boundary of its intersection with $A$.
This isomorphism carries the algebraic triple intersection of surfaces to a trilinear antisymmetric form $\CI_{AA^{\prime}}$ on $\CL_A$.
$$\CI_{AA^{\prime}}(a_{i},a_{j},a_{k})=\langle \partial_{MV}^{-1}(a_i), \partial_{MV}^{-1}(a_j), \partial_{MV}^{-1}(a_k)\rangle_{A \cup
-A^{\prime}}$$

Let $(a_1,a_2,\dots,a_g)$ be a basis of $\CL_A$, and let $z_1,\dots,z_g$ be homology classes of $\partial A$, such that $(z_1,z_2,\dots,z_g)$ is dual to $(a_1,a_2,\dots,a_g)$ with respect to $\langle,\rangle_{\partial A}$ 
($\langle a_i,z_j \rangle_{\partial A}=\delta_{ij}=\left\{\begin{array}{ll} 1 & \mbox{if}\; i=j\\ 0 & \mbox{if}\; i\neq j\end{array}\right.$). Note that $(z_1,\dots,z_g)$
is a basis of $H_1(A;\QQ)$.

Let $A$ be a rationally null-homologous $\QQ$-handlebody in $\EXT$. Then $A$ lifts as $\ZZ$ homeomorphic copies in $\tilde{\EXT}$. Fix one of these lifts $\tilde{A}$, and for any $i$, let $\tilde{z}_i$ denote the lift of $z_i$ in $\tilde{A}$.

Represent $\CI_{AA^{\prime}}$
by the following combination $\tilde{T}(\CI_{AA^{\prime}})$ of tripods 
whose three univalent vertices form an ordered set:

$$\tilde{T}(\CI_{AA^{\prime}})=\sum_{\left\{\{i,j,k\}  \subset \{1,2,\dots ,g\} ; i <j <k\right\}}\CI_{AA^{\prime}}(a_{i},a_{j},a_{k})
\begin{tikzpicture} \useasboundingbox (-.2,-.4) rectangle (1.5,0);
\begin{scope}[yshift=-.3cm]
\draw [thick] (.4,0) -- (0,0) -- (.4,.3) (0,0) -- (.4,-.3);
\draw (.4,.3) node[right]{\footnotesize $\tilde{z}_k$} (.4,0) node[right]{\footnotesize $\tilde{z}_j$} (.4,-.3) node[right]{\footnotesize $\tilde{z}_i$};
\fill (0,0) circle (1.5pt);
\end{scope}
\end{tikzpicture}$$

Let $G$ be a graph with $2n$ oriented trivalent vertices and with univalent vertices. Assume that the univalent vertices of $G$ are decorated with lifts in $\tilde{\EXT}$ of disjoint rationally null-homologous
curves of $\EXT$.
Let $P(G)$ be the set of partitions of the set of univalent vertices of $G$ in disjoint pairs.

For $p \in P(G)$, identifying the two vertices of each pair provides a vertex-oriented trivalent diagram. Orient each edge of this diagram arbitrarily and bead it with the equivariant linking number of the two curves associated with its two half-edges, ordered by the edge orientation (\begin{tikzpicture} \useasboundingbox (-.5,-.1) rectangle (.5,.2);
\draw [thick] (.4,-.1) -- (.3,0) -- (.4,.1) (-.4,-.1) -- (-.3,0) -- (-.4,.1);
\draw [thick,->](.3,0) -- (0,0)  (-.3,0) -- (0,0);
\draw (.15,.15) node{\tiny $2$} (-.15,.15) node{\tiny $1$};
\fill (.3,0) circle (1.5pt) (-.3,0) circle (1.5pt);
\end{tikzpicture}). This yields $[G_p] \in \CA^{h}_n(\delta(K))$.

Define the contraction $\langle \langle G \rangle \rangle_n$ of $G$ 

$$\langle \langle G \rangle \rangle_n=\sum_{p \in P(G)} [G_p].$$ 

The contraction $\langle \langle . \rangle \rangle$ is linearly extended to linear combination
of graphs, and the disjoint union of combinations of graphs is bilinear.

A {\em $k$--component null Lagrangian-preserving surgery datum\/} in a knot complement $\EXT$ is a datum $(\EXT;(A^{(i)},A^{(i)\prime})_{i \in \{1,\dots,k\}})$
of $k$ disjoint rationally null-homologous $\QQ$--handlebodies $A_i$, for $i \in \{1,\dots,k\}$, in $\EXT$, and $k$ associated LP--surgeries $(A^{(i)},A^{(i)\prime})$.

The main theorem of this preprint is the following one.

\begin{theorem}
\label{thmflag}
Set $N=\{1,\dots,2n\}$.
Let  $$(\EXT;(A^{(i)},A^{(i)\prime})_{i \in N})$$
be a \/ $2n$--component null Lagrangian-preserving surgery datum in the exterior $\EXT$ of a null-homologous knot $K$ in a rational homology sphere. For $S \subset N$, let $\EXT_S$ denote the manifold obtained from $\EXT$ by replacing $A^{(i)}$ by $A^{(i)\prime}$ for all $i\in S$, and let $K_S$ denote the knot with exterior $\EXT_S$ and with the same meridian as $K$. Then
$$\sum_{S \subset N} (-1)^{\sharp S}\tilde{Z}_n(K_S)=\langle \langle \bigsqcup_{i \in N} \tilde{T}(\CI_{A^{(i)}A^{(i)\prime}}) \rangle \rangle_n.$$
\end{theorem}
The sum of the left-hand side runs over all the subsets of $N$ (including $\emptyset$ and $N$).
Note that the right-hand side of the above equality is independent of the chosen lifts of the 
 $A^{(i)}$ thanks to the holonomy relation and to the conjugation relation in $\CA^h(\delta(K))$.

\subsection{Organization of the article}

In order to prove Theorem~\ref{thmflag}, we use the definition of $\tilde{Z}$ from equivariant intersections in configuration spaces that is given in Section~\ref{secrecdef}. The definition will start with a construction of invariants $\tilde{Z}(\EXT,\tau)$ of knot exteriors $\EXT$ equipped with some parallelization $\tau$. Then Pontrjagin classes allow us to associate a numerical invariant $p_1(\tau_{\pi},\tau)$ with such a parallelization so that 
$$\tilde{Z}(\EXT)= \tilde{Z}(\EXT,\tau)\exp\left(-\frac{p_1(\tau_{\pi},\tau)}{4} \ansothree\right)$$ 
for a constant  
$$\ansothree=\sum_{n \in \NN} \ansothree_{2n+1} \in \left(\CA=\prod_{k\in \NN}\CA_k\right)$$ 
that is called the {\em anomaly.\/}

In order to prove Theorem~\ref{thmflag} in general, we need to make our construction more flexible and use generalisations $\tilde{\tau}$ of parallelizations that we call pseudo-parallelizations, that are also equipped with a numerical invariant $p_1(\tau_{\pi},\tilde{\tau})$
and we define  $\tilde{Z}(\EXT,\tilde{\tau})$ so that
$$\tilde{Z}(\EXT)= \tilde{Z}(\EXT,\tilde{\tau})\exp\left(-\frac{p_1(\tau_{\pi},\tilde{\tau})}{4} \ansothree\right)$$

These pseudo-parallelizations defined in Section~\ref{secpseudotriv} are not needed for the proof of Theorem~\ref{thmflag} when the $A^{(i)}$ and the $A^{(i)\prime}$ are integral homology handlebodies, as in the case of Borromean surgeries.
Sections~\ref{secprelproof} and \ref{secproof} give the proof of Theorem~\ref{thmflag} assuming natural extensions of parallelization properties to our pseudo-parallelizations.
Then the last sections are devoted to the justification of these properties that relies on earlier work in the setting of configuration space integrals instead of the dual equivariant intersections in configuration spaces.
Section~\ref{secaltdefzform} contains a definition of $\tilde{Z}$ based on configuration space integrals, and is interesting in its own.
Sections~\ref{secdefpseudointform} and \ref{secdefpseudoint} contain the more flexible definitions of $\tilde{Z}$ with respect to pseudo-parallelizations rather than parallelizations, in the settings of differential forms and equivariant intersections, respectively.
Section~\ref{secrelPont} contains the precise definition of the invariant $p_1(\tau_{\pi},.)$ and some of its properties.
The anomaly is completely defined in Section~\ref{secanomaly}.

Finally, we study the behaviour of $\tilde{Z}$ under connected sums in Section~\ref{secconnsum}.

\begin{remark}
This preprint won't be submitted in its present form.
Sections \ref{secrelPont} to \ref{secdefpseudoint} can be thought of as appendices. Their results will be included in a monograph about \emph{Invariants of links and $3$--manifolds from graph configurations} that I am compiling.
I added Section~\ref{secconnsum} because its results are needed in Delphine Moussard's Ph.D. thesis \cite{moussardthese}. 
\end{remark}

\section{Definition of the invariant $\tilde{Z}$}
\label{secrecdef}
\setcounter{equation}{0}

\subsection{Introduction}

Let $K$ be a knot in a rational homology sphere $R$ and let $\EXT$ be its exterior, this exterior is equipped with the meridian $m(K)$ of $K$ so that it determines $K$.
Let $p_{\EXT} \colon \tilde{\EXT} \rightarrow \EXT$ denote the infinite cyclic covering of $\EXT$, and let $\delta=\delta(K)$ be the annihilator of $H_1(\tilde{\EXT};\QQ)$.

In this subsection, we define the invariant $\tilde{Z}(K)=\tilde{Z}(\EXT)$, following 
\cite{lesbonn}, except for the differences listed below.

\begin{itemize}
\item Let $\Mbet=R(K,0)$ be the manifold obtained from $R$ by $0$-surgery on $K$, let $\KK$ be the core of the surgery framed by the meridian of $K$, and let $\overline{\tilde{z}_n(\Mbet,\KK)}$ be obtained from the invariant $\tilde{z}_n(\Mbet,\KK) \in \CA^h_n(\delta(K))$ in
\cite{lesbonn} by {\em conjugating\/} all the beads, that is by replacing $t$ by $t^{-1}$.
Here, we call $\tilde{z}_n(K)$ or $\tilde{z}_n(\EXT)$ the configuration space invariant $\overline{\tilde{z}_n(\Mbet,\KK)}$. It lives in the subspace
of $\CA^{h}_n(\delta(K))$ made of connected beaded diagrams. This convention is more natural and leads to nicer formulae.
\item Here, we consider non-connected trivalent graphs in order to define $$\tilde{Z}(K)=(\tilde{Z}_n(K))_{n\in \NN}=\exp\left(\sum_{n\in \NN}\tilde{z}_n(K)\right) \in\CA^{h}(\delta(K)).$$
\item We do not repeat an error in \cite{lesbonn}, where $\tilde{M}_{[1,k]}$ should be replaced by $\overline{\tilde{M}_{]1,k]}}$ and ${M}_{[1,k]}$ should be replaced by 
$\overline{{M}_{]1,k]}}$ in Section~1.8.
\item The map $\funcE$ was denoted by $f_M$ in \cite{lesbonn}.
\item Parallelization maps $\tau$ go from $M \times \RR^3$ to $TM$, here, while they go from $TM \times \RR^3$ to $M\times \RR^3$ in \cite{lesconst,lessumgen,lesbonn}. The anomaly $\xi$ is replaced by the anomaly $\ansothree=-\xi$. This convention is also more natural.
\item The proofs of \cite{lesbonn} are not repeated here, but some notions are presented with more details here.
\end{itemize}

The invariant $\tilde{Z}(K)$ will be defined from equivariant intersections in configuration spaces. 
More precisely, we shall define a $\ZZ$-covering $\tilde{C}_2(\EXT)$ of a compactification
${C}_2(\EXT)$ of $\EXT^2 \setminus \mbox{diagonal}$. This compactification ${C}_2(\EXT)$ will be a $6$-manifold with boundary and ridges that has the same homotopy type as $\EXT^2 \setminus \mbox{diagonal}$.
Then we shall define the equivariant intersection $I_n(\{H_i\}_{i=1,\dots,3n})$ of $3n$ generic cycles $H_i$ of $(\tilde{C}_2(\EXT),\partial \tilde{C}_2(\EXT))$ as an element of $\CA^h_n(\delta(K))$.
Next, we shall fix a parallelization $\tau$ of $\EXT$, and fix cycles $\deleqprop_i(\tau)$ of $(\tilde{C}_2(\EXT),\partial \tilde{C}_2(\EXT))$ on a neighborhood of $\partial \tilde{C}_2(\EXT)$ so that
$$\tilde{Z}_n(K,\tau)=I_n(\{\frac1{\delta(t)}\deleqprop_i(\tau)\}).$$
The fundamental $\eqprop_i(\tau)=\frac1{\delta(t)}\deleqprop_i(\tau)$ will be called \emph{equivariant propagators}.
Finally, $\tilde{Z}(K)$ will be defined as $\tilde{Z}(K)= \tilde{Z}(K,\tau)\exp\left(-\frac{p_1(\tau_{\pi},\tau)}{4} \ansothree \right)$ for the \emph{anomaly} $$\ansothree =-\xi=\sum_n \ansothree_{2n+1}\in \CA,$$ that is a constant described in Definition~\ref{defxin} or in \cite[Subsection 1.6]{lesconst}, and for the Pontrjagin number $p_1(\tau_{\pi},\tau)$ associated with $\tau$ described in Subsection~\ref{subpont}.
We first describe the involved configurations spaces.

\subsection{Involved configuration spaces}
The {\em configuration space\/} $C_2(\Mbet)$ is the configuration space obtained from $\Mbet^2$ by {\em blowing up\/} the diagonal, that is by replacing it by its unit normal bundle (in this article). (In a trivialized tubular neighborhood $\RR^3 \times \mbox{diag}(M^2)$ of the diagonal, the fiber $\RR^3= \{0\} \cup (\RR^{+\ast}\times S^2)$ is replaced by $\RR^{+}\times S^2$, so that $(0,(m,m))$ is replaced by the fiber $\{0\}\times S^2 \times \{(m,m)\}$ of the unit normal bundle to the diagonal at $(m,m)$.) In particular, the boundary $\partial C_2(\Mbet)$ of $C_2(\Mbet)$ is the unit normal bundle to the diagonal that is identified with the unit tangent bundle $U\Mbet$ to $\Mbet$ by the map induced by the following one
$$\begin{array}{lll} T_{(m,m)}(\Mbet^2)=(T_m\Mbet)^2 & \rightarrow& T_m \Mbet\\
   (x,y) &\mapsto & y-x.
  \end{array}$$

The points of the interior of $C_2(\Mbet)$ are pairs $(m_1,m_2)$ of distinct points 
of $\Mbet$, and the points of $\partial C_2(\Mbet)$ are pairs $(m,m)$ equipped with a unit tangent vector to $\Mbet$ at $m$.

Consider the quotient $\tilde{C}_2(\Mbet)$ of the space of paths $\gamma$ from a point $\gamma(0)$ of $\Mbet$ to another such $\gamma(1)$ - equipped with a tangent vector of $T_{\gamma(0)} \Mbet$, when $\gamma(1)=\gamma(0)$ - by the relation that identifies two paths $\gamma$ and $\gamma^{\prime}$ if and only if $\gamma(0)=\gamma^{\prime}(0)$, $\gamma(1)=\gamma^{\prime}(1)$ and the path composition $\gamma\overline\gamma^{\prime}$ vanishes in $H_1(\Mbet;\QQ)$.
Equip $\tilde{C}_2(\Mbet)$ with the natural topology that makes 
$$\begin{array}{llll}p_C \colon &\tilde{C}_2(\Mbet)& \rightarrow& {C}_2(\Mbet)\\
   & \gamma &\mapsto & (\gamma(0),\gamma(1))
  \end{array}$$
become an infinite cyclic covering.
The transformation $\thetap$ that maps the class of a path $\gamma$ to its path composition with a loop
conjugated to the oriented meridian of $K$ generates the transformation group of the covering.

Note that the boundary $U\Mbet$ of $C_2(\Mbet)$ has a canonical lift $\{0\} \times U\Mbet$ that corresponds to constant paths $\gamma$ in the description above.
This boundary $U\Mbet$ lifts as $\ZZ \times U\Mbet=\coprod_{n\in \ZZ}\thetap^n(\{0\} \times U\Mbet)$ in $\tilde{C}_2(\Mbet)$.

When $A$ is a $3$-dimensional submanifold of $\Mbet$, $C_2(A)$ denotes the preimage of $A^2$ under the blow-down map
from $C_2(\Mbet)$ to $\Mbet^2$, and $\tilde{C}_2(A)$ is the preimage of $C_2(A)$ under the covering map $p_C$ from $\TCM$ to $C_2(\Mbet)$. This applies to $A=\EXT$.

The interior of $\tilde{C}_2(\EXT)$ is the quotient of $\tilde{\EXT}^2\setminus (p_{\EXT} \times p_{\EXT})^{-1}(\mbox{diag}(\EXT^2))$ by the equivalence relation such that the equivalence class of $(\tilde{m}_1, \tilde{m}_2)$ is $\{(\theta_{\EXT}^n(\tilde{m}_1),\theta_{\EXT}^n(\tilde{m}_2));n\in \ZZ\}$, where $\theta_{\EXT}$ is the covering transformation of $\tilde{\EXT}$ induced by the action of the oriented meridian of $K$.

The covering transformation $\thetap$ of $\TCM$ maps the class of $(\tilde{m}_1, \tilde{m}_2)$ 
to the class of $(\tilde{m}_1, \theta_{\EXT}(\tilde{m}_2))$.

We also consider the configuration spaces $C_{2n}(\Mbet)$ that are obtained from $\Mbet^{2n}$ by blowing up all the partial diagonals of $\Mbet^{2n}$, starting with the diagonals where the number of coinciding points is maximal, successively. This configuration space is denoted by $C_{2n}[\Mbet]$ in \cite{Sinha} where it is presented in details. See also \cite[Section 5]{axelsingII} and \cite[Section 3]{lesconst}.
Again, for a submanifold $A$ of $\Mbet$, $C_{2n}(A)$ will denote the preimage of $A^{2n}$ under the blow-down map $p_b$
from $C_{2n}(\Mbet)$ to $\Mbet^{2n}$.

\subsection{Equivariant intersections valued in spaces of beaded graphs}
\label{subeqint}

Let $\CS^u_n$ (resp. $\CS_n$) be the set of trivalent graphs (resp. connected trivalent graphs) $\Gamma$ with $2n$ vertices numbered from $1$ to $2n$ and with $3n$ oriented edges numbered from $1$ to $3n$ without loops (edges with coinciding ends).
Formally, such a $\Gamma$ is an ordered set of $3n$ ordered pairs of elements of $\{1,2,\dots,2n\}$ such that any element of $\{1,2,\dots,2n\}$ belongs to three pairs.
The \emph{degree} of a graph of $\CS^u_n$ is $n$.

Any $\Gamma \in \CS^u_n$ has the following {\em canonical vertex-orientation\/} $o(\Gamma)$ (up to an even number of changes).
Orient the vertices of $\Gamma$ so that the permutation of the half edges from (first half of first edge, second half of first edge, \dots, second half of last edge) to (half-edges of the first vertex $v(1)$
ordered in a way compatible with the vertex orientation, half-edges of $v(2)$
ordered in a way compatible with the vertex orientation, \dots) is even.

Consider the {\em generalized holonomy relation\/} generated by the identification of a beaded graph with the one obtained by multiplying the beads of incoming edges by $t$ and the beads of the outcoming edges by $t^{-1}$ at some vertex. This relation is usually a consequence of the holonomy relation and the conjugation relation, but here we do not want to reverse the orientation of the edges that induces the vertex-orientation of $\Gamma$.
\begin{lemma}
\label{lemgenhol}
Let $\Gamma \in \CS^u_n$. An element $\rho$ of $H^1(\Gamma;\ZZ)$ can be identified with a beaded graph $\Gamma_{\rho}$ with support $\Gamma$ (equipped with its canonical vertex-orientation)
and with beads in $\{t^n;n\in \ZZ \}$, up to the generalized holonomy relation.
\end{lemma}
\bp 
Assume that a graph with support $\Gamma$ and with beads in $\{t^n;n\in \ZZ \}$ is given. Consider a cycle of $\Gamma$. This cycle is an algebraic sum of oriented edges.
Map it to the corresponding algebraic sum of the edge bead exponents. The obtained morphism of $H^1(\Gamma;\ZZ)$ does not depend on the representative of the beaded graph in his class modulo generalized holonomy. 
Conversely, when an element of $H^1(\Gamma;\ZZ)$ is given, consider a maximal tree in each component of $\Gamma$ and bead its edges by $1$, then bead the remaining edges so that they induce the given morphism of $H^1(\Gamma;\ZZ)$ by the above process. The class of the obtained beaded graph modulo generalized holonomy relation is independent of the chosen maximal tree, since the generalized holonomy relation allows one to modify the edge beads in any subtree of $\Gamma$, arbitrarily.
\eop

Let $(H_i)_{i=1,\dots,3n}$ be a family of integral cycles of $(\tilde{C}_2(\EXT),\partial \tilde{C}_2(\EXT))$.

We first define the intersection $I_{\Gamma}(\{H_i\}_{i=1,\dots,3n})=I_{\Gamma}(\{H_i\})$ of the 
$H_i$ with respect to a graph $\Gamma$ of $\CS^u_n$.
Consider $\Gamma \in \CS^u_n$.

Let $e=e(i)$ be the edge of $\Gamma$ numbered by $i$ that goes from  $v(j)$ to $v(k)$. Consider the map $p(\Gamma,i)=p(\Gamma,e(i))\colon C_{2n}(\EXT) \rightarrow C_2(\EXT)$ that lifts 
$(m_1, \dots, m_{3n}) \mapsto (m_j,m_k)$, continuously.
Let $H_i(\Gamma) \subset C_{2n}(\EXT)$ be defined as
$$H_i(\Gamma)=p(\Gamma,i)^{-1}\left(p_C(H_i)\right).$$ 
Then $H_i(\Gamma)$ is a codimension $2$ {\em chain\/} (here, a combination of simplices, or of manifolds with boundaries) whose points (that project to) $(m_1, \dots, m_{3n})$ are continuously equipped with rational homology classes of paths from $m_j$ to $m_k$. It is cooriented in $C_{2n}(\EXT)$ (that is oriented like $\EXT^{2n}$) by the coorientation of $H_i$ in $C_2(\EXT)$.

Then the intersection of the $H_i(\Gamma)$, for $i \in \{1,2,\dots, 3n\}$, is a compact subspace $I(\Gamma,\{H_j\})$ of $C_{2n}(\EXT)$. Its image $p(\Gamma,i)(I(\Gamma,\{H_j\}))$ in $H_i$ is a compact subset of $H_i$.

\begin{definition}
 \label{defcyclegen}
The $H_i$ are said to be in {\em general position with respect to $\Gamma$\/} if
\begin{itemize}
\item $I(\Gamma,\{H_j\})$ is finite,
\item $p(\Gamma,i)(I(\Gamma,\{H_j\}))$ is made of points in the interiors of the $4$-manifolds with boundaries of $H_i$, and, 
\item all the $H_i(\Gamma)$ intersect transversally at the points of $I(\Gamma,\{H_j\})$ that are in the interior of $C_{2n}(\EXT)$.\end{itemize}
\end{definition}

Assume that our $H_i$ are in general position with respect to $\Gamma$.
Then
the cooriented $H_i(\Gamma)$ only intersect transversally at distinct points 
in the interior of $C_{2n}(\EXT)$.
We define their equivariant algebraic intersection $I_{\Gamma}(\{H_i\}) \in \CA^{h}_n(\delta)$ as follows.

Consider an intersection point $m$, such that $p_b(m)=(m_1, \dots, m_{2n}) \in \EXT^{2n}$, it is equipped with a sign $\varepsilon(m)$ as usual, and it is also equipped with the following additional data:
Associate $m_j$ with $v(j)$. Since $m$ belongs to $H_i(\Gamma)$, the edge $e(i)$ from $v(j)$ to $v(k)$ is equipped with a rational homology class of paths from $m_j$ to $m_k$. Each edge of $\Gamma$ is equipped with a rational homology class of paths between its ends in this way. By composition, each cycle of $\Gamma$ is equipped with an element of $H_1(\EXT;\ZZ)/\mbox{Torsion}$. Since $H_1(\EXT;\ZZ)/\mbox{Torsion}=\ZZ m(K)$, the isomorphism from $H_1(\EXT;\ZZ)/\mbox{Torsion}$ to $\ZZ$ that maps $m(K)$ to $1$ associates an integer with each cycle of $\Gamma$. Thus $m$ defines a map from $H_1(\Gamma;\ZZ)$ to $\ZZ$, that is an element $\rho(m)$ of $H^1(\Gamma;\ZZ)$.

Assign $\varepsilon(m)[\Gamma_{\rho(m)}] \in \CA^{h}_n(\delta)$ to $m$, with the notation of Lemma~\ref{lemgenhol}.

Then define $I_{\Gamma}(\{H_i\}) \in \CA^{h}_n(\delta)$ as the sum over the intersection points of the $\varepsilon(m)[\Gamma_{\rho(m)}]$.

Note that $I_{\Gamma}(\{H_i\})$ does not depend on the numbering of the vertices of $\Gamma$. 
Also note that when some $H_i$ is changed to $P(\thetap) H_i$ for a Laurent polynomial $P$ with integral coefficients, the bead on Edge number $i$ in every $\Gamma_{\rho(m)}$ is multiplied by $P(t)$.

The definition of $I_{\Gamma}(\{H_i\})$ naturally extends to cycles $H_i$ of $(C_2(\EXT),\partial C_2(\EXT))$ with coefficients in $\frac1{\delta(t)}\QQ[t^{\pm 1}]$ as follows.
First recall that a cycle of $(C_2(\EXT),\partial C_2(\EXT))$ with coefficients in $\QQ[t^{\pm 1}]$ is nothing but a rational cycle of
$(\tilde{C}_2(\EXT),\partial \tilde{C}_2(\EXT))$, where the multiplication by $t$ stands for the action of $\thetap$. 
Cycles with coefficients in $\frac1{\delta(t)}\QQ[t^{\pm 1}]$ are obtained from these by tensoring
this $\QQ[t^{\pm 1}]$-module by $\frac1{\delta(t)}\QQ[t^{\pm 1}]$ over $\QQ[t^{\pm 1}]$.
For any collection of cycles $H_i$ of $(C_2(\EXT),\partial C_2(\EXT))$ with coefficients in $\frac1{\delta(t)}\QQ[t^{\pm 1}]$, there is  an integral number $k$ such that $k\delta(t) H_i$ is an integral chain of $\tilde{C}_2(\EXT)$ for all $i$. 
Cycles $H_i$ of $(C_2(\EXT),\partial C_2(\EXT))$ with coefficients in $\frac1{\delta(t)}\QQ[t^{\pm 1}]$ are said to be in {\em general $3n$-position with respect to $\Gamma$\/} if for some $k \in \NN$, the $k\delta(t)H_i$ are integral cycles of $(\tilde{C}_2(\EXT),\partial \tilde{C}_2(\EXT))$ in general position with respect to $\Gamma$.
For such a collection of cycles, define $I_{\Gamma}(\{H_i\})$ as the element of $\CA^{h}_n(\delta)$ obtained from $I_{\Gamma}(\{k\delta(t) H_i\})$ by dividing the bead of each edge by $k\delta(t)$. 

Cycles $H_i$ of $(C_2(\EXT),\partial C_2(\EXT))$ with coefficients in $\frac1{\delta(t)}\QQ[t^{\pm 1}]$ are said to be in {\em general $3n$-position\/} if for some $q\in \QQ$, the $q\delta(t)H_i$ are integral cyles of $(\tilde{C}_2(\EXT),\partial \tilde{C}_2(\EXT))$ in general position with respect to $\Gamma$, for any $\Gamma \in \CS^u_n$.
For such a collection of cycles, define the {\em $\CA^{h}_n$-intersection\/} of the $H_i$ as

$$I_n(\{H_i\})=\sum_{\Gamma \in \CS^u_n}\frac{I_{\Gamma}(\{H_i\})}{2^{3n}(3n)!(2n)!}.$$

Considering only the connected graphs as in $\cite{lesbonn}$, define the {\em connected
$\CA^{h}_n$-intersection\/} of the $H_i$ as
$$I^c_n(\{H_i\})=\sum_{\Gamma \in \CS_n}\frac{I_{\Gamma}(\{H_i\})}{2^{3n}(3n)!(2n)!}.$$

\subsection{Equivariant propagators of $\EXT$}
\label{subcyc}

The \emph{equivariant propagators} $\eqprop_i(\tau)=\frac1{\delta(t)}\deleqprop_i(\tau)$ of $(C_2(\EXT),\partial C_2(\EXT))$ with coefficients in $\frac1{\delta(t)}\QQ[t^{\pm 1}]$ that will be used to define $\tilde{Z}_n(K,\tau)$ as $I_n(\{\eqprop_i(\tau)\})$ can be thought of as parallel to each other.
We are going to give conditions on $\deleqprop_1(\tau)=\deleqprop(\tau)$ depending on parameters $\cvarM=\cvarM_1 \in S^2$
and $\sK=\sK_1 \in S^1$ to constrain $\deleqprop(\tau)$, and the $\deleqprop_i(\tau)$ will be defined by replacing the parameter $\cvarM$ by $\cvarM_i \in S^2$ and $\sK$ by $s_i \in S^1$.

For $\alpha \subset [0,4]$, let $D_{\alpha}=\{z \in \CC;|z| \in \alpha\}$. Fix a collar neighborhood $\overline{\EXT_{]1/2,4]}}$ of $ \partial \EXT$ that reads as the product of the annulus $D_{[1/2,4]}$ by the circle $S^1$ so that, for
$s$ in the circle $D_{\{1\}}$, $s \times S^1$ is a meridian of $K$ and, for $z \in S^1$, $D_{\{1\}} \times \{z\}$ is a parallel of $K$ that bounds in $\EXT$.

The product $D_{[0,4]} \times \RR$ embeds in $\CC \times \RR$, naturally
where $\CC$ is thought of as horizontal in $\RR^3=\CC\times \RR$ where $\RR$ is vertical.

The standard parallelization of $\RR^3$ induces a parallelization $\tau_{\pi}$ of $D_{[0,4]} \times \RR$. Since  $\tau_{\pi}$ is equivariant
under the translations by the upward unit vector of $\RR^3$, it induces a parallelization on the quotient $D_{[0,4]} \times S^1$ that is still denoted by $\tau_{\pi}$, and that is the natural product parallelization.

Embed $p_{\EXT}^{-1}(\overline{\EXT_{]1/2,4]}})=D_{[1/2,4]} \times \RR$ in $\CC \times \RR$, naturally.

Fix a parallelization $\tau \colon \EXT \times \RR^3 \rightarrow T\EXT$ of $\EXT$ that coincides with the parallelization $\tau_{\pi}$ on $T\EXT_{]1/2,4]}$. (There is no obstruction to doing so since $D_{\{1\}} \times \{z\}$ vanishes in $H_1(\EXT;\ZZ)$.)

The vector $\cvarM$ of $S^2$ defines a section $s_{\tau}(\EXT;\cvarM)=\tau(\EXT \times \cvarM)$ of the unit tangent bundle $U\EXT$ to
$\EXT$ that sits in $\partial C_2 (\EXT)$.

Consider a map $\funcE \colon \EXT \rightarrow S^1$ that coincides with the projection
to $S^1$ on $D_{[1/2,4]} \times S^1$, and a lift of this map $\functE \colon \tilde{\EXT} \rightarrow \RR$, through the standard covering $(t \mapsto \exp(2i\pi t))$ such that
the projection to $\RR$ of $\tilde{\EXT}_{]1/2,4]}=D_{]1/2,4]} \times \RR$ is $\functE$.

Let 
$$\begin{array}{llll}
   r \colon & \EXT & \rightarrow & [1/2,4]\\
  & x \in \EXT\setminus (D_{]1/2,4]} \times S^1) & \mapsto & 1/2\\
& (z_D,z) \in D_{[1/2,4]} \times S^1&\mapsto &|z_D|.\\
\end{array}
$$
When $\alpha \subset [0,4]$, set $\EXT_{\alpha}=r^{-1}(\alpha)$
and $\tilde{\EXT}_{\alpha}=p_{\EXT}^{-1}(\EXT_{\alpha})$ so that $\EXT_{[0,1]}=\EXT\setminus (D_{]1,4]} \times S^1)$.

Fix $\varepsilon \in ]0,1/8]$. 
Let $$\begin{array}{llll}\chi \colon &[-4,4]& \rightarrow& [0,1]\\
&t \in [-\varepsilon,4] &\mapsto &1\\
&t \in [-4,-2\varepsilon] &\mapsto &0
\end{array}$$
be a smooth map. Recall that $\tilde{\EXT}_{]1,4]}$ is embedded in $\RR^3$ that is seen as $\CC \times \RR$.
When $(u,v) \in \left(\tilde{\EXT}_{[0,4]}^2 \setminus \tilde{\EXT}_{[0,2[}^2\right)$, set
$$U(u,v)=(1-\chi(r(u) -r(v)))(0,\functE(u)) + \chi(r(u)-r(v))u$$
$$V(u,v)=(1-\chi(r(v) -r(u)))(0,\functE(v)) + \chi(r(v)-r(u))v$$
where $0u=0 \in \RR^3$, for any $u \in \EXT$,
so that $(U(u,v),V(u,v)) \in (\RR^3)^2$.
Define $$\begin{array}{llll}\pi \colon &p_C^{-1}\left((\EXT_{[0,4]}^2 \setminus \EXT_{[0,2[}^2)\setminus \mbox{diag}({\EXT}_{[2,4]}^2) \right)&\rightarrow& S^2\\
 &\overline{(u,v)}&\mapsto & \frac{V(u,v)-U(u,v)}{\parallel V(u,v)-U(u,v)\parallel}.
         \end{array}$$

The map $\pi$ naturally extends to $\tilde{C}_2(\EXT_{[0,4]})\setminus \tilde{C}_2(\EXT_{[0,2[})$.

The following proposition is proved in \cite[Section 12, Propositions 12.2 and 12.4]{lesbetaone}.

\begin{proposition}
 \label{propdefbord}
Set $$J_{\Delta}=J_{\Delta}(t)=\frac{t\Delta^{\prime}(t)}{\Delta(t)}.$$
Let $S^2_H$ denote the subset of $S^2$ made of the vectors whose vertical coordinate is in $]0,\frac{1}{50}[$.
Let $\cvarM \in S^2_H$, let $\sK \in D_{\{1\}}$,
then
there exists a $4$-dimensional rational chain $C(\cvarM,\sK,\tau)$ of $\tilde{C}_2(\EXT_{[0,2]})$ whose boundary is
$$\delta(K)(\thetap)\left(\pi_{|\partial \tilde{C}_2(\EXT_{[0,2]}) \setminus \partial \tilde{C}_2(\EXT_{[0,2[})}^{-1}(\cvarM) \cup s_{\tau}(\EXT_{[0,2]};\cvarM) \cup J_{\Delta}(\thetap) U\EXT_{|{\sK \times S^1}}\right),$$
and that is transverse to $\partial \tilde{C}_2(\EXT_{[0,2]})$.
\end{proposition}

\begin{lemma}
\label{lemgenpos}
There exist
\begin{itemize}
\item $\varepsilon \in ]0,\frac{1}{2n}[$, $\funcE$ (involved in the above definition of the map $\pi$), $k\neq 0 \in \NN$,
\item $3n$ distinct elements $\sK_1, \dots, \sK_{3n}$ of $D_{\{1\}}$,
\item a {\em regular\/} $3n$-tuple $(\cvarM_1, \dots,\cvarM_{3n})$  of $(S^2_H)^{3n}$, where {\em regular\/} means in some open dense subset of $(S^2_H)^{3n}$, that is specified in \cite[Definition 2.7]{lesbonn}. (A more restrictive notion of {\em regular\/} is defined in Subsection~\ref{subuplegen} below, it still defines a dense open subset.)
\item $3n$ integral chains (integral combinations of $C^{\infty}$ properly embedded $4$-manifolds with boundaries) $C(\cvarM_i,\sK_i,\tau)$ of $\tilde{C}_2(\EXT_{[0,2]})$  in {\em general $3n$-position\/} whose boundaries are $$k\delta(K)(\thetap)
\left(\pi_{|\partial \tilde{C}_2(\EXT_{[0,2]}) \setminus \partial \tilde{C}_2(\EXT_{[0,2[})}^{-1}(\cvarM_i) \cup s_{\tau}(\EXT_{[0,2]};\cvarM_i) \cup J_{\Delta}(\thetap) U\EXT_{|{\sK_i \times S^1}}\right),$$ respectively.
\end{itemize}
\end{lemma}
\bp
It is proved in \cite[Sections 2.3, 2.4, 2.5]{lesbonn}, with the difference that general $3n$-position
only means general position with respect to connected graphs in \cite{lesbonn}.
However, if any subfamily of cardinality $3k$ of the $C_i(\tau)$ is in general position with respect to all connected graphs of $\CS_k$, then the $C_i(\tau)=C(\cvarM_i,\sK_i,\tau)$ will be in general position with respect to all graphs in $\CS_n^u$, and it is proved in \cite[Proposition 2.9]{lesbonn} that this is true when $(C_i(\tau))_{i\in \{1,2,\dots,3n\}}$ belongs to a finite intersection of open dense sets.
\eop

\begin{theorem}
\label{thminvwithtau}
Under the assumptions of Lemma~\ref{lemgenpos}, set $$\eqprop_i(\tau)=\frac{1}{k\delta(K)(t)}C(\cvarM_i,\sK_i,\tau) + \overline{\pi^{-1}(\cvarM_i)} \subset C_2(\EXT).$$
Then neither
$ I_n(\{\eqprop_i(\tau)\})$ nor $ I^c_n(\{\eqprop_i(\tau)\})$ depends on the equivariant propagators $\eqprop_i(\tau)$ satisfying the conditions above. Therefore, they are invariants of 
$(K,\tau)$, that are denoted by $\tilde{Z}_n(K,\tau)$ and $\tilde{z}_n(K,\tau)$, respectively.
Furthermore, $$(\tilde{Z}_n(K,\tau))_{n\in \NN}=\exp\left(\sum_{k\geq 1}\tilde{z}_k(K,\tau)\right)$$ in $\CA^h(\delta(K))$.
\end{theorem}
\bp According to \cite[Theorem 2.4]{lesbonn}, $ I^c_n(\{\eqprop_i(\tau)\})$ is an invariant $$\tilde{z}_n(K,\tau)=\sum_{\Gamma \in \CS_n}\frac{I_{\Gamma}(\{\eqprop_i(\tau)\})}{2^{3n}(3n)!(2n)!}$$ 
of $(K,\tau)$. 
Thus it is enough to prove that
$$\left(I_n(\{\eqprop_i(\tau)\})\right)_n=\exp\left(1+\sum_{k\geq 1}\tilde{z}_k(K,\tau)\right).$$
Let $\CP(n)$ be the set of partitions of $n$, that is the set of uples $((n_1,k_1);(n_2,k_2);\dots;(k_r,n_r))$ where the $k_j$ and the $n_j$ are positive integers
such that $\sum_{j=1}^rk_jn_j=n$. $0<n_1 < n_2 <\dots < n_r$. Then
$I_n(\{\eqprop_i(\tau)\})$ splits as $$I_n(\{\eqprop_i(\tau)\})=\sum_{p \in \CP(n)}I_n(p,\{\eqprop_i(\tau)\})$$
where $I_n(p,\{\eqprop_i(\tau)\})$ is a combination of beaded graphs with $k_j$ connected components with $2n_j$ vertices for every $j$. Furthermore, for such a given $p$, 
$I_n(p,\{\eqprop_i(\tau)\})$ splits as $$I_n(p,\{\eqprop_i(\tau)\})=\sum_{q \in \CQ(p)}I_n(p,q,\{\eqprop_i(\tau)\})$$ over the set $\CQ(p)$ of partitions $q$ of $\{1,\dots,3n\}$ that have $k_j$ subsets $A(\ell,n_j)$ of cardinality $(3n_j)$, for every $j$. 

Let $\CS^{\prime}_n$ (resp. $\CS^{u\prime}_n$) be the set of graphs of $\CS_n$ (resp. of $\CS^{u}_n$) where the vertices are not numbered anymore. Then
$\tilde{z}_n(K,\tau)=\sum_{\Gamma \in \CS^{\prime}_n}\frac{I_{\Gamma}(\{\eqprop_i(\tau)\})}{2^{3n}(3n)!}$.
For $p=((n_1,k_1);(n_2,k_2);\dots;(k_r,n_r)) \in \CP(n)$, and for $q \in \CQ(p)$,
let $\CS^{\prime}_n(p,q)$ be the set of graphs $\Gamma$ of $\CS^{u\prime}_n$ such that two edges $e(k)$ and $e(\ell)$ are in the same connected component of $\Gamma$ if and only if $k$ and $\ell$ are in the same component of $q$.

According to \cite[Theorem 2.4]{lesbonn}, 
$\tilde{z}_{n_j}(K,\tau)=I^c_{n_j}(\{\eqprop_i(\tau)\}_{i \in A(\ell,n_j)} )$, for every $A(\ell,n_j)$. Furthermore, when $\Gamma =\Gamma_1 \coprod \Gamma_2$, $I_{\Gamma}(\{\eqprop_i(\tau)\})=I_{\Gamma_1}(\{\eqprop_i(\tau)\}_{e(i)\in \Gamma_1})I_{\Gamma_2}(\{\eqprop_i(\tau)\}_{e(i)\in \Gamma_2})$.
Thus
$$I_n(p,q,\{\eqprop_i(\tau)\})=\sum_{\Gamma \in \CS^{\prime}_n(p,q)}\frac{I_{\Gamma}(\{\eqprop_i(\tau)\})}{2^{3n}(3n)!}=\frac{\prod_{j=1}^r (3n_j)!^{k_j}}{(3n)!}\prod_{j=1}^r\tilde{z}_{n_j}(K,\tau)^{k_j}.$$
In particular, $I_n(p,q,\{\eqprop_i(\tau)\})$ does not depend on $q$, and since the number of partitions $q$ of $\CQ(p)$ is $\frac{(3n)!}{\prod_{j=1}^r\left(k_j!(3n_j)!^{k_j}\right)}$,
we find that
$$I_n(p,\{\eqprop_i(\tau)\})=\frac{1}{\prod_{j=1}^rk_j!}\prod_{j=1}^r\tilde{z}_{n_j}(K,\tau)^{k_j}.$$
\eop

\begin{remark}
 As in \cite[Proposition 3.3 and Theorem 4.8]{lesbetaone}, the $\eqprop_i(\tau)$ represent the equivariant linking form in the sense that for a two-component link $(J,L)$ of $\tilde{\EXT}$ such that
$p_{\EXT}(J) \cap p_{\EXT}(L)=\emptyset$,
$$lk_e(J,L)= \langle L \times J, \eqprop_i(\tau) \rangle_e =\langle \eqprop_i(\tau),J \times L  \rangle_e.$$
\end{remark}

Consider a $4$-dimensional manifold $W$ with signature $0$ and with boundary and ridges such that $$\partial W = -(D_{[0,4]} \times S^1) \cup_{\{0\}\times \partial D_{[0,4]} \times S^1} (-[0,1] \times \partial D_{[0,4]} \times S^1) \cup_{\{1\}\times \partial \EXT} (\{1\}\times \EXT)$$
and $W$ is identified with an open subspace of one of the products $[0,1[ \times D_{[0,4]} \times S^1$ or $]0,1] \times \EXT$ near $\partial W$.
Recall that $\tau_{\pi}$ is the standard product parallelization of $D_{[0,4]} \times S^1$. Then the parallelizations $\tau$ and $\tau_{\pi}$ together induce a trivialization $\tau_{\partial W}$ of $TW$ on $\partial W$ and $p_1(\tau_{\pi},\tau)$ is the relative first Pontrjagin class of $\tau_{\partial W}$ in $(W,\partial W)$. The definition and properties of $p_1$ are given in details in Subsection~\ref{subpont}.

\begin{theorem}
\label{thmhighloop}
Let $\ansothree_n=-\xi_n$ be the element of $\CA_n$ defined in Definition~\ref{defxin} or in \cite[Subsection 1.6]{lesconst}. It
is zero when $n$ is even and $\ansothree_1= \frac{1}{12}[\theta]$.
Then $$\tilde{z}_n(K)=\tilde{z}_n(K,\tau)-\frac{p_1(\tau_{\pi},\tau)}{4}\ansothree_n$$ is an invariant of $K$ for any $n>0$. Set $\tilde{z}_0(K)=0$ and $\tilde{z}(K)=\sum_{n\in \NN}\tilde{z}_n(K)$. Then
$\tilde{Z}(K)=\exp(\tilde{z}(K))$
is an invariant of $K$ such that $$\tilde{Z}(K)=\tilde{Z}(K,\tau)\exp(-\frac{p_1(\tau_\pi,\tau)}{4}\ansothree).$$
\end{theorem}
\bp This is a consequence of \cite[Theorem 2.4]{lesbonn}. According to \cite[Proposition 2.45]{lesconst}, $\xi_1= -\frac{1}{12}[\theta]$, and 
according to \cite[Proposition 1.10]{lesconst}, $\xi_n=0$ if $n$ is even. \eop

\section{Preliminaries for the proof of Theorem~\ref{thmflag}}
\setcounter{equation}{0}
\label{secprelproof}

We refer to Section~\ref{secrecdef} above for the construction of $(\tilde{Z}_n)_n$.
What is actually constructed is an invariant $\tilde{Z}(K,\tau)=(\tilde{Z}_n(K,\tau))_{n \in \NN}$ for some parallelization $\tau$ of $\EXT$ that is fixed on a collar of $\partial \EXT$ such that
$$\tilde{Z}(K)= \tilde{Z}(K,\tau)\exp\left(-\frac{p_1(\tau_{\pi},\tau)}{4} \ansothree\right)$$ 
where 
$$\ansothree=\sum_{n \in \NN} \ansothree_{2n+1} \in \CA.$$

Fix $\EXT$ and $n$. Since the data of $\EXT$ and $m(K)$ are equivalent to the datum of $K \subset R$, we write 
$\tilde{z}_n(\EXT)=\tilde{z}_n(K)$.
Under the assumptions of Theorem~\ref{thmflag},
we compute $$\tilde{Z}_n([\EXT;(A^{(i)\prime}/ A^{(i)})_{i\in N}])=\sum_{S \subset\{1,2,\dots,2n\}}(-1)^{\sharp S} \tilde{Z}_n(\EXT_S).$$

In order to do this, we use $3n$ parallel cycles $\deleqprop_r(\EXT_S)$ of $(\tilde{C}_2(\EXT_S),\partial \tilde{C}_2(\EXT_S))$  similar to the chains $\delta(t)\eqprop_V(M_S)$ constructed in \cite[Section 11]{lesbetaone}.
We first review the properties of the \emph{equivariant propagators} $\eqprop_r(\EXT_S)=\frac1{\delta(t)} \deleqprop_r(\EXT_S)$.
In order to construct them, we assume that the $A^{(i)}$ are inside $\EXT_{[0,1/2]}$.
Recall that $p_{\EXT}\colon \tilde{\EXT} \rightarrow \EXT$ is the covering map of the infinite cyclic covering for $\EXT$. Note that $\delta(K_S)$ is independent of $S$. Denote it by $\delta$ or $\delta(K)$.

\subsection{The equivariant propagator $\eqprop(\EXT)$}

Let $(a^{(i)}_j,z^{(i)}_j)_{j=1, \dots g(A^{(i)})}$ be a basis of $H_1(\partial A^{(i)}) $ such that $a^{(i)}_j=\partial \Sigma(a^{(i)}_j)$ where $\Sigma(a^{(i)}_j)$ is a rational chain of $A^{(i)}$ and $\langle a^{(i)}_j,z^{(i)}_k \rangle=\delta_{jk}$.

Let $[-7,4] \times \partial A^{(i)} \subset \EXT$ be a neighborhood of $\partial A^{(i)}= 0 \times\partial A^{(i)}$. ($[-7,0] \times \partial A^{(i)} \subset A^{(i)}$.)

For $s \in [-7,4]$, set
$$A^{(i)}_s=\left\{\begin{array}{ll} A^{(i)} \cup  ([0,s] \times \partial A^{(i)} ) & \mbox{if}\;\; s \geq 0\\
A^{(i)} \setminus ( ]s,0]\times \partial A^{(i)} ) & \mbox{if} \;\;s \leq 0\end{array} \right.$$
and $(\EXT \setminus A^{(i)})_s = \EXT \setminus \mathring{A^{(i)}}_s$.
$$ \partial A^{(i)}_s =\{s\}  \times\partial A^{(i)}=-\partial (\EXT \setminus A^{(i)})_s.$$

Assume that $\Sigma(a^{(i)}_j)$ intersects $ [-4,0] \times \partial A^{(i)} $ as
$[-4,0]\times a^{(i)}_j $ and construct a rational chain 
$\Sigma_4(a^{(i)}_j)=\Sigma(a^{(i)}_j)\cup [0,4] \times a^{(i)}_j $.
Set $\Sigma_s(a^{(i)}_j)=A^{(i)}_s \cap \Sigma_4(a^{(i)}_j)$.

Recall that we have fixed preferred lifts for the $A^{(i)}$ in $\tilde{\EXT}$, the preferred lifts of the curves $a_j^{(i)}$, and the chains $\Sigma(a_j^{(i)})$ in the preferred lift of $A^{(i)}$ are also denoted by 
$a_j^{(i)}$  or $\Sigma(a_j^{(i)})$, while the preferred lifts of the curves $z_j^{(i)}$ will still be denoted by $\tilde{z}_j^{(i)}$.

Assume that the map $\funcE \colon \EXT \rightarrow S^1$ is constant on each $A^{(i)}$ and that the $\funcE(A^{(i)})$ are pairwise distinct. There is no loss since the restriction of $\funcE$ to $A^{(i)}$ is null-homotopic. The functions $\func_{\EXT_S}$ will be assumed to coincide with $\funcE$ outside $\cup_{i \in S}A^{(i)\prime}$ and to map $A^{(i)\prime}$ to $\funcE(A^{(i)})$.

Let $p(A^{(i)})$ be a point of $\partial A^{(i)}$ outside the $a^{(i)}_j$ and the $z^{(i)}_j$.
Let $\gamma^+(A^{(i)})$ (resp. $\gamma^-(A^{(i)})$) denote a path in $p_{\EXT}^{-1}((\EXT \setminus A^{(i)})_{-4} \setminus \coprod_{j\neq i}A_4^{(j)})$ that intersects $\left(p_{\EXT}^{-1}(A^{(i)}_4)\subset \tilde{\EXT}\right)$ as $( [-4,4] \times p(A^{(i)}))$, and that intersect $\tilde{\EXT}_{[1,4]}\subset \CC \times \RR$ as the half-line through $\{0\}\times \tilde{\funcE}(A^{(i)})$ directed by $\cvarM$ (resp. by $-\cvarM$) does.

Any path $\gamma^{\pm}(A^{(i)})$ is supposed to be disjoint from all the paths $\gamma^{\pm}(A^{(j)})$ for $j \neq i$.
For $s\in [-4,4]$, $\gamma^{\pm}_s(A^{(i)}) =\gamma^{\pm}(A^{(i)})\cap p_{\EXT}^{-1}((\EXT \setminus A^{(i)})_s)$.

For $k=1, \dots, g(A^{(i)})$, $\delta(K)(\theta_{\EXT})\tilde{z}^{(i)}_k$ bounds a rational chain in $\tilde{\EXT}$. Therefore, it cobounds a rational chain $\delta(K)(t)\Sigma(\check{z}^{(i)}_k)$ in $\left(\tilde{\EXT} \setminus p_{\EXT}^{-1}(\mathring{A}^{(i)})\right)$ with a combination of $a^{(i)}_{\ell}$ that has its coefficients in $\QQ[t,t^{-1}]$.
$$\partial \delta(K)(t)\Sigma(\check{z}^{(i)}_k)= \delta(K)(\theta_{\EXT})\tilde{z}^{(i)}_k - \sum_{\ell=1}^{g(A^{(i)})} \left(\delta(K) lk_e(\tilde{z}^{(i)}_k,\tilde{z}^{(i)-}_{\ell})\right)(\theta_{\EXT})\left(a^{(i)}_{\ell}\right)=\delta(K)(t)\check{z}^{(i)}_k.$$
Furthermore assume that, for $i\neq j$,
$$\delta(t)\Sigma(\check{z}^{(i)}_k) \cap p_{\EXT}^{-1}(A_4^{(j)})=\sum_{\ell=1}^{g(A^{(j)})}\left(\delta(K) lk_e(\tilde{z}^{(i)}_k,\tilde{z}^{(j)}_{\ell})\right)(\theta_{\EXT})\left(\Sigma_4(a^{(j)}_{\ell})\right).$$
Again use the product structures to extend $\delta(t)\Sigma(\check{z}^{(i)}_k)$ on $[-4,0] \times \partial A^{(i)}$ by setting
$$\delta(t)\Sigma_u(\check{z}^{(i)}_j)= \delta(t)\Sigma(\check{z}^{(i)}_j) \cup -([u,0] \times \partial \delta(t)\Sigma(\check{z}^{(i)}_j)) \subset \left(\widetilde{\EXT} \setminus p_{\EXT}^{-1}( \mathring{A_u^{(i)}})\right)$$
for $u \in[-4,0]$.
Note that $C_2(A_4^{(j)})$ has a preferred lift in $\tilde{C}(\EXT)$, namely the lift that contains $\{0\} \times UA_4^{(j)}$.
In the following proposition that is proved in \cite[Proposition 11.8]{lesbetaone}, $A^{(j)}_t \times p(A^{(j)})$, $\Sigma_t(a_i^{(j)})\times {z}_k^{(j)}$, $p(A^{(j)}) \times A^{(j)}_t$ and ${z}^{(j)}_k \times \Sigma_t(a_i^{(j)})$ are in the preferred lift of $C_2(A^{(j)}_4)$.

\begin{proposition}
The chain $\deleqprop(\EXT)$ can be assumed to satisfy the following properties.
\begin{itemize}
\item For any $j=1,\dots, 2n$, for any $t\in[-4,0]$,
$\deleqprop(\EXT)$ intersects 
$$A^{(j)}_t \times_{e} (\EXT\setminus A^{(j)})_{t+3} =p_{C}^{-1}(A^{(j)}_t \times (\EXT\setminus A^{(j)})_{t+3})$$ as 

$$\delta(\thetap)(A^{(j)}_t \times \gamma^+_{t+3}(A^{(j)}))+\sum_{i \in \{1, \dots, g(A^{(j)})\}} \Sigma_t(a_i^{(j)})\times \delta(t)\Sigma_{t+3}(\check{z}_i^{(j)}).$$

\item  For any $j=1,\dots, 2n$, for any $t\in[-4,0]$,
$\deleqprop(\EXT)$ intersects $(\EXT\setminus A^{(j)})_{t+3} \times_{e} A^{(j)}_t$ as
$$\sum_{i \in \{1, \dots, g(A^{(j)})\}} \delta(t^{-1}) \Sigma_{t+3}(\check{z}^{(j)}_i)\times \Sigma_t(a_i^{(j)})
-\delta(\thetap)(\gamma^-_{t+3}(A^{(j)}) \times A^{(j)}_t).$$
\item  For any $j=1,\dots, 2n$, for any $i =1,\dots, g(A^{(j)})$,
$\langle \deleqprop(\EXT), \Sigma_3(a^{(j)}_i)\times p(A^{(j)}) \rangle_e=0$ and
$$\langle \deleqprop(\EXT), p(A^{(j)})\times \Sigma_3(a^{(j)}_i) \rangle_e=0.$$
\item $\deleqprop(\EXT)=\partial \deleqprop(\EXT) \times [0,1]$ in an equivariant neighborhood $\partial \tilde{C}_2( \mathring{\EXT} ) \times [0,1]$ of $\partial \tilde{C}_2( \mathring{\EXT} )$ in $\TCE$, where the coordinate in $[0,1]$ is thought of as the distance between two points in a configuration near $U\EXT$.
\end{itemize}
\end{proposition}
We therefore assume that $\deleqprop(\EXT)$ satisfies all the above properties.

\subsection{Defining $\deleqprop(\EXT_S)$ from $\deleqprop(\EXT)$}

We consider parallelizations $\tau$ of $\EXT$ that are standard on $\EXT_{[1,4]}$. See Subsection~\ref{subcyc}.
When the $A^{(i)}$ and the $A^{(i)\prime}$ are integral homology handlebodies, the restriction of such a parallelization $\tau$ of $\EXT$ to $\EXT \setminus \coprod_{i=1}^{2n}\mathring{A}{(i)}$ extends as a parallelization $\tau(\EXT_S)$
on $\EXT_S$ for every $S \subset N$ where $N=\{1,\dots,2n\}$.
It does not necessarily extend when $A^{(i)}$ and $A^{(i)\prime}$ are rational homology handlebodies. 
In order to prove Theorem~\ref{thmflag} in general, we use the pseudo-parallelizations introduced in \cite[Sections 4.2 and 4.3]{lessumgen} as in \cite[Section~10]{lesbetaone}.

These {\em pseudo-parallelizations\/} $\tilde{\tau}$ defined in Definition~\ref{defpseudotriv}
\begin{itemize}
\item generalize parallelizations,
\item always extend to rational homology handlebodies (Lemma~\ref{lem102lesbetaone}),
\item induce genuine trivializations $\tilde{\tau}_{\CC}$ of $T\EXT \otimes_{\RR} \CC$ that have Pontrjagin numbers $p_1(\tau_{\pi},\tilde{\tau}_{\CC})$ (see Definition~\ref{defpseudotrivpone}),
\item define $3$-dimensional pseudo-sections $s_{\tilde{\tau}}(\EXT;\cvarM)$ of $U\EXT$, for $\cvarM \in S^2$, so that $$\langle s_{\tilde{\tau}}(\EXT;\fvarM) \cap U\EXT_{|\Sigma}, s_{\tilde{\tau}}(\EXT;\svarM)\rangle_{U\EXT} =0$$ for a two-dimensional cycle $\Sigma$ when $\fvarM \neq \svarM$ (see Definition~\ref{defpseudosec} and \cite[Lemma 10.5]{lesbetaone}),
\item provide $3$-dimensional bounding cycles $\partial(\cvarM,\sK,\tilde{\tau})$ of $\partial \tilde{C}_2(\EXT_{[0,2]})$, for $\sK \in D_{\{1\}}$ and $\cvarM \in S^2_H$,
$$\partial(\cvarM,\sK,\tilde{\tau})=\delta(K)(\thetap)\left(\pi_{|\partial \tilde{C}_2(\EXT_{[0,2]}) \setminus \partial \tilde{C}_2(\EXT_{[0,2[})}^{-1}(\cvarM) \cup s_{\tilde{\tau}}(\EXT_{[0,2]};\cvarM) \cup J_{\Delta}(\thetap) U\EXT_{|{\sK \times S^1}}\right),$$
\item thus allow us 
to define equivariant propagators $\eqprop_i(\tilde{\tau})=\frac1{\delta(t)}\deleqprop_i(\tilde{\tau})$ of $\EXT$ such that $$\eqprop_i(\tilde{\tau})=\frac{1}{\delta(K)(t)}C(\cvarM_i,\sK_i,\tilde{\tau}) + \overline{\pi^{-1}(\cvarM_i)} \subset C_2(\EXT)$$
where $C(\cvarM_i,\sK_i,\tau)$ is a rational $4$--chain such that
$\partial C(\cvarM_i,\sK_i,\tau)=\partial(\cvarM_i,\sK_i,\tilde{\tau})$,
so that $$\tilde{Z}_n(K,\tilde{\tau})=I_n(\{\eqprop_i(\tau)\})$$ is an invariant of $(K,\tilde{\tau})$
and
$$\tilde{Z}(K)=\tilde{Z}(K,\tilde{\tau})\exp(-\frac{p_1(\tau_\pi,\tilde{\tau}_{\CC})}{4}\ansothree).$$
See Theorem~\ref{thmdefinvpseudo}.
\end{itemize}
These pseudo-parallelizations are genuine parallelizations outside a link tubular neighborhood in the interior of $\EXT$, and they are ``averages`` of parallelizations on this tubular neighborhood.

The reader who is only interested by the proof in the case where the $LP$-surgeries involve pairs of integral homology handlebodies -that includes the Borromean surgery case, and thus the clasper case-, or more generally in the cases where the parallelizations extend to the replacing rational homology handlebodies, does not need to read Sections~\ref{secaltdefzform} to \ref{secdefpseudoint} that justify all the above assertions, in order to get the proof of Theorem~\ref{thmflag}. She/he just needs to replace the word pseudo-parallelization by parallelization in the proof below.

In any case, the restriction of our parallelization $\tau$ of $\EXT$ to $\EXT \setminus \coprod_{i=1}^{2n}\mathring{A}^{(i)}$ extends as a pseudo-parallelization $\tau(A^{(i)\prime})$ on each $A^{(i)\prime}$. For every $S \subset N$, define the pseudo-parallelization $\tau(\EXT_S)$ so that it coincides with $\tau(\EXT)$ on $\EXT \setminus \coprod_{i\in S}\mathring{A}^{(i)\prime}$ and with $\tau(A^{(i)\prime})$ on $A^{(i)\prime}$ when $i \in S$.

The equivariant propagators $\eqprop(\EXT_S)=\frac1{\delta(t)}\deleqprop(\EXT_S)$ associated with these knot exteriors $\EXT_S$ will be associated with these pseudo-parallelizations so that
\begin{itemize}
\item $\partial \deleqprop(\EXT_S)=\partial \deleqprop(\tau(\EXT_S))$ satisfies the same conditions as before,
\item $p(i)=p_1(\tau_{\pi},\tau(\EXT_{\{i\}})_{\CC})-p_1(\tau_{\pi},\tau=\tau(\EXT))=p_1(\tau_{\pi},\tau(\EXT_{S \cup \{i\}})_{\CC})-p_1(\tau_{\pi},\tau(\EXT_S)_{\CC})$ for any $i$ of $N$ and for any $S \subset N \setminus \{i\}$, thanks to Proposition~\ref{proppont} (2), so that
$$p_1(\tau_{\pi},\tau(\EXT_S)_{\CC})=p_1(\tau_{\pi},\tau) +\sum_{i \in S}p(i).$$
\end{itemize}

\begin{lemma}
\label{lemgetridtau}
In order to prove Theorem~\ref{thmflag}, it suffices to prove that 
$$\sum_{S \subset N} (-1)^{\sharp S}\tilde{Z}_n(\EXT_S,\tau(\EXT_S))=\langle \langle \bigsqcup_{i \in N} \tilde{T}(\CI_{A^{(i)}A^{(i)\prime}}) \rangle \rangle_n.$$
(The Pontrjagin classes won't appear in the alternate sum.)
\end{lemma}
\bp
By definition, we have that
$$\sum_{S \subset N}(-1)^{\sharp S}\tilde{Z}_n(\EXT_S)=\sum_{S \subset N}(-1)^{\sharp S}\left(\tilde{Z}_n(\EXT_S,\tau(\EXT_S)) + \sum_{j<n}\tilde{Z}_j(\EXT_S,\tau(\EXT_S)) P_{n-j}(S)\right)$$
where $P_{n-j}(S)$ stands for an element of $\CA_{n-j}$ that is a combination of $m[\Gamma]$ where the $m$ are monomials in $p_1(\tau_{\pi},\tau)$ and in the $p(i)$ of degree
at most $(n-j)$, for degree $(n-j)$ Jacobi diagrams $\Gamma$. Furthermore, such an $m[\Gamma]$ appears in $P_{n-j}(S)$ if and only if $m$ is a monomial in the variables $p_1(\tau_{\pi},\tau)$ and $p(i)$ for $i \in S$.
Therefore, we can rewrite the sum of the annoying terms by factoring out the $m[\Gamma]$. Let $I \subset N$ be the subset of the $i$ such that $p(i)$ appears in $m$.
$(\sharp I \leq n-j)$. The factor of $m[\Gamma]$ reads
$$\sum_{I \subset S \subset N}(-1)^{\sharp S}\tilde{Z}_j(\EXT_S,\tau(\EXT_S))=(-1)^{\sharp I}\sum_{J \subset (N \setminus I)}(-1)^{\sharp J}\tilde{Z}_j(\EXT_{I \cup J},\tau(\EXT_{I \cup J}))$$
where $ \sharp (N \setminus I) \geq \sharp N +j - n$. Since $\sharp N +j - n=n+j$ and $j<n$, $ \sharp (N \setminus I) > 2j$. Thus since the right-hand side of the equality in the statement of the lemma is obviously invariant under a disjoint null LP-surgery, the factor of $m[\Gamma]$ vanishes if this equality holds true. 
\eop

Let $A^{(i)}_S$ denote $A^{(i)\prime}$ if $i \in S$ and $A^{(i)}$ if $i \notin S$.

For $j=1,\dots,g(A^{(i)})$, let $\Sigma^{\prime}(a^{(i)}_j)$ be a rational chain in $A^{(i)\prime}$ that intersects $ [-4,0] \times \partial A^{(i)\prime} $ as
$[-4,0]\times a^{(i)}_j $. Identify a collar $[-4,0] \times \partial A^{(i)\prime}$ of $\partial A^{(i)\prime}$ in $A^{(i)\prime}$ with 
$[-4,0] \times \partial A^{(i)}$, naturally.
Define $\delta(t)\Sigma^{S}(\check{z}^{(i)}_j) \subset \left(\widetilde{\EXT_S} \setminus p_{\EXT_S}^{-1}( \mathring{A^{(i)}_S})\right)$ from $\delta(t)\Sigma(\check{z}^{(i)}_j)$
by replacing the pieces $\Sigma(a^{(k)}_{\ell})$ by pieces $\Sigma^{S}(a^{(k)}_{\ell})$, where $\Sigma^{S}(a^{(k)}_{\ell})=\Sigma^{\prime}(a^{(k)}_{\ell})$ if $k \in S$ and $\Sigma^{S}(a^{(k)}_{\ell})=\Sigma(a^{(k)}_{\ell})$ if $k \notin S$.

Define $\check{\deleqprop}(\EXT_i=\EXT_{\{i\}})$ on $(\tilde{C}_2(\EXT_i) \setminus \mbox{Int}(\tilde{C}_2(A^{(i)\prime}_{-1})))$ so that
\begin{itemize}
 \item $\check{\deleqprop}(\EXT_i)= \deleqprop(\EXT)$ on 
$\tilde{C}_2((\EXT \setminus A^{(i)\prime})_{-4})$
\item
For any $t\in[-4,0]$, $\check{\deleqprop}(\EXT_i)$ intersects $A^{(i)\prime}_t \times_{e} (\EXT \setminus A^{(i)})_{t+3}$  as
$$ \sum_{j \in \{1, \dots, g(A^{(i)})\}} \Sigma^{\prime}_t(a^{(i)}_j)\times \delta(t)\Sigma_{t+3}(\check{z}^{(i)}_j) +\delta(\thetap)\left(A^{(i)\prime}_t \times \gamma^+_{t+3}(A^{(i)})\right).$$
\item  
For any $t\in[-4,0]$, $\check{\deleqprop}(\EXT_i)$ intersects $(\EXT\setminus A^{(i)})_{t+3} \times_{e} A^{(i)\prime}_t$ as 
$$\sum_{j \in \{1, \dots, g(A^{(i)})\}} \delta(t^{-1})\Sigma_{t+3}(\check{z}^{(i)}_j)\times \Sigma_t^{\prime}(a_j^{(i)}) - \delta(\thetap)\left(\gamma^-_{t+3}(A^{(i)}) \times A^{(i)\prime}_t\right).$$
\item $\check{\deleqprop}(\EXT_i)$ intersects a neighborhood $\ZZ \times U\EXT_i \times [0,1]$ of $\ZZ \times U\EXT_i$ in $\TCE$ as $$\delta(\thetap)\left(\{0\} \times s_{\tau}(\EXT_i;\fvarM)  \times [0,1]\right).$$
\end{itemize}

Recall that $\delta(t)\Sigma^{S}_u(\check{z}^{(i)}_j)= \delta(t)\Sigma^{S}(\check{z}^{(i)}_j) \cup -([u,0] \times \partial \delta(t)\Sigma^{S}(\check{z}^{(i)}_j)) \subset \left(\widetilde{\EXT_S} \setminus p_{\EXT_S}^{-1}(\mathring{A_u^{(i)}})\right),$ for $u \in[-4,0]$.
Set $\overline{\Sigma}^{S}_u(\check{z}^{(i)}_j)=p_{\EXT}(\delta(t)\Sigma^{S}_u(\check{z}_{j}^{(i)})) \subset \EXT_S$.
For a $2$--chain $\Sigma$ of $\EXT_S$, let $U(\Sigma)$ abusively denote $U\EXT_{S|\Sigma}$ and define $\deleqpropb^{(i)}_{S,u} \subset \partial \tilde{C}_2(\EXT_S)$ as 
$$\deleqpropb^{(i)}_{S,u}=\sum_{(j,k,\ell)\in \{1,\dots g(A^{(i)})\}^3} \CI_{A^{(i)}A^{(i)\prime}}(a^{(i)}_j, a^{(i)}_k, a^{(i)}_{\ell})\delta(\thetap)lk_e(\tilde{z}_k^{(i)},\tilde{z}^{(i)}_{\ell})(\thetap) U(\overline{\Sigma}^{S}_u(\check{z}_{j}^{(i)})). $$

Define $\check{\deleqprop}(\EXT_S)$ on $\tilde{C}_2(\EXT_S) \setminus \coprod_{i \in S}\mbox{Int}(\tilde{C}_2(A^{(i)\prime}))$ as 
\begin{itemize}
 \item $\deleqprop(\EXT)$ on $\tilde{C}_2\left(\EXT \setminus \coprod_{i \in S}\mathring{A}^{(i)}_{-4}\right)$,
\item $\sum_{j \in \{1, \dots, g(A^{(i)})\}} \Sigma_t^{\prime}(a^{(i)}_j)\times \delta(t)\Sigma^{S}_{t+3}(\check{z}^{(i)}_j) + \delta(\thetap)\left(A^{(i)\prime}_t \times \gamma^+_{A^{(i)},t+3}\right)$ on $A^{(i)\prime}_{t} \times_e (\EXT_S \setminus A^{(i)\prime})_{t+3}$, for $i\in S$ and $t\in [-4,0]$,
\item $\sum_{j \in \{1, \dots, g(A^{(i)})\}}
 \delta(t^{-1})\Sigma^{S}_{t+3}(\check{z}^{(i)}_j)\times \Sigma_t^{\prime}(a_j^{(i)}) - \delta(\thetap)\left(\gamma^-_{A^{(i)},t+3} \times A_t^{(i)\prime}\right)$ on $(\EXT_S \setminus A^{(i)\prime})_{t+3} \times_e A^{(i)\prime}_t$, for $i\in S$ and $t\in [-4,0]$.
\end{itemize}

\begin{proposition}
\label{propprodbord}
Let $\eta \in ]0,\frac{1}{6n^2}[$.
Choose an injective map $g \colon \{1,2,\dots,2n\} \rightarrow ]0,1[$.

There exist rational chains $\deleqprop(A^{(i)\prime})$ in $\tilde{C}_2(A^{(i)\prime})$ such that
\begin{itemize}
\item $\deleqprop(A^{(i)\prime})$ intersects the equivariant neighborhood $\partial \tilde{C}_2(\mathring{A}_{-1}^{(i)\prime}) \times [0,1]$ 
as the intersection of  $\partial \tilde{C}_2(\mathring{A}_{-1}^{(i)\prime}) \times [0,1]$ with
$$\left((\deleqprop(A^{(i)\prime}) \cap (\ZZ \times UA^{(i)\prime})) \times [0,1]\right) \cup 
\left(-\partial \deleqpropb^{(i)}_{\{i\},-4+g(i)\eta} \times [g(i)\eta,1] \right)$$
$$\cup \left( \deleqpropb^{(i)}_{\{i\},-4+g(i)\eta} \times \{g(i)\eta\} \right).$$
(Here, we use our equivariant neighborhood $\partial \tilde{C}_2( \mathring{\EXT_i} ) \times [0,1]$ of $\partial \tilde{C}_2( \mathring{\EXT_i} )$ in $\tilde{C}_2(\EXT_i)$ where the coordinate in $[0,1]$ is thought of as the distance between two points in a configuration near $U\EXT_i$.)
\item The chains $$\begin{array}{ll}\deleqprop(\EXT_S)=&\check{\deleqprop}(\EXT_S) \cap \left(\tilde{C}_2(\EXT_S) \setminus \coprod_{i \in S}\mathring{\tilde{C}}_2(A^{(i)\prime}) \right)+\sum_{i\in S}\deleqprop(A^{(i)\prime})\\&
+\sum_{i\in S}\left(\left(\deleqpropb^{(i)}_{S,-4+g(i)\eta} \times \{g(i)\eta\}\right)
\cap \left(\tilde{C}_2(\EXT_S) \setminus \mathring{\tilde{C}}_2(A^{(i)\prime}) \right)\right)
\end{array}$$
are of the prescribed form on $\left(\tilde{C}_2(\EXT) \setminus \tilde{C}_2(\EXT_{[0,2]})\right) \cup U\EXT_{[0,2]}$ and their boundaries are the wanted ones.

\end{itemize}
\end{proposition}
\bp See \cite[Lemma 11.11]{lesbetaone}. \eop

Then the chains $\deleqprop(\EXT_S)$ can be used to compute $\tilde{Z}_n(\EXT_S,\tau(\EXT_S))$.

\subsection{The parallel chains $\deleqprop_r(\EXT_S)$}
\label{subparchain}

The construction of the chains $\deleqprop(\EXT_S)$ involves choices of transverse chains  $\Sigma^S(a_j^{(i)})$, $\delta(t)\Sigma^{S}(\check{z}_{j}^{(i)})$ and $\gamma^{\pm}(A^{(i)})$.
In order to construct the  $\deleqprop_r(\EXT_S)$, we use parallel copies $\Sigma^S(a_j^{(i)},r)$, $\delta(t)\Sigma^{S}(\check{z}_{j}^{(i)},r)$ and $\gamma^{\pm}(A^{(i)},r)$ of these chains 
with the same properties.
We furthermore fix consistent Riemannian metrics on the $E_S$ and assume
that the $\gamma^{\pm}(A^{(i)},r)$ are at distance at least $\eta$ from each other for a small positive $\eta$ for distinct $(i,r)$.
We also assume that the neighborhoods of the intersections of three $\Sigma^S(a_j^{(i)},r)$ are identified with euclidean cubes that are intersected by the surfaces as parallel or orthogonal planes, and that the pseudo-parallelizations match the standard parallelization of these cubes.

\begin{lemma}
\label{lemmino}
The $\Sigma^S(a_j^{(i)},r)$ can be chosen so that
there exists $\nu \in ]0,1/2[$ such that for any $i$, for any $S \subset N$,
\begin{itemize}
\item for any $4$ distinct pairs $(j_k,r_k)$ in $\{1,\dots,g(A^{(i)})\} \times \{1,\dots,3n\}$, for any $x \in A^{(i)}_S$, there exists $(j_{\ell},r_{\ell})$ among the $4$ pairs such that
$d(x,\Sigma^S(a_{j_{\ell}}^{(i)},r_{\ell})) \geq \nu$,
\item if there exist $3$ distinct pairs $(j_k,r_k)$ in $\{1,\dots,g(A^{(i)})\} \times \{1,\dots,3n\}$ and $x \in A^{(i)}_S$ such that $d(x,\Sigma^S(a_{j_k}^{(i)},r_k)) < \nu$ for any $k$,
then the $j_k$ are pairwise distinct, and $x$ is in a standard neighborhood of an intersection point of the three $\Sigma^S(a_{j_k}^{(i)},r_k)$.
\end{itemize}
\end{lemma}
\bp
When a curve $a_j^{(i)}$ is null-homologous, for different $r$, the
chains $\Sigma^{S}(a_j^{(i)},r)$ (can and) will (be assumed to) be disjoint parallel surfaces. 
In general, let $k$ be the order of $a_j^{(i)}$ in $H_1(A^{(i)}_S;\ZZ)$, the $\Sigma^S(a_j^{(i)},r)$ will intersect $A^{(i)}_{S,-6}$
as disjoint parallel surfaces --with coefficient $1/k$-- 
and they will
intersect a collar $[-4,0] \times \partial A^{(i)}$ as $[-4,0] \times  \{t_r\} \times a_j^{(i)}$ where $[-1,1] \times a_j^{(i)}$ is a collar of $a_j^{(i)}$ in $\partial A^{(i)}$, and the $t_r$ are different elements of $[-1,1]$.
Then $\Sigma^S(a_j^{(i)},r)$ and $\Sigma^S(a_j^{(i)},s)$ will intersect inside $]-6,-4] \times \partial A^{(i)}$ as $k(k-1)/2$ curves $\{x_{r,s,i,j,\ell}\} \times a_j^{(i)}$, these curves will be disjoint for different $(r,s,\ell)$.
See the figure where $k=3$.

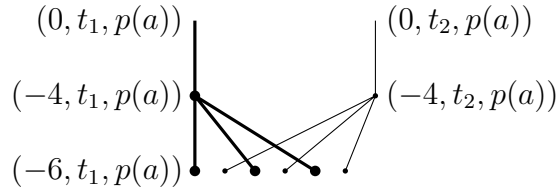
\begin{center}
\begin{tikzpicture} \useasboundingbox (0,-.7) rectangle (2.4,2.9);
\draw [very thick] (0,1) -- (1.6,0) (0.8,0) -- (0,1) node[left]{$(-4,t_1,p(a))$} (0,0) -- (0,2) node[left]{$(0,t_1,p(a))$} (0,0) node[left]{$(-6,t_1,p(a))$};
\draw (2,0) -- (2.4,1) -- (1.2,0) (0.4,0) -- (2.4,1) -- (2.4,2) (2.4,1) node[right]{$(-4,t_2,p(a))$} (2.4,2) node[right]{$(0,t_2,p(a))$};
\fill (0.4,0) circle (1pt) (1.2,0) circle (1pt) (2.4,1) circle (1pt) (2,0) circle (1pt);
\fill (0,1) circle (2pt) (0.8,0) circle (2pt) (0,0) circle (2pt) (1.6,0) circle (2pt);
\end{tikzpicture}
\end{center}

There is no loss in assuming that the intersections $\Sigma^S(a_u^{(i)},r) \cap \Sigma^S(a_v^{(i)},s)$ for a given $(u,v)$, with $u\neq v$, sit outside the collar $[-7,0] \times \partial A^{(i)}$ of $\partial A^{(i)}$, either, and we do.

For any $i$, and for any set of $4$ elements $(j_k,r_k)$ in $\left(\{1,\dots,g(A^{(i)})\} \times \{1,\dots,3n\}\right)$, the continuous function from $A^{(i)}_S$ to $\RR^+$ $$x \mapsto \mbox{max}_{k=1}^4d(x,\Sigma^S(a_{j_k}^{(i)},r_k))$$ has a non zero minimum because $A^{(i)}_S$ is compact and there are no quadruple intersections.
Similarly, for any set of $3$ elements $(j_k,r_k)$ in $\{1,\dots,g(A^{(i)})\} \times \{1,\dots,3n\}$ such that $j_1=j_2$, $$x \mapsto \mbox{max}_{k=1}^3d(x,\Sigma^S(a_{j_k}^{(i)},r_k))$$ has a non zero minimum.

The minimum over all the above finitely many minima satisfies the wanted properties.
\eop

Note the easy lemma.
\begin{lemma}
\label{lemconsC}
 Fix $3n$ distinct vectors $\cvarM_1$, \dots, $\cvarM_{3n}$ (in generic position) in $S^2$. Then there exists
$\CONS >3n$ such that for any vector $\vec{u}$ of norm less than $2/\CONS$, for any $i\neq j$, $(\cvarM_i+\vec{u})\neq \parallel \cvarM_i+\vec{u} \parallel \cvarM_j$,
and the map $(\cvarM \mapsto \frac1{\parallel \cvarM+\vec{u} \parallel}(\cvarM+\vec{u}))$ is a diffeomorphism of $S^2$.
\end{lemma}
\eop

Fix $\eta \in ]0,\frac{1}{\CONS^2}[$ such that the $\gamma^{\pm}(A^{(i)})$ are at distance at least $\eta$ from each other.

For $i \in N$, and for $r \in \{1,2,\dots,3n\}$, set $$g(i,r)=\nu\frac{i}{(\CONS^2)^r}.$$

Set $\deleqprop_{1}(\EXT_S)=\deleqprop(\EXT_S)$ with $g(i)=g(i,1)$.
Then define $\check{\deleqprop}_{r}(\EXT_S)$ from $\check{\deleqprop}(\EXT_S)$ and $\deleqprop_{r}(A^{(i)\prime})$ from $\deleqprop(A^{(i)\prime})$ by replacing $\Sigma^{S}(.)$ by $\Sigma^{S}(.,r)$, $\gamma^{\pm}(A^{(i)})$ by $\gamma^{\pm}(A^{(i)},r)$ and $g(i)$ by $g(i,r)$ everywhere in the construction.

Then set
$$\begin{array}{ll}\deleqprop_{r}(\EXT_S)=&\check{\deleqprop}_{r}(\EXT_S) \cap \left(\tilde{C}_2(\EXT_S) \setminus \coprod_{i \in I}\mathring{\tilde{C}}_2(A^{(i)\prime}) \right)\\&
+\sum_{i\in S}\left(\left(\deleqpropb^{(i)}_{S,-4+g(i,r)\eta} \times \{g(i,r)\eta\}\right)
\cap \left(\tilde{C}_2(\EXT_S) \setminus \mathring{\tilde{C}}_2(A^{(i)\prime}) \right)\right)
+\sum_{i\in S}\deleqprop_{r}(A^{(i)\prime}).\end{array}$$

Furthermore assume that for a pair of points in $U\EXT_{S[0,2]} \times \{t\} \subset U\EXT_{S[0,2]} \times [0,1]$, the distance between the two points is $t$ if $t<\nu$.

\section{Proof of Theorem~\ref{thmflag}}
\setcounter{equation}{0}
\label{secproof}

\subsection{Associating colorings with contributing intersections}

Recall $N=\{1,2,\dots,2n\}$. We want to compute 
$$\sum_{S\subset N} (-1)^{\sharp S}\tilde{Z}_n(\EXT_S,\tau(\EXT_S))=\sum_{S\subset N;\Gamma \in \CS^u_n} (-1)^{\sharp S}\frac{\CI_{\Gamma}(\{\eqprop_r(\EXT_S)\})}{2^{3n}(3n)!(2n)!}.$$
Recall that $\CS^u_n$ is defined in Subsection~\ref{subeqint}.
Consider an intersection point $m=(m_i)_{i \in N} \in \EXT_S^N$ that contributes to $\CI_{\Gamma}(\{\eqprop_r(\EXT_S)\})$ for $\Gamma \in \CS^u_n$.
Let $\eta$ be the parameter of Proposition~\ref{propprodbord}.
When $m_i$ belongs to $A^{(j)}_{-\eta}$ or $A^{(j)\prime}_{-\eta}$, color the corresponding vertex $v(i)$ of $\Gamma$ by $j$. Then color the remaining vertices by $0$.
If the intersection point is in $p(\Gamma,e(r))^{-1}\left(p_C\left(\deleqpropb^{(i)}_{S,-4+g(i,r)\eta} \times \{g(i,r)\eta\}\right)\right)$, then color the edge $e(r)$ by $i$.
Then the two vertices of $e(r)$ must have the same color unless one of them is in $A^{(j)}\setminus A^{(j)}_{-\eta}$. In this case, switch its color from $0$ to $j$. Iterate such switches so that two vertices of a colored edge always have the same color, and vertices colored by $j$ are in $A^{(j)}$ or $A^{(j)\prime}$. Let $\CaC$ denote the obtained coloring of the diagram. 
This process assigns a coloring $\CaC$ to each intersection point. Let $\CI_{\Gamma}(\{\eqprop_r(\EXT_S)\};\CaC)$ be the contribution of all the intersection points colored by $\CaC$ to $\CI_{\Gamma}(\{\eqprop_r(\EXT_S)\})$ so that $\CI_{\Gamma}(\{\eqprop_r(\EXT_S)\})=\sum_{\CaC}\CI_{\Gamma}(\{\eqprop_r(\EXT_S)\};\CaC)$. For this proof,
we denote the alternate sum of the intersections associated with $(\Gamma,\CaC)$ by $I(\Gamma,\CaC)$.
 $$I(\Gamma,\CaC)=\sum_{S\subset N} (-1)^{\sharp S}\CI_{\Gamma}(\{\eqprop_r(\EXT_S)\};\CaC).$$
$$\sum_{S\subset N} (-1)^{\sharp S}\tilde{Z}_n(\EXT_S,\tau(\EXT_S))=\sum_{(\Gamma,\CaC);\;\Gamma \in \CS^u_n,\;\CaC \;\mbox{\scriptsize coloring of}\;\Gamma}\frac{I(\Gamma,\CaC)}{2^{3n}(3n)!(2n)!}.$$
We first fix $\Gamma$ and $\CaC$ and we get rid of colorings $\CaC$ such that $I(\Gamma,\CaC)=0$, in order to be left with the {\em admissible colorings\/} of Definition~\ref{def5cases}.
Consider the graph $T(\Gamma,\CaC)$ obtained from $\Gamma$ by removing the uncolored edges.

\begin{lemma}
$T(\Gamma,\CaC)$ is a disjoint union of trees (where an isolated vertex is a tree).
\end{lemma}
\bp 
Otherwise, consider a simple cycle in $T(\Gamma,\CaC)$, and the two (points associated with the) vertices of the edge $e(r)$ with minimal label $r$.
The distance between these two vertices related by the edge $e(r)$ colored by $i$ is $g(i,r) \eta$ where
$$g(i,r) \eta > 2n g(2n,r^{\prime}) \eta$$ for any $r^{\prime}> r$. Since there are at most $(2n-1)$ edges different from $e(r)$ in that cycle, the distance between our two vertices is then less than $\frac{2n-1}{2n}{g(i,r) \eta}$, and the cycle cannot exist.
\eop

\begin{lemma}
If $I(\Gamma,\CaC)\neq 0$, then the color $0$ does not appear and every color of $N$ appears either once as an edge color or as a color of the vertices of a tree with colored edges.
\end{lemma}
\bp Let us first prove that if at least one color of $N$ does not appear in the coloring $\CaC$, then $I(\Gamma,\CaC)=0$.
Let $i$ be the smallest color that does not appear. Let $S$ be a part of $N$.
If $S$ contains $i$, then a contributing intersection point $m$ to $\CI_{\Gamma}(\{\eqprop_r(\EXT_S)\};\CaC)$ is in $(\EXT_S \setminus A^{(i)\prime}_{-\eta})^{2n}$ outside the parts $p(\Gamma,e(r))^{-1}\left(p_C\left(\deleqpropb^{(i)}_{S,-4+g(i,r)\eta} \times \{g(i,r)\eta\}\right)\right)$. 
Then the same point $s(m)$ of $(\EXT_{S \setminus \{i\}} \setminus A^{(i)}_{-\eta})^{2n}$ will contribute to 
$\CI_{\Gamma}(\{\eqprop_r(\EXT_{S \setminus \{i\}})\};\CaC)$ like $m$
contributes to $\CI_{\Gamma}(\{\eqprop_r(\EXT_S)\};\CaC)$, so that the contributions of $m$ and $s(m)$ cancel in the alternate sum.
Setting $s(s(m))=m$ allows us to define an involution $s$ of the set of intersection points such that one color of $N$ does not appear, so that the contribution of $m$ and $s(m)$ cancel each other. Thus if $I(\Gamma,\CaC)\neq 0$, then all the colors of $N$ appear in $\CaC$.

There is one vertex color in each component of $T(\Gamma,\CaC)$. Therefore the number of colors in $\Gamma$ is at most the sum of the number of components of $T(\Gamma,\CaC)$ and the number of edges of $T(\Gamma,\CaC)$, that is the number $2n$ of vertices of $\Gamma$,
because of the previous sublemma.
Therefore, the number of colors in $\Gamma$ is at most $2n$.
Finally, since the $2n$ colors of $N$ must appear, the color $0$ cannot appear and every color of $N$ appears either once as an edge color or as a color of the vertices of a tree.
\eop

Fix $\nu$ and the $\Sigma^S(a_j^{(i)},r)$ so that the conditions of Lemma~\ref{lemmino} are satisfied.

\begin{lemma}
\label{lemtreethree}
Let $T$ be a component of $\;T(\Gamma,\CaC)$.
Let $E(T)$ denote the set of edges of \/$T$, let \/$V(T)$ be the set of vertices of \/$T$, and let $E_e(T)$ be the set of edges of \/$\Gamma$ that have one end in $V(T)$. If $I(\Gamma,\CaC)\neq 0$, then
$$\sharp(E(T))+\sharp(E_e(T))\leq 3.$$
\end{lemma}

\bp
The vertices of $V(T)$ are colored by $\ell \in N$.
An edge of $E_e(T)$ constrains its vertex in $T$ to live in some $\Sigma(a_{i(r)}^{(\ell)},r)$.
The edges of $T$ constrain both their vertices to live at a distance less than $g(i,r)\eta$ from some $\Sigma(a_{i(r)}^{(\ell)},r)$.
All the vertices of $T$ are at distance less than $(2n-1)g(2n,1)\eta$ from each other where $(2n-1)g(2n,1)\eta< \nu\eta < \frac{\nu}{2}$. 
Then, since $g(i,r) \eta < \frac{\nu}{2}$, a vertex of $T$ will be at distance less than $\nu$ from all the $\Sigma(a_{i(r)}^{(\ell)},r)$ associated with the edges of $E(T)$ or $E_e(T)$ that are distinct because the $r$ are distinct.
According to Lemma~\ref{lemmino}, this finishes the proof. \eop

\begin{definition}
\label{def5cases}
Let $\CS_n^{\ell}$ be the set of trivalent graphs with $2n$ numbered vertices and with $3n$ numbered and oriented edges, with possible loops.
 Let $\Gamma \in \CS_n^{\ell}$. An {\em admissible coloring of\/} $\Gamma$ is a coloring $\CaC$ of all the vertices of $\Gamma$ and some edges of $\Gamma$ by elements of $N$ such that
\begin{itemize}
 \item all the elements of $N$ color either an edge or at least one vertex,
\item two vertices have the same color if and only if they are connected by a path of colored edges,
\item the graph $T(\Gamma,\CaC)$ obtained from $\Gamma$ by removing the uncolored edges is made of components $T$ that satisfy one of the following set of hypotheses.
\begin{itemize}
\item{Case V:} $T$ is a vertex.
\item{Case $\theta$:} $T$ is an edge, and the component of $T$ in $\Gamma$ is a $\theta$-graph.
\item{Case D1:} $T$ is an edge, and there is one uncolored edge in $\Gamma$ that has the same ends as $T$.
\item{Case D2:} $T$ is made of two edges $e(r_1)$ and $e(r_2)$ sharing one vertex $v(b)$, the other vertex $v(a)$ of $e(r_1)$ belongs to the only edge of $E_e(T)$, the third vertex $v(c)$ of $T$ is connected both to $v(a)$ and $v(b)$ by uncolored edges, and $r_1<r_2$.
\item{Case D3:} $T$ is a linear tree with $4$ vertices, and its component in $\Gamma$ is as in the figure below with $r_1<r_2$ and $r_1<r_3$.
\end{itemize}
\end{itemize}

These $5$ cases are drawn below, the colored edges are dotted.\\
\begin{center}
Case V: \begin{tikzpicture} \useasboundingbox (-.3,-.1) rectangle (.9,.4);
\begin{scope}[yshift=-.1cm]
\draw (0,0) -- (.3,.2) -- (.3,.5) (.6,0) -- (.3,.2);
\fill (.3,.2) circle (1.5pt);
\end{scope}
\end{tikzpicture}\\
\medskip
Case $\theta$: 
\begin{tikzpicture}\useasboundingbox (-.3,0) rectangle (.9,.4);
\begin{scope}[yshift=-.1cm]
\draw (0.3,0.2) circle (.3);
\draw [very thick, dotted] (0,.2) -- (.6,.2) (.55,.25) node[right]{\scriptsize $\ell$} (.05,.25) node[left]{\scriptsize $\ell$} (.3,.1) node[above]{\scriptsize $i$};
\fill (0,.2) circle (1.5pt) (.6,.2) circle (1.5pt);\end{scope}
\end{tikzpicture}\\
Case D1: 
\begin{tikzpicture}\useasboundingbox (-.3,-.6) rectangle (2,.5);
\begin{scope}[yshift=-.5cm]
\draw (-.25,0) -- (.25,0) (1.25,0) -- (1.75,0) (.25,0) .. controls (.25,-.4) and (1.25,-.4) .. (1.25,0) (1.15,.15) node[right]{\scriptsize $\ell$} (.35,.15) node[left]{\scriptsize $\ell$} (.75,.2) node[above]{\scriptsize $i$};
\draw [very thick, dotted] (.25,0) .. controls (.25,.4) and (1.25,.4) .. (1.25,0);
\fill (.25,0) circle (1.5pt) (1.25,0) circle (1.5pt);\end{scope}
\end{tikzpicture}\\
Case D2: \begin{tikzpicture}\useasboundingbox (-.5,-.2) rectangle (3.4,1);
\begin{scope}[yshift=-.1cm]
\draw (-.4,0) -- (0,0) (0,0) .. controls (0,.6) and (3,.6) .. (3,0) (1.5,0) .. controls (1.5,.4) and (3,.4) .. (3,0) (2.9,.15) node[right]{\scriptsize $\ell$} (1.6,.15) node[left]{\scriptsize $\ell$} (2.25,-.1) node[above]{\scriptsize $i$} (.75,-.1) node[above]{\scriptsize $j$} (2.25,.1) node[below]{\scriptsize $e(r_2)$} (.75,.1) node[below]{\scriptsize $e(r_1)$};
\draw [very thick, dotted] (0,0) -- (3,0);
\fill (3,0) circle (1.5pt) (0,0) circle (1.5pt) (1.5,0) circle (1.5pt);\end{scope}
\end{tikzpicture}
where $r_1 < r_2$\\
Case D3: \begin{tikzpicture}\useasboundingbox (-.4,-.6) rectangle (3.4,.5);
\begin{scope}[yshift=-.5cm]
\draw (0,0) .. controls (0,-.6) and (3,-.6) .. (3,0) (0,0) .. controls (0,.4) and (1,.4) .. (1,0) (2.9,.15) (2,0) .. controls (2,.4) and (3,.4) .. (3,0) (2.9,.15) node[right]{\scriptsize $\ell$} (2.1,.15) node[left]{\scriptsize $\ell$} (.9,.15) node[right]{\scriptsize $\ell$} (.1,.15) node[left]{\scriptsize $\ell$} (.5,-.1) node[above]{\scriptsize $i$} (1.5,-.1) node[above]{\scriptsize $j$} (2.5,-.1) node[above]{\scriptsize $k$} (.5,.1) node[below]{\scriptsize $e(r_2)$} (1.5,.1) node[below]{\scriptsize $e(r_1)$} (2.5,.1) node[below]{\scriptsize $e(r_3)$};
\draw [very thick, dotted] (0,0) -- (3,0);
\fill (3,0) circle (1.5pt) (0,0) circle (1.5pt) (1,0) circle (1.5pt) (2,0) circle (1.5pt);\end{scope}
\end{tikzpicture}
where $r_1 < r_2$ and $r_1<r_3$.
\end{center}
\end{definition}

\begin{lemma}
\label{lem5cases}
Let $\Gamma \in\CS_n^u$. If $I(\Gamma,\CaC)\neq 0$, then $\CaC$ is an admissible coloring of $\Gamma$.
\end{lemma}
\bp If $\sharp(E(T))=0$, then $T$ is a vertex. If $\sharp(E(T))=1$, then $T$ is an edge, and $\sharp(E_e(T))\leq 2$, then either $\sharp(E_e(T))=0$, and the component of $T$ in $\Gamma$ is a $\theta$-graph, or $\sharp(E_e(T))=2$, and there is one uncolored edge in $\Gamma$ that has the same ends as $T$.

If $\sharp(E(T))=2$, then $T$ is made of two consecutive edges sharing a vertex $v(b)$, and $\sharp(E_e(T))$ that is necessarily odd is one.
\begin{center}
\begin{tikzpicture}\useasboundingbox (-.4,-.1) rectangle (2.4,.1);
\begin{scope}[yshift=-.1cm]
\draw (.05,0) node[left]{\scriptsize $v(a)$} (1.1,.15) node[left]{\scriptsize $v(b)$} (1.95,0) node[right]{\scriptsize $v(c)$};
\draw [very thick, dotted] (0,0) -- (2,0);
\fill (1,0) circle (1.5pt) (0,0) circle (1.5pt) (2,0) circle (1.5pt);\end{scope}
\end{tikzpicture}
\end{center}

In particular, in all cases but Case $\theta$, $\sharp(E(T))+\sharp(E_e(T))= 3$. Therefore the vertices of $T$ must be at distance less than $\nu$ from three $\Sigma(a_{i(r)}^{(\ell)},r)$ associated with distinct $r$. Thus they are near triple intersection points of surfaces, according to Lemma~\ref{lemmino}, and at a distance less that $\frac{\nu}2$ from each other (see the proof of Lemma~\ref{lemtreethree}).

Recall that the pseudo-parallelizations are standard on the neighborhoods of these points.
Then the constraint associated with an uncolored edge $e(r)$ of $\Gamma$ from a vertex $v(x)$ of $T$ to a vertex $v(y)$ of $T$ will sit inside $U(A_{-7}^{(\ell)}) \times [0,1]$ and since our chains $\deleqprop$ read as products by $[0,1]$ there, according to Proposition~\ref{propprodbord}, (where this $[0,1]$ is thought of as an interval of distances between two points), such an uncolored edge $e(r)$ constrains the direction of the vector from $v(x)$ to $v(y)$ to be positively colinear to $\cvarM_r$.

Back to the case where $\sharp(E(T))=2$, there cannot be a double uncolored edge between the two ends of our linear tree with two edges because the $\cvarM_r$ are distinct.

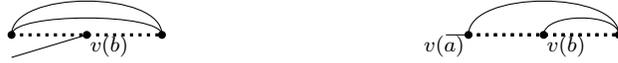
\begin{center}
\begin{tikzpicture}\useasboundingbox (-.4,-.4) rectangle (2.4,.4);
\begin{scope}[yshift=-.1cm]
\draw (0,-.3) -- (1,0) (.9,-.15) node[right]{\scriptsize $v(b)$} (0,0) .. controls (0,.3) and (2,.3) .. (2,0) (0,0) .. controls (0,.6) and (2,.6) .. (2,0);
\draw [very thick, dotted] (0,0) -- (2,0);
\fill (1,0) circle (1.5pt) (0,0) circle (1.5pt) (2,0) circle (1.5pt);\end{scope}
\end{tikzpicture} \hspace{3cm} \begin{tikzpicture}\useasboundingbox (-.4,-.4) rectangle (2.4,.4);
\begin{scope}[yshift=-.1cm]
\draw (-.3,0) -- (0,0) (.1,-.15) node[left]{\scriptsize $v(a)$} (.9,-.15) node[right]{\scriptsize $v(b)$} (1,0) .. controls (1,.3) and (2,.3) .. (2,0) (0,0) .. controls (0,.6) and (2,.6) .. (2,0);
\draw [very thick, dotted] (0,0) -- (2,0);
\draw [very thick, dotted] (0,0) -- (2,0);
\fill (1,0) circle (1.5pt) (0,0) circle (1.5pt) (2,0) circle (1.5pt);\end{scope}
\end{tikzpicture} 
\end{center}

Therefore the edge of $E_e(T)$ meets $T$ at a vertex $v(a) \neq v(b)$, and the edges of $\Gamma$ involving $T$ are as in Case D2. Now, if $r_1>r_2$, then the length of the edge $e(r_2)$ (or of a geodesic segment joining the images of its ends) is greater than $\CONS$ times the length of $e(r_1)$, then according to our assumption on $\CONS$ in Lemma~\ref{lemconsC}, the directions of the two uncolored edges are too close to each other to be colinear to different $\cvarM_r$.

We are left with the case when $\sharp(E(T))=3$. In this case $\sharp(E_e(T))=0$, and $T$ contains the $4$ vertices of a connected component of $\Gamma$.
Then either $T$ has a trivalent vertex and we say that $T$ is a \emph{tripod}, or $T$ is a linear tree.
If $T$ is a tripod, \begin{tikzpicture}\useasboundingbox (-.2,0) rectangle (.8,.3);
\begin{scope}[yshift=-.1cm]
\draw [very thick, dotted] (.6,0) -- (.3,.2) (0,0) -- (.3,.2) -- (.3,.5);
\fill (.3,.5) circle (1.5pt) (0,0) circle (1.5pt) (.3,.2) circle (1.5pt) (.6,0) circle (1.5pt);\end{scope}
\end{tikzpicture}, the component of $T$ in $\Gamma$ is $\tetra$. One of the edges of $E(T)$ is far longer than the two other ones : the length of the edge indexed by the smallest $r$ is greater than $\CONS$ times the length of any of the two other edges. Then the directions of the two edges of $\tetra \setminus T$ that are connected to the exterior vertex of this one are too close to each other in $S^2$ to be mapped to different $\cvarM_i$, because of the condition on $\CONS$ and the $\cvarM_i$ in Lemma~\ref{lemconsC}.

Then $T$ is a linear tree with $4$ vertices 
\begin{tikzpicture}\useasboundingbox (-.7,-.1)(3.7,.4);
\draw (.05,0) node[left]{\scriptsize $v(a)$} (1.1,.15) node[left]{\scriptsize $v(b)$} (1.95,.15) node[right]{\scriptsize $v(c)$} (2.95,0) node[right]{\scriptsize $v(d)$};
\draw [thick, dotted] (0,0) -- (3,0);
\fill (1,0) circle (1.5pt) (0,0) circle (1.5pt) (2,0) circle (1.5pt) (3,0) circle (1.5pt);
\end{tikzpicture}
and the extra edges that connect them in $\Gamma$ are 
\begin{itemize}
 \item either $\{v(b),v(c)\}$ and twice
$\{v(a),v(d)\}$, \begin{tikzpicture}\useasboundingbox (-.4,-.1)(3.4,.4);
\draw 
(0,0) .. controls (0,.6) and (3,.6) .. (3,0) (0,0) .. controls (0,-.4) and (3,-.4) .. (3,0) (1,0) .. controls (1,.3) and (2,.3) .. (2,0);
\draw [very thick, dotted] (0,0) -- (3,0);
\fill (1,0) circle (1.5pt) (0,0) circle (1.5pt) (2,0) circle (1.5pt) (3,0) circle (1.5pt);
\end{tikzpicture}
\item or $\{v(b),v(d)\}$, $\{v(a),v(c)\}$ and $\{v(a),v(d)\}$, \begin{tikzpicture}\useasboundingbox (-.4,-.1)(3.4,.4);
\draw 
(0,0) .. controls (0,.6) and (3,.6) .. (3,0) (0,0) .. controls (0,-.4) and (2,-.4) .. (2,0) (1,0) .. controls (1,.3) and (3,.3) .. (3,0);
\draw [very thick, dotted] (0,0) -- (3,0);
\fill (1,0) circle (1.5pt) (0,0) circle (1.5pt) (2,0) circle (1.5pt) (3,0) circle (1.5pt);
\end{tikzpicture}
\item  or in the last case
$\{v(b),v(a)\}$, $\{v(d),v(c)\}$ and $\{v(a),v(d)\}$
\begin{tikzpicture}\useasboundingbox (-.4,-.1)(3.4,.4);
\draw 
(0,0) .. controls (0,.4) and (1,.4) .. (1,0) (3,0) .. controls (3,.4) and (2,.4) .. (2,0) (0,0) .. controls (0,-.4) and (3,-.4) .. (3,0);
\draw [very thick, dotted] (0,0) -- (3,0);
\fill (1,0) circle (1.5pt) (0,0) circle (1.5pt) (2,0) circle (1.5pt) (3,0) circle (1.5pt);
\end{tikzpicture}
\end{itemize}
The first two cases can be ruled out by the longest edge argument, and the third case can be ruled out by the same argument unless $\{v(b),v(c)\}$ is the longest edge of $T$. Thus we have no contributions unless we are in Case $D3$.
\eop

\subsection{Sketch of the proof of Theorem~\ref{thmflag}}

The definition of $I(\Gamma,\CaC)$ can be extended to any graph of $\CS_n^{\ell}$ equipped with an admissible coloring as the following alternate sum of $\CA_n^h$ intersections.
$$I(\Gamma,\CaC)=\sum_{S \subset N}(-1)^{\sharp S}I_{\Gamma}(S;\CaC)$$
where $I_{\Gamma}(S;\CaC)$ is the intersection $I_{\Gamma}$ defined by the following constraints when $\Gamma$ has no $\theta$-case.
\begin{itemize}
 \item The vertices colored by $j$ must be in $A_S^{(j)}$ (that is  $A^{(j)\prime}$ if $ j \in S$ and $A^{(j)}$ if $ j \notin S$).
\item The constraint associated with an uncolored edge $e(r)$ from a vertex colored by $i$ to a vertex colored by $j$ is 
$$\eqprop_{(i,j)}=\sum_{(u,v)\in \CE(i,j)} lk_e(\tilde{z}_{u}^{(i)},\tilde{z}_{v}^{(j)})(t)\Sigma^S(a_{u}^{(i)},r) \times
\Sigma^S (a_{v}^{(j)},r)$$
where $\CE(i,j)=\{1,\dots,g(A^{(i)})\} \times \{1,\dots,g(A^{(j)})\}$ if $j \neq i$ and
$\CE(i,i)=\{1,\dots,g(A^{(i)})\}^2 \setminus \mbox{diag}$.
\item The constraint associated with an edge $e(r)$ colored by $i \in S$ between two vertices colored by $m$ is 
$$\sum \CI_{A^{(i)}A^{(i)\prime}}(a^{(i)}_j, a^{(i)}_k, a^{(i)}_{\ell})lk_e(\tilde{z}_k^{(i)},\tilde{z}^{(i)}_{\ell})(t)lk(z_j^{(i)},z_v^{(m)})U_{g(i,r)\eta}(\Sigma^S(a_v^{(m)},r)) $$
where the sum runs over the $((j,k,\ell),v) \in \{1,\dots, g(A^{(i)})\}^3\times \{1,\dots, g(A^{(m)})\}$ and\\ $U_{g(i,r)\eta}(\Sigma^S(a_v^{(m)},r))$ denotes $\{0\} \times  U\EXT_{S|\Sigma^S(a_v^{(m)},r)} \times \{g(i,r)\eta\}$ in the neighborhood $\{0\} \times U\EXT_S \times [0,1]$ of $\{0\} \times U\EXT_S$ in $\tilde{C}_2(M)$.
\item The constraint associated with an uncolored edge $e(r)$ between two vertices colored by $\ell$ reads
$$s_{\tau_S}(A_S^{(\ell)};\cvarM_r) \times [0,1]$$
and determines the direction of the vector associated with its ends.
\end{itemize}

It is easy to check that this definition of $I(\Gamma,\CaC)$ is consistent with the existing one when $\Gamma \in \CS^u_n$.

For a graph $\Gamma$ of $\CS^u_n$ equipped with a coloring $\CaC$, $I(\Gamma,\CaC)$ is the product of the intersections $I(\Gamma_i;\CaC_{|\Gamma_i})$, where the product runs over the connected components of $\Gamma$. Similarly, for a graph $\Gamma$ of $\CS^{\ell}_n$ equipped with an admissible coloring $\CaC$, $I(\Gamma,\CaC)$ is the product of the intersections $I(\Gamma_i;\CaC_{|\Gamma_i})$, where the product runs over the connected components of $\Gamma$. This finishes the general definition of $I(\Gamma,\CaC)$, since $I(\Gamma,\CaC)$ has already been defined when $\Gamma$ is a $\theta$-case.

A {\em simple coloring\/} of $\Gamma \in \CS_n^{\ell}$ is an admissible coloring $\CaC$ such that no edge is colored.

The following subsections will be devoted to the proof of the following proposition.

\begin{proposition}
\label{propkeythmflag}
 $$\sum_{\Gamma \in \CS^u_n}\sum_{\CaC\; \mbox{\scriptsize \textrm admissible coloring of}\;\Gamma}I(\Gamma,\CaC)=\sum_{\Gamma \in \CS^{\ell}_n}\sum_{\CaC\; \mbox{\scriptsize \textrm simple coloring of}\;\Gamma}I(\Gamma,\CaC)$$
\end{proposition}

Assuming Proposition~\ref{propkeythmflag}, since
$$\sum_{S\subset N} (-1)^{\sharp S}\tilde{Z}_n(\EXT_S, \tau(\EXT_S))=\sum_{(\Gamma,\CaC);\;\Gamma \in \CS^u_n,\;\CaC \;\mbox{\scriptsize admissible coloring of}\;\Gamma}\frac{I(\Gamma,\CaC)}{2^{3n}(3n)!(2n)!},$$
Theorem~\ref{thmflag} is the direct consequence of the following proposition that we prove in this subsection.

\begin{proposition}
\label{propfinthmflag}
 $$\sum_{\Gamma \in \CS^{\ell}_n}\sum_{\CaC\; \mbox{\scriptsize \textrm simple coloring of}\;\Gamma}\frac{I(\Gamma,\CaC)}{2^{3n}(3n)!(2n)!}=\langle \langle \bigsqcup_{i \in N} \tilde{T}(\CI_{A^{(i)}A^{(i)\prime}}) \rangle \rangle_n$$
\end{proposition}

\begin{definition}
\label{defcdadmissible} 
For $i\in N$, let $\CD(i)$ be the set of maps $d(i,.) \colon  \{1,2,3\} \rightarrow \NN$ such that $$1\leq d(i,1) <d(i,2)<d(i,3) \leq g(A^{(i)}).$$ For $d(i,.) \in \CD(i)$, set $$\CI(d(i,.))=\CI_{A^{(i)}A^{(i)\prime}}(a_{d(i,1)},a_{d(i,2)},a_{d(i,3)}).$$
Let ${\cal D}$ be the set of maps $d \colon N \times \{1,2,3\} \rightarrow \NN$ such that
for all $i \in N$, $d(i,.) \in \CD(i)$. For $d\in \CD$, set
$$\CI(d)=\prod_{i \in N}\CI(d(i,.)).$$

Note the obvious lemma:

\begin{lemma}
\label{lemdecd}
For $d \in \CD$, set
$$\tilde{T}(i,d(i,.))= \begin{tikzpicture}\useasboundingbox (-.55,-.3) rectangle (3,.9);
\begin{scope}[yshift=-.25cm]
\draw (0.7,0) -- (1.5,0) (1.1,0) -- (1.1,.4) (1.5,0) node[right]{$\tilde{z}^{(i)}_{d(i,1)}$} (.7,0) node[left]{$\tilde{z}^{(i)}_{d(i,3)}$} (1.1,.4) node[above]{$\tilde{z}^{(i)}_{d(i,2)}$};
\fill (1.1,0) circle (1.5pt) (1.2,.1) node[below]{\scriptsize $i$};
\fill (1.1,0) circle (1.5pt);\end{scope}
\end{tikzpicture}\;\;\mbox{and} \;\;
\tilde{T}(d)=\bigsqcup_{i \in N}\tilde{T}(i,d(i,.)).$$

Then $$\langle \langle \bigsqcup_{i \in N} \tilde{T}(\CI_{A^{(i)}A^{(i)\prime}}) \rangle \rangle_n=
\sum_{d \in {\cal D}}\CI(d)\langle \langle \tilde{T}(d) \rangle \rangle_n.$$
\end{lemma}
\eop

A simple coloring of $\Gamma \in \CS_n^{\ell}$ induces a bijection from the set of vertices of $\Gamma$ to $N$.
For a simple coloring $\CaC$ of a graph $\Gamma$, and for $d \in \CD$, define $I(\Gamma,\CaC,d)$ as follows.
Denote the three half-edges of the vertex colored by $i$ of $\Gamma$ by $h(i,1)$, $h(i,2)$ and $h(i,3)$, respectively, arbitrarily.
The bijection $h$ from $N\times \{1,2,3\}$ to the set $H(\Gamma)$ of half-edges of $\Gamma$ induces an orientation $o(h)$ of $\Gamma$ that is called its \emph{$h$-vertex-orientation.}
An edge $e(r)$ in our graph $\Gamma$ goes from a vertex colored by $i(r)$ to a vertex colored by $j(r)$.
There is a natural bijection $b$ from $\{1,2,\dots,3n\} \times \{1,2\}$ to $H(\Gamma)$ that maps $(r,1)$ to the first half-edge of $e(r)$ and $(r,2)$ to the second one.
Let ${\cal G}(\Gamma,\CaC,d)$ be the set of maps $\gamma\colon \{1,2,\dots,3n\} \times \{1,2\}\rightarrow \NN$ such that
$\gamma \circ b^{-1}(h(i,k))=d(i,\rho_{\gamma}(i,k))$ for a permutation $\rho_{\gamma}(i,.)$ of $\{1,2,3\}$, for all $i$. Let $\mbox{sign}(\rho_{\gamma}(i,.))=\pm 1$ be the signature of $\rho_{\gamma}(i,.)$, and set $\varepsilon(\gamma,h)=\prod_{i \in N}\mbox{sign}(\rho_{\gamma}(i,.))$. Note that $\rho_{\gamma}$ is determined by $\gamma$ and that $\gamma$ is determined by $\rho_{\gamma}$.

Equip $\Gamma$ with its $h$-vertex-orientation. For $\gamma \in {\cal G}(\Gamma,\CaC,d)$, define $[\Gamma_{\gamma}(o(h))]$ as the (class in $\CA_n^h$) of the graph obtained by beading the oriented edge $e(r)$ of $\Gamma(o(h))$ by $lk_e(\tilde{z}_{\gamma(r,1)}^{(i(r))},\tilde{z}_{\gamma(r,2)}^{(j(r))})(t)$.
Define $$I(\Gamma,\CaC,d)=\sum_{\gamma \in {\cal G}(\Gamma,\CaC,d)}\varepsilon(\gamma,h)\CI(d)[\Gamma_{\gamma}(o(h))].$$
\end{definition}

\begin{proposition}
\label{propthmflagstone}
If a graph $\Gamma \in \CS_n^{\ell}$ is equipped with a simple coloring $\CaC$, then 
$$I(\Gamma,\CaC)=\sum_{d\in \CD}I(\Gamma,\CaC,d).$$
\end{proposition}
\bp
For an edge $e(r)$ from a vertex $v(x)$ colored by $i(r)$ to a vertex
$v(y)$ colored by $j(r)$, the constraint is a sum running over the elements $(t,u)$ of $\CE(i(r),j(r))$ of constraints that read
$$(m_x,m_y) \in  \Sigma^S(a_{t}^{(i(r))},r) \times
\Sigma^S(a_{u}^{(j(r))},r)$$ with a beading coefficient $lk_e(\tilde{z}_{t}^{(i(r))},\tilde{z}_{u}^{(j(r))})(t)$,
where
the normal bundle to $\Sigma^S(a_{t}^{(i)},r) \times
\Sigma^S(a_{u}^{(j)},r)$ is oriented by the normal $V(r,1)$ to $\Sigma^S(a_{t}^{(i)},r)$ in Copy $x$ of $\EXT_S$ followed by the normal $V(r,2)$ to $\Sigma^S(a_{u}^{(j)},r)$ in Copy $y$ of $\EXT_S$. Associate $V(r,1)$ with the half-edge
$(v(x),e(r))$ and $V(r,2)$ with the half-edge
$(v(y),e(r))$.

\begin{lemma}
 \label{lemorientation} Let $\Gamma \in \CS_n^u$.
Consider a point $m=(m_1,\dots,m_{2n})$ in a transverse intersection $\cap_{r \in \{1,2,\dots,3n\}} p(\Gamma,r)^{-1}(p_C(\deleqprop_r))$. Assume that the oriented normal bundle to $\deleqprop_r$ is generated by two vectors whose preimages $V(r,1)$ and $V(r,2)$ live in $T_{m_{k(V(r,1))}}\EXT$ and $T_{m_{k(V(r,2))}}\EXT$, respectively, so that exactly three independent normals live in every $T_{m_{j}}\EXT$, for every $j$, ($k(V(r,1))=x$ and $k(V(r,2))=y$, above). Let $\varepsilon(m)[\Gamma_{\rho(m)}]$ denote the contribution of $m$ to $I_{\Gamma}(\{\eqprop_r\})$ (where $\Gamma$ is equipped with its canonical vertex-orientation $o(\Gamma)$). Then the sign of the contribution $\varepsilon(m)[\Gamma_{\rho(m)}]$ may be computed as follows. 

Assign the normal $V(r,1)$ (resp. $V(r,2)$) to the first half-edge of $e(r)$ (resp. to the second one). Then use a permutation $\sigma$ of the set of normals to permute the normals so that the normals in $T_{m_{j}}\EXT$ are assigned to the three half-edges adjacent to $v(j)$. Equip $\Gamma$ with an arbitrary vertex-orientation $o_a$ and let $[\Gamma(o_a)]$ be the corresponding element of $\CA_n^h$.
Let $\varepsilon(m,a)$ be the product of the sign of the permutation of the normals by all the signs associated with the vertices, where the sign associated with the vertex $v(j)$ is $1$ if the normals assigned to the half-edges of $v(j)$ in the cyclic order induced by $o_a$ induce the orientation of $T_{m_{j}}\EXT$, and it is $(-1)$ otherwise.
Then $\varepsilon(m)[\Gamma_{\rho(m)}]=\varepsilon(m,o_a)[\Gamma_{\rho(m)}(o_a)]$.
\end{lemma}
Observe that the sign of $\varepsilon(m,o_a)[\Gamma_{\rho(m)}(o_a)]$ is independent of $o_a$ and of the numbering of vertices.

\noindent{\sc Proof of Lemma~\ref{lemorientation}:}
Equip $\Gamma$ with its canonical vertex-orientation $o(\Gamma)$ (induced by the edge-orientation and the order of factors, see Subsection~\ref{subeqint}). Then the sign in front of the beaded $\Gamma$ associated with the intersection point $m$ is determined by comparing the orientation of $\EXT^{2n}$ with the orientation given by the sequence of normals associated with the half-edges ordered by the order of factors of $\EXT^{2n}$ and by the vertex-orientation (up to even permutations).
\eop

In our current proof of Proposition~\ref{propthmflagstone}, there is no need of permuting the half-edges and the contribution of the normals to the sign in front of $[\Gamma]$, where $\Gamma$ is now equipped with the arbitrary vertex-orientation $o(h)$ that does not necessarily coincide with the canonical one, is the product over the vertices of the signs associated with each oriented vertex by comparing the orientation induced by the three normals to the orientation of $\EXT_S$.
The original numbering of the vertices of $\Gamma$ is unimportant, but the vertices inherit a new numbering from the coloring of the vertices.

Numbering the vertices of $\Gamma$ with the coloring, our intersection point $m$ reads $(m_1,m_2,\dots, m_{2n})$ where
each $m_i$ is an intersection point of three surfaces $\Sigma(a_{d(i,k)}^{(i)},r(i,k))$ in $A^{(i)}$ if $i \notin S$ and it is an intersection point of three surfaces $\Sigma^{\prime}(a_{d(i,k)}^{(i)},r(i,k))$ in $A^{(i)\prime}$ if $i \in S$. In this latter case, the contribution
is multiplied by $(-1)$ in the alternate sum. Then we can get rid of the $S$ and of the alternate sum by considering $m_i$ as an intersection point of
three surfaces $\Sigma(a_{d(i,k)}^{(i)},r(i,k)) \cup (-\Sigma^{\prime}(a_{d(i,k)}^{(i)},r(i,k)))$ in $A^{(i)} \cup -A^{(i)\prime}$.

Here, we fix the $d(i,k)$, for $k=1,2,3,$ so that $1\leq d(i,1) <d(i,2)<d(i,3) \leq g(A^{(i)})$.
Thus our intersection point determines an element $d$ of ${\cal D}$ and a permutation $\rho_{\gamma}(i,.)$ of $\{1,2,3\}$ such that the normal to $\Sigma(a_{d(i,\rho_{\gamma}(i,k))}^{(i)},r(h(i,k)))$ is associated with $h(i,k)$, and the sum of the contributions to $I(\Gamma,\CaC)$ of the intersection points
that determine $d$ and $\gamma$ is 
$\varepsilon(\gamma,h)\CI(d)[\Gamma_{\gamma}(o(h))].$ This concludes the proof of Proposition~\ref{propthmflagstone}.
\eop

In particular, for a graph $\Gamma$ of $\CS^{\ell}_n$ equipped with a simple coloring $\CaC$, $I(\Gamma,\CaC)$ is independent of the numbering of the edges, of their orientations and of the initial numbering of the vertices.

Thanks to Lemma~\ref{lemdecd} and Proposition~\ref{propthmflagstone}, Proposition~\ref{propfinthmflag} is now the direct consequence of the following lemma.
\begin{lemma}
For any $d \in \CD$,
$$\CI(d)\langle \langle \tilde{T}(d) \rangle \rangle_n=\sum_{\Gamma \in \CS^{\ell}_n}\sum_{\CaC\; \mbox{\scriptsize \textrm simple coloring of}\;\Gamma}\frac{I(\Gamma,\CaC,d)}{2^{3n}(3n)!(2n)!}.$$
\end{lemma}
\bp 
For any graph $\Gamma$ of $\CS_n^{\ell}$, the simple coloring induced by the numbering of vertices of $\Gamma$ in $\CS_n^{\ell}$ is called the {\em canonical coloring\/} of $\Gamma$, it is denoted by $\CaC_{\tiny can}(\Gamma)$. Note that
$$\sum_{\Gamma \in \CS^{\ell}_n}\sum_{\CaC \; \mbox{\scriptsize \textrm simple coloring of}\;\Gamma}\frac{I(\Gamma,\CaC,d)}{2^{3n}(3n)!(2n)!}
=\sum_{\Gamma \in \CS^{\ell}_n}\frac{I(\Gamma,\CaC_{\tiny can}(\Gamma),d)}{2^{3n}(3n)!}.$$

Consider the set $\CB$ of bijections from $\{1,2,\dots,3n\} \times \{1,2\}$ to $N\times \{1,2,3\}$. Let $p_N \colon N\times \{1,2,3\} \rightarrow N$ be the natural projection. With each bijection $b \in \CB$, associate the maps $i\colon \{1,2,\dots,3n\} \rightarrow N$ and $j\colon \{1,2,\dots,3n\} \rightarrow N$ such that $i(r)=p_N(b(r,1))$ and $j(r)=p_N(b(r,2))$.
Then the bijection $b$ determines the graph $\Gamma(b) \in \CS_n^{\ell}$ such that the edge $e(r)$ of $\Gamma(b)$ goes from $v(i(r))$ to $v(j(r))$.
Furthermore, $b$ equips the half-edge $b^{-1}(i,k)$ of the vertex $v(i)$ of $\Gamma(b)$ 
with $d(i,k)$, it also equips the vertex $v(i)$ with the orientation induced by $b^{-1}(i,.)$. This vertex-orientation of $\Gamma(b)$ will be called its $d$-orientation and denoted by $o(d)$. Let $\Gamma(b)_{b}(o(d))$ be the graph obtained from $\Gamma(b)$ by beading each edge $e(r)$ by $lk_e(\tilde{z}^{(i(r))}_{d(b(r,1))},\tilde{z}^{(j(r))}_{d(b(r,2))})$.
 
Observe that for any $\Gamma \in \CS_n^{\ell}$, $$I(\Gamma,\CaC_{\tiny can}(\Gamma),d)=\sum_{b \in \CB; \Gamma(b)=\Gamma}\CI(d)[\Gamma(b)_{b}(o(d))].$$
Thus
$$\sum_{\Gamma \in \CS^{\ell}_n}\sum_{\CaC \; \mbox{\scriptsize \textrm simple coloring of}\;\Gamma}\frac{I(\Gamma,\CaC,d)}{2^{3n}(3n)!(2n)!}
=\sum_{b \in \CB}\frac{\CI(d)[\Gamma(b)_{b}(o(d))]}{2^{3n}(3n)!}.$$

Now, a bijection $b$ of $\CB$ determines a pairing of $P\left(\tilde{T}(d) \right)$. See the definitions before Theorem~\ref{thmflag}.  Namely, if $b(r,1)=(i,k)$ and $b(r,2)=(j,\ell)$, the leg labeled by $d(i,k)$ is paired with the leg labeled by $d(j,\ell)$. Two bijections $b$ and $b^{\prime}$ determine the same pairing if and only if 
they differ by a permutation of $\{1,2,\dots,3n\} \times \{1,2\}$ that preserves the partition in pairs. Such a permutation permutes the labels of the edges and reverses some edge orientations. There are $2^{3n}(3n)!$ such permutations.
\eop

\begin{remark}
 Let $\CS_n^{vc}$ denote the set of trivalent graphs (with possible loops) with $2n$ numbered vertices, whose edges are unordered and unoriented.
Let $\Gamma \in \CS^{\ell}_n$ be a graph equipped with a simple coloring $\CaC$.
Let $\Gamma(\CaC) \in\CS_n^{vc} $ be the graph obtained from $(\Gamma,\CaC)$ by forgetting the numberings and the orientations of $\Gamma$ and by numbering the vertices of $\Gamma$ by the coloring.
For a graph $\Gamma \in \CS^{\ell}_n$, equipped with a simple coloring $\CaC$,
$I(\Gamma,\CaC,d)$ depends neither on the numbering of edges and vertices of $\Gamma \in \CS^{\ell}_n$, nor on the orientations of the edges (because of the conjugation relation for the edge reversals), thus it only depends on $\Gamma(\CaC)$ and will be denoted by $I(\Gamma(\CaC),d)$.
There are $\frac{(3n)!2^{3n}(2n)!}{\sharp \mbox{\small Aut}(\Gamma(\CaC))}$ colored graphs $(\Gamma,\CaC)$ that give rise to a given $\Gamma(\CaC)$, where $\mbox{\small Aut}(\Gamma(\CaC))$ denotes the set of automorphisms of $\Gamma(\CaC)$. Such an automorphism must preserve
the colored vertices pointwise. It only permutes multiple edges with identical ends and reverses looped edges.
$$\sum_{\Gamma \in \CS^{\ell}_n}\sum_{\CaC \;\mbox{\scriptsize \textrm simple coloring of}\;\Gamma}\frac{I(\Gamma,\CaC,d)}{2^{3n}(3n)!(2n)!}
=\sum_{\Gamma_c \in \CS^{vc}_n}\frac{I(\Gamma(\CaC),d)}{\mbox{\small Aut}(\Gamma(\CaC))}.$$
\end{remark}

\subsection{The $\theta$ case}

First note the following two easy lemmas.

\begin{lemma}
\label{lembmbbar}

 In $\CA^h$, for any beading $b(t)$ we have

$$\begin{tikzpicture}\useasboundingbox (-.5,-.6) rectangle (2,.6);
\begin{scope}[yshift=-.5cm]
\draw (-.25,0) -- (.25,0) (1.25,0) -- (1.75,0) (.25,0) .. controls (.25,-.4) and (1.25,-.4) .. (1.25,0) (.75,.2) node[above]{\scriptsize $b(t)-\overline{b}(t)$} (1.25,0) .. controls (1.25,.15) and (.9,.3) .. (.75,.3);
\draw [->] (.25,0) .. controls (.25,.15) and (.6,.3) .. (.75,.3);
\fill (.25,0) circle (1.5pt) (1.25,0) circle (1.5pt);\end{scope}
\end{tikzpicture}=\begin{tikzpicture}\useasboundingbox (-.5,-.6) rectangle (2,.6);
\begin{scope}[yshift=-.5cm]
\draw (-.25,0) -- (1.75,0) (.75,0) -- (.75,.4) arc (-90:270:.2) (.75,.7) node[above]{\scriptsize $b(t)$};
\draw [-<] (.75,.4) arc (-90:90:.2);
\fill (.75,0) circle (1.5pt) (.75,.4) circle (1.5pt);\end{scope}
\end{tikzpicture} = - \begin{tikzpicture}\useasboundingbox (-.5,-.6) rectangle (2,.6);
\begin{scope}[yshift=-.5cm]
\draw (-.25,0) -- (1.75,0) (.75,0) -- (.75,.4) arc (-90:270:.2) (.75,.7) node[above]{\scriptsize $\overline{b}(t)$};
\draw [-<] (.75,.4) arc (-90:90:.2);
\fill (.75,0) circle (1.5pt) (.75,.4) circle (1.5pt);\end{scope}
\end{tikzpicture}
=\frac12 \begin{tikzpicture}\useasboundingbox (-.5,-.6) rectangle (2,.6);
\begin{scope}[yshift=-.5cm]
\draw (-.25,0) -- (1.75,0) (.75,0) -- (.75,.4) arc (-90:270:.2) (.75,.7) node[above]{\scriptsize $b(t)-\overline{b}(t)$};
\draw [-<] (.75,.4) arc (-90:90:.2);
\fill (.75,0) circle (1.5pt) (.75,.4) circle (1.5pt);\end{scope}
\end{tikzpicture}$$
where $\overline{b}(t)=b(t^{-1})$.
\end{lemma}
\bp These are easy consequences of the Jacobi relation, the conjugation relation and the antisymmetry relation.
\eop

\begin{lemma}
\label{lembeadvan}
For any $P(t) \in \frac1{\delta(t)}\QQ[t,t^{-1}]$,

$$\begin{tikzpicture}\useasboundingbox (-.2,-.6) rectangle (2,.6);
\begin{scope}[yshift=-.5cm]
\draw (-.25,0) -- (1.75,0) (.75,0) -- (.75,.7) arc (-90:270:.2) (.9,.95) node[right]{\scriptsize $b(t)$} (.7,.3) node[right]{\scriptsize $P(t)$};
\draw [->] (.75,.7) arc (-90:0:.2);
\draw [->] (.75,0) -- (.75,.4);
\fill (.75,0) circle (1.5pt) (.75,.7) circle (1.5pt);\end{scope}
\end{tikzpicture} =\begin{tikzpicture}\useasboundingbox (-.2,-.6) rectangle (2,.6);
\begin{scope}[yshift=-.5cm]
\draw (-.25,0) -- (1.75,0) (.75,0) -- (.75,.7) arc (-90:270:.2) (.9,.95) node[right]{\scriptsize $b(t)$} (.7,.3) node[right]{\scriptsize $P(1)$};
\draw [->] (.75,.7) arc (-90:0:.2);
\draw [->] (.75,0) -- (.75,.4);
\fill (.75,0) circle (1.5pt) (.75,.7) circle (1.5pt);\end{scope}
\end{tikzpicture}$$

\end{lemma}
\eop

Let $\Xi(i,\ell)$ be the following component of a colored graph of $\CS^{\ell}_n$
$$ \Xi(i,\ell)=\begin{tikzpicture}\useasboundingbox (-.2,-.4) rectangle (1.4,.2);
\begin{scope}[yshift=-.45cm]
\draw (.3,.2) -- (.9,.2) (.1,.2) circle (.2) (1.1,.2) circle (.2) (.2,.35) node[right]{\scriptsize $i$} (1,.35) node[left]{\scriptsize $\ell$};
\fill (.3,.2) circle (1.5pt) (.9,.2) circle (1.5pt);\end{scope}
\end{tikzpicture}$$
where the coloring is the only represented data. See Definition~\ref{def5cases}.

Let $\mathfrak{A}_3$ denote the set of even permutations of $\{1,2,3\}$. For $i \in N$, $\sigma \in \mathfrak{A}_3$, and $d(i,.) \in \CD(i)$, set
$$B(i,\sigma,d(i,.))=lk_e(\tilde{z}^{(i)}_{d(i,\sigma(2))},\tilde{z}^{(i)}_{d(i,\sigma(3))})(t)-lk_e(\tilde{z}^{(i)}_{d(i,\sigma(3))},\tilde{z}^{(i)}_{d(i,\sigma(2))})(t).$$

Note the two immediate lemmas.
\begin{lemma}
\label{lemIxiell}
$$I(\Xi(i,\ell))=\sum_{(d(i,.),d(\ell,.)) \in \CD(i) \times \CD(\ell)}\sum_{(\sigma,\tau) \in \mathfrak{A}_3^2}\CI(d(i,.))\CI(d(\ell,.))lk(z_{d(i,\sigma(1))}^{(i)},z_{d(\ell,\tau(1))}^{(\ell)})\begin{tikzpicture}\useasboundingbox (0,-.3) rectangle (2.7,.8);
\begin{scope}[yshift=-.75cm]
\draw (.2,.3) -- (.2,.9) arc (-90:270:.2) (.2,-.1) arc (-90:270:.2) (.35,1.1) node[right]{\scriptsize $B(\ell,\tau,d(\ell,.))(t)$} (.35,.1) node[right]{\scriptsize $B(i,\sigma,d(i,.))(t)$};
\draw [->] (.2,.9) arc (-90:0:.2);
\draw [->] (.2,-.1) arc (-90:0:.2);
\fill (.2,.3) circle (1.5pt) (.2,.9) circle (1.5pt);\end{scope}
\end{tikzpicture} .$$
\end{lemma}

\begin{lemma}
\label{lemdeccoledge}
$$\deleqpropb^{(i)}_{S,u} \times \{g(i,r)\eta\}=\sum_{d(i,.) \in \CD(i)}\sum_{\sigma \in \mathfrak{A}_3} \CI(d(i,.))(\delta B(i,\sigma,d(i,.)))(\thetap) U_{g(i,r)\eta}(\overline{\Sigma}^{S}_u(\check{z}_{d(i,\sigma(1))}^{(i)}))$$
where $U_{g(i,r)\eta}(\overline{\Sigma})$ denotes $\{0\} \times U\EXT_{S|\overline{\Sigma}} \times \{g(i,r)\eta\}$ in the neighborhood $\{0\} \times U\EXT_S \times [0,1]$ of $\{0\} \times U\EXT_S$ in $\tilde{C}_2(M)$.
\end{lemma}
\eop

\begin{proposition}
\label{propthetaiellone}
$$ I\left(\begin{tikzpicture}\useasboundingbox (-.6,-.6) rectangle (2.4,.2);
\begin{scope}[yshift=-.5cm]
\draw (-.1,.1) node[left]{\scriptsize $\ell$} (1.9,.1) node[right]{\scriptsize $\ell$} (.9,-.1) node[above]{\scriptsize $i$} (.9,.1) node[below]{\scriptsize $e(r)$} (.9,.6) node{\scriptsize $e(t)$} (.9,-.6) node{\scriptsize $e(s)$} (-.2,0) .. controls (-.2,.6) and (2,.6) .. (2,0) (-.2,0) .. controls (-.2,-.6) and (2,-.6) .. (2,0);
\draw [very thick, dotted] (-.2,0) -- (2,0);
\fill (-.2,0) circle (1.5pt) (2,0) circle (1.5pt);\end{scope}
\end{tikzpicture}  \right)
=\left\{\begin{array}{ll} 0 & \mbox{if} \;r=\mbox{max}(r,s,t),\\
       I(\Xi(i,\ell)) \; & \mbox{if}\; r=\mbox{min}(r,s,t),\\    
\frac12 I(\Xi(i,\ell)) \; & \mbox{otherwise.}
\end{array}\right.$$

(In particular, it is independent of the edge-orientation.)
\end{proposition}
\bp
The vertices $v(x)$ and $v(y)$ of an intersection configuration are in $A^{(\ell)}$ or $A^{(\ell)\prime}$, and, according to Lemma~\ref{lemdeccoledge}, their pair is in $U_{g(i,r)\eta}(\Sigma^S(a_{v}^{(\ell)},r))$ equipped with the coefficient $$\CI(d(i,.))B(i,\sigma,d(i,.))(\thetap)lk(z_{d(i,\sigma(1))}^{(i)},z_v^{(\ell)}),$$ 
for some $d(i,.) \in \CD(i)$, for some $\sigma \in \mathfrak{A}_3$, and for some $v \in \{1,\dots,g(A^{(\ell)})\}$.
This pair of vertices is in $U(A_S^{(\ell)}) \times [0,1]$ and in two other chains $p_C(\deleqprop_t(\EXT_S))$ and $p_C(\deleqprop_s(\EXT_S))$. There,
according to Proposition~\ref{propprodbord}, up to parts $\left(\deleqpropb^{(j)}_{S,-4+g(j,t)\eta} \times \{g(j,t)\eta\}\right)$
that cannot intersect because $g(j,t) \neq g(i,r)$,
$\eqprop_t(\EXT_S)$ reads $\left(s_{\tau_S}(A^{(\ell)};\cvarM_t) \times [0,1]\right)$
if $\ell \notin S$ and $\deleqprop_t(\EXT_S)$ reads $$\delta(\thetap)\left(s_{\tau_S}(A^{(\ell)\prime};\cvarM_t) \times [0,1]\right) \cup $$
$$\left(-\partial \deleqpropb^{(\ell)}_{S,-4+g(\ell,t)\eta}  \times [g(\ell,t)\eta,1] \right)
\cup \left(\deleqpropb^{(\ell)}_{S,-4+g(\ell,t)\eta} \times \{g(\ell,t)\eta\}\right)$$
if $\ell \in S$.
If we are dealing with actual parallelizations $s_{\tau_S}(A^{(\ell)};\cvarM_t) \cap s_{\tau_S}(A^{(\ell)};\cvarM_s)=\emptyset$, and
since $g(i,r)\neq g(\ell,s)$, up to exchanging $t$ and $s$, our pair is in $p_{C}( \deleqpropb^{(\ell)}_{S,-4+g(\ell,s)\eta} \times [g(\ell,s)\eta,1])$.
If we are dealing with pseudo-parallelizations, since the algebraic intersection $$\langle s_{\tau_S}(\Sigma^S(a_{v}^{(\ell)},r);\cvarM_t), s_{\tau_S}(\Sigma^S(a_{v}^{(\ell)},r);\cvarM_s)\rangle_{U(\Sigma^S(a_{v}^{(\ell)},r))}$$ vanishes, according to \cite[Lemma 10.5]{lesbetaone}, the possible intersections in
$s_{\tau_S}(A^{(\ell)\prime};\cvarM_t) \cap s_{\tau_S}(A^{(\ell)\prime};\cvarM_s)$ do not contribute to $\tilde{Z}_n(\EXT_S)$ so that we may assume
in any case that our pair is in 
$$ p_{C}(\partial \deleqpropb^{(\ell)}_{S,-4+g(\ell,s)\eta} \times [g(\ell,s)\eta,1])\cap (s_{\tau_S}(\Sigma^S(a_{v}^{(\ell)},r);\cvarM_t) \times \{g(i,r)\eta\})$$
and therefore that $g(i,r)>g(\ell,s)$, $r<s$, and $\{i,\ell\}\subset S$.
Recall that $$\deleqpropb^{(\ell)}_{S,u}=\sum_{d(\ell,.)\in \CD(\ell)}
\sum_{\tau \in \mathfrak{A}_3}\CI(d(\ell,.))(\delta B(\ell,\tau,d(\ell,.)))(\thetap)U(\overline{\Sigma}_{u}^S(\check{z}_{d(\ell,\tau(1))}^{(\ell)})).$$
Now, the third constraining propagator $\eqprop_t(\EXT_S)$ reads $s_{\tau_S}(A^{(\ell)\prime};\cvarM_t) \times [0,1]$ near our intersection point and $\tau_S$ is a genuine parallelization there, so that this third propagator constrains the direction of the vector from $m_x$ to $m_y$, and the vector $(m_y-m_x)$ is fixed because its length is fixed by $e(r)$.

Assume that the three edges go from $v(x)$ to $v(y)$. 
Then Edge $e(r)$ constrains $m_x$ to live in $\Sigma^S(a_{v}^{(\ell)},r)$, and Edge $e(s)$ constrains $m_x$ to live in $\partial\overline{\Sigma}_{u}^S(\check{z}_{d(\ell,\tau(1))}^{(\ell)})$, so that $v=d(\ell,\tau(1))$.

Associate the normal $z=z_{v}^{(\ell)}$ to $\Sigma^S(a_{v}^{(\ell)},r)$ with the half-edge $(v(x),e(r))$ and the positive normal $N(S)=N(S_x)$ to the sphere $S_x$ around $v(x)$ with radius $g(i,r)\eta$ with $(v(y),e(r))$.
Note that $(-\partial U(\overline{\Sigma}_{u}^S(\check{z}_{d(\ell,\tau(1))}^{(i)})) \times [g(\ell,s)\eta,1])$ is cooriented by the positive normal $a=a_{d(\ell,\tau(1))}^{(\ell)}$ to $\overline{\Sigma}_{u}^S(\check{z})$ followed by the outward normal $N(\partial A^{(\ell)})$ to $\partial A^{(\ell)}_{u}$.
These normals are respectively associated with the half-edges $(v(x),e(s))$ and $(v(y),e(s))$.
Associate two tangent vectors $v_1$ and $v_2$ to $S^2$ at $\cvarM_t$ that orient $S^2$ with the half-edges $(v(x),e(t))$ and $(v(y),e(t))$.
Using Lemma~\ref{lemorientation}, we compute the signs with the following picture:

$$
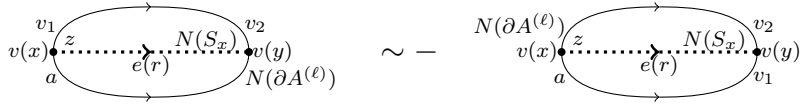
\begin{tikzpicture}\useasboundingbox (-1,-.5) rectangle (3.8,.5);
\begin{scope}[yshift=-.4cm]
\draw [->] (-.2,0) .. controls (-.2,.6) and (.8,.6) .. (1.1,.6);
\draw [->] (-.2,0) .. controls (-.2,-.6) and (.8,-.6) .. (1.1,-.6);
\draw [very thick, dotted,->] (-.2,0) -- (1.1,0);
\draw [very thick, dotted] (1.1,0) -- (2.4,0);
\draw (1.1,.6) .. controls (1.4,.6) and (2.4,.6) .. (2.4,0) (1.1,-.6) .. controls (1.4,-.6) and (2.4,-.6) .. (2.4,0) (1.1,.1) node[below]{\scriptsize $e(r)$} (-.2,.15) node[right]{\scriptsize $z$} (2.4,.15) node[left]{\scriptsize $N(S_x)$} (2.3,0) node[right]{\scriptsize $v(y)$} (2.2,.35) node[right]{\scriptsize $v_2$} (2.2,-.35) node[right]{\scriptsize $N(\partial A^{(\ell)})$} (-.1,0) node[left]{\scriptsize $v(x)$} (0,.35) node[left]{\scriptsize $v_1$} (0,-.35) node[left]{\scriptsize $a$};
\fill (-.2,0) circle (1.5pt) (2.4,0) circle (1.5pt);\end{scope}
\end{tikzpicture}\;\;\sim - \;\;
\begin{tikzpicture}\useasboundingbox (-1.6,-.5) rectangle (3.2,.5);
\begin{scope}[yshift=-.4cm]
\draw [->] (-.2,0) .. controls (-.2,.6) and (.8,.6) .. (1.1,.6);
\draw [->] (-.2,0) .. controls (-.2,-.6) and (.8,-.6) .. (1.1,-.6);
\draw [very thick, dotted,->] (-.2,0) -- (1.1,0);
\draw [very thick, dotted] (1.1,0) -- (2.4,0);
\draw (1.1,.6) .. controls (1.4,.6) and (2.4,.6) .. (2.4,0) (1.1,-.6) .. controls (1.4,-.6) and (2.4,-.6) .. (2.4,0) (1.1,.1) node[below]{\scriptsize $e(r)$} (-.2,.15) node[right]{\scriptsize $z$} (2.4,.15) node[left]{\scriptsize $N(S_x)$} (2.3,0) node[right]{\scriptsize $v(y)$} (2.2,.35) node[right]{\scriptsize $v_2$} (2.2,-.35) node[right]{\scriptsize $v_1$} (-.1,0) node[left]{\scriptsize $v(x)$} (0,.35) node[left]{\scriptsize $N(\partial A^{(\ell)})$} (0,-.35) node[left]{\scriptsize $a$};
\fill (-.2,0) circle (1.5pt) (2.4,0) circle (1.5pt);\end{scope}
\end{tikzpicture}
$$
\medskip

\noindent so that we get positive signs at the two vertices of the right-hand side vertex-oriented $\theta$.

Fix $d(i,.)$ and $d(\ell,.)$. For $(\sigma,\tau) \in \mathfrak{A}_3^2$, set $B(i,\sigma)=B(i,\sigma,d(i,.))(t)$ and $B(\ell,\tau)=B(\ell,\tau,d(\ell,.))(t)$.
Then the intersection associated with this graph with its extra-structure given by $d(i,.)$ and $d(\ell,.)$, the coloring and the choice of $s>r$ is
$$\sum_{(\sigma,\tau) \in \mathfrak{A}_3^2}\CI(d(i,.))\CI(d(\ell,.))lk(z_{d(i,\sigma(1))}^{(i)},z^{(\ell)}_{d(\ell,\tau(1))})(-1)
\thetabibell$$
where there is no possible choice of $s$ if $r=\mbox{max}(r,s,t)$, there are two choices if $r=\mbox{min}(r,s,t)$, and there is one otherwise.

Since $\overline{B(i,\sigma)}=-B(i,\sigma)$, Lemma~\ref{lembmbbar} implies that 
$$\thetabibell \;\; = -\frac{1}2\;\;
\begin{tikzpicture}\useasboundingbox (-1.1,-.4) rectangle (2.3,.4);
\begin{scope}[yshift=-.5cm]
\draw [->] (.9,.2) arc (-180:0:.2);
\draw [->] (.3,.2) arc (0:180:.2);
\draw (.3,.2) -- (.9,.2) (1.3,.2) arc (0:180:.2) (-.1,.2) arc (-180:0:.2) (1.25,.2) node[right]{\scriptsize $B(i,\sigma)$} (-.05,.2) node[left]{\scriptsize $B(\ell,\tau)$};
\fill (.3,.2) circle (1.5pt) (.9,.2) circle (1.5pt);\end{scope}
\end{tikzpicture} \;\; \mbox{in} \; \CA^h.$$

Thus Proposition~\ref{propthetaiellone} is proved when the three edges are oriented from $v(x)$ to $v(y)$.

\begin{lemma}
\label{lemoredgedir}
Reversing the orientation of an edge with fixed length associated with a constraint $s_{\tau}(\EXT;\cvarM_t) \times [0,1]$ does not change the contribution
of the intersection point.
\end{lemma}
\bp The constraint on Edge $e(t)$ from $v(x)$ to $v(y)$ is induced by the $S^2$--valued map ''direction of $m_y-m_x$``. Assume that $(\cvarM_t,v_1,v_2)$ is a direct orthogonal basis of $\RR^3$. Then $m_y-m_x=\lambda\cvarM_t+ \alpha v_1+ \beta v_2$, and the constraint reads $(\alpha,\beta)=0$.

When this constraint is seen as a constraint on $m_y$ for a given $m_x$, let $v_{1,y}$ (resp. $v_{2,y}$) be the vector
$v_1$ (resp. $v_2$) seen in $T_{m_y}S_x$ as a vector tangent to the sphere $S_x$ centered at $m_x$ with radius the length $\lambda$ of $(m_y-m_x)$. Then the normal vectors associated with $e(t)$ are $v_{1,y}$ assigned to the first half-edge $b(t,1)=b(t,x)$ and $v_{2,y}$ assigned to $b(t,2)=b(t,y)$.

When the constraint is seen as a constraint on $m_x= m_y - \lambda\cvarM_t+ \alpha (-v_1)+ \beta (-v_2)$ for a given $m_y$, let $(-v_{1,x})$ (resp. $(-v_{2,x})$) be the vector
$(-v_1)$ (resp. $(-v_2)$) seen in $T_{m_x}S_y$ as a vector tangent to the sphere $S_y$ centered at $m_y$ with radius $\lambda$. Then the normals associated with $e(t)$ are $(-v_{1,x})$ assigned to the first half-edge $b(t,1)$ and $(-v_{2,x})$ assigned to $b(t,2)$.

Thus when the orientation of $e(t)$ is reversed, $e(t)$ goes from $v(y)$ to $v(x)$, the associated constraint $(\alpha,\beta)=0$ with $m_x-m_y=\lambda\cvarM_t+ \alpha v_1+ \beta v_2$, seen as a constraint on $m_y$ will be represented by the normals $(-v_{1,y}) \in T_{m_y}S_x$ on the new initial half-edge
$b(t,y)$ that contains $v(y)$ and $(-v_{2,y}) \in T_{m_y}S_x$ on $b(t,x)$.

In a next step, the normal vectors will be permuted so that the normal vectors in $T_{m_y}S_x$ are assigned to the half-edges that contain $v(y)$.
Next, the orientation of $T_{m_y}S_x$ induced by the normal vectors around $v(y)$ and by the vertex-orientation of $v(y)$ will be compared with the orientation of $E_S$. The involved normal vectors are $\pm v_{1,y}$, $\pm v_{2,y}$ and the normal $N(S_x)$ associated with the fixed distance between $m_x$ and $m_y$,
that is a positive multiple of $\cvarM_t$ for the first orientation of $e(t)$ and a negative one, otherwise.
\eop

Let us now study the effect of reversing the orientation of the edge $e(s)$. 
This exchanges the normals $a$ and $N^-(\partial A^{(\ell)})$ and thus reverses the sign, and the bead $B(\ell,\tau)$ is now on the edge
with the opposite orientation, so that if we keep the initial orientation, the bead is changed to $\overline{B(\ell,\tau)}=-B(\ell,\tau)$.

Finally, the contribution does not depend on the edge orientations since changing all the orientations changes the graph by permuting its vertices.
This concludes the proof of Proposition~\ref{propthetaiellone}. 
\eop

\subsection{Proof of Proposition~\ref{propkeythmflag} }

We now take care of the colored edges outside the $\theta$-components by proving the following proposition.

\begin{proposition}
\label{propcruc}
 Let $\Gamma$ be a graph of $\CS^{\ell}_n$ equipped with an admissible coloring $\CaC$. If there exists an edge $e(r)$ colored by $i$ from a vertex $v(x)$ to a vertex $v(y)$ colored by $\ell$ outside a $\theta$-component, and an uncolored edge $e(s)$ between $v(x)$ and $v(y)$, define $(\phi_i(\Gamma),\phi_i(\CaC))$ as the colored graph obtained from $(\Gamma,\CaC)$ as follows:
The vertex $v_{\phi}(y)$ of $\phi_i(\Gamma)$ is colored by $i$, and the edge $e_{\phi}(r)$ of $\phi_i(\Gamma)$ is looped at $v_{\phi}(y)$.
The edge $e_{\phi}(s)$ of $\phi_i(\Gamma)$ goes from $v_{\phi}(x)$ to $v_{\phi}(y)$ if $e(s)$ goes from $v(x)$ to $v(y)$, and it goes from 
$v_{\phi}(y)$ to $v_{\phi}(x)$, otherwise. The only change in the other parts of the graph is that $v_{\phi}(x)$ replaces $v(y)$ in the other edge that contains $v(y)$ in $\Gamma$. This can be pictured as follows:
$$
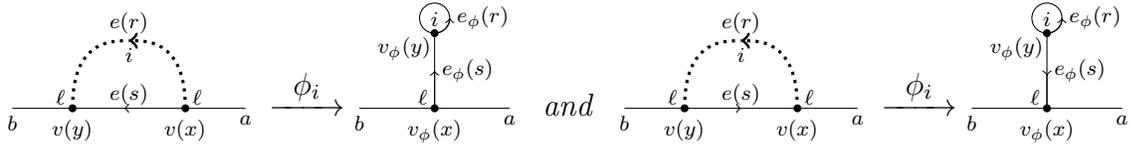
\begin{tikzpicture}\useasboundingbox (0,-.1) rectangle (3.5,.8);
\draw [-<] (.2,0) -- (1.75,0);
\draw (1.7,0) -- (3.3,0) (.2,.05) node[below]{\scriptsize $b$} (3.3,.05) node[below]{\scriptsize $a$} (1,0) node[below]{\scriptsize $v(y)$} (2.5,0) node[below]{\scriptsize $v(x)$} (1.75,-.1) node[above]{\scriptsize $e(s)$}  (1.75,.85) node[above]{\scriptsize $e(r)$} (1.75,.95) node[below]{\scriptsize $i$} (1.05,.15) node[left]{\scriptsize $\ell$} (2.45,.15) node[right]{\scriptsize $\ell$};
\draw [very thick, dotted, ->] (1,0) .. controls (1,.9) and (1.5,.9) .. (1.75,.9)(2.5,0) .. controls (2.5,.9) and (2,.9) .. (1.75,.9);
\fill (1,0) circle (1.5pt) (2.5,0) circle (1.5pt);
\end{tikzpicture}\; \hfl{\phi_i}
\begin{tikzpicture}\useasboundingbox (-.2,-.1) rectangle (2.2,1.5);
\draw (0,0) -- (2,0) (1,0) -- (1,1) arc (-90:270:.2) (1.15,1.2) node[right]{\scriptsize $e_{\phi}(r)$} (.95,.5) node[right]{\scriptsize $e_{\phi}(s)$} (0,.05) node[below]{\scriptsize $b$} (2,.05) node[below]{\scriptsize $a$} (1,0) node[below]{\scriptsize $v_{\phi}(x)$} (1,.95) node[above]{\scriptsize $i$} (1.05,.15) node[left]{\scriptsize $\ell$} (1.1,.8) node[left]{\scriptsize $v_{\phi}(y)$};
\draw [->] (1,1) arc (-90:0:.2);
\draw [->] (1,0) -- (1,.5);
\fill (1,0) circle (1.5pt) (1,1) circle (1.5pt);
\end{tikzpicture}
\;\; \mbox{and}\;\;\begin{tikzpicture}\useasboundingbox (0,-.1) rectangle (3.5,.8);
\draw [->] (.2,0) -- (1.75,0);
\draw (1.7,0) -- (3.3,0) (.2,.05) node[below]{\scriptsize $b$} (3.3,.05) node[below]{\scriptsize $a$} (1,0) node[below]{\scriptsize $v(y)$} (2.5,0) node[below]{\scriptsize $v(x)$} (1.75,-.1) node[above]{\scriptsize $e(s)$}  (1.75,.85) node[above]{\scriptsize $e(r)$} (1.75,.95) node[below]{\scriptsize $i$} (1.05,.15) node[left]{\scriptsize $\ell$} (2.45,.15) node[right]{\scriptsize $\ell$};
\draw [very thick, dotted, ->] (1,0) .. controls (1,.9) and (1.5,.9) .. (1.75,.9)(2.5,0) .. controls (2.5,.9) and (2,.9) .. (1.75,.9);
\fill (1,0) circle (1.5pt) (2.5,0) circle (1.5pt);
\end{tikzpicture}\; \hfl{\phi_i} 
\begin{tikzpicture}\useasboundingbox (-.2,-.1) rectangle (2.2,1.5);
\draw (0,0) -- (2,0) (1,0) -- (1,1) arc (-90:270:.2) (1.15,1.2) node[right]{\scriptsize $e_{\phi}(r)$} (.95,.5) node[right]{\scriptsize $e_{\phi}(s)$} (0,.05) node[below]{\scriptsize $b$} (2,.05) node[below]{\scriptsize $a$} (1,0) node[below]{\scriptsize $v_{\phi}(x)$} (1,.95) node[above]{\scriptsize $i$} (1.05,.15) node[left]{\scriptsize $\ell$} (1.1,.8) node[left]{\scriptsize $v_{\phi}(y)$};
\draw [->] (1,1) arc (-90:0:.2);
\draw [-<] (1,0) -- (1,.5);
\fill (1,0) circle (1.5pt) (1,1) circle (1.5pt);
\end{tikzpicture}$$

Then $$I(\Gamma,\CaC)=\frac12 I(\phi_i(\Gamma);\phi_i(\CaC)).$$
\end{proposition}

The edge-orientations do not affect the $I(.,.)$, but they will be used later in combinatorial arguments.

\begin{lemma}
\label{lemtinyconf}
Under the assumptions of Proposition~\ref{propcruc}, the vector $(m_y-m_x)$ is a tiny vector determined by the two edges, and the number and the nature of the intersection points corresponding to $(\Gamma,\CaC)$ are independent of this tiny subconfiguration.
Such a $(v(x),v(y),e(r),e(s))$ will be called an {\em independent tiny subconfiguration.\/}
\end{lemma}
\bp Our two vertices $v(x)$ and $v(y)$ are colored by $\ell$. Since there is no $\theta$-case, the point $m_x$ is very close to a triple intersection of surfaces $\Sigma^S(a_{d(\ell,k)}^{(\ell)})$. See the proof of Lemma~\ref{lemtreethree}. The colored edge determines the length of $(m_y-m_x)$ that is the very small $g(i,r)\eta$, so that our pair is in $U(A^{(\ell)}) \times [0,1]$. Near a triple intersection of surfaces, the contributing part of $\eqprop_s(\EXT_S)$ reads $\left(s_{\tau_S}(A^{(\ell)};\cvarM_s) \times [0,1]\right)$ and determines the direction of the edge, so that $(m_y-m_x)$ is fixed as soon as we have a double-edge in our cases D1, D2 and D3. See Definition~\ref{def5cases}.
Then in case D1, imposing $m_x \in \Sigma(a)$ or $(m_y \in \Sigma(a) \Leftrightarrow m_x \in \Sigma(a)+(m_x-m_y))$ will not
change the nature of the solution $m_x$ that will be at the intersection of three surfaces $\Sigma(a^{(\ell)}_{d(\ell,k)})$ or at the intersection of three parallel surfaces close to these surfaces.
In cases D2 and D3, the vectors corresponding to the edges colored by $i$ and $k$ are similarly determined.
The length of the remaining colored edge is also determined so that the corresponding vector $\cvarM$ should be found in a sphere, with a constraint on the direction of $\cvarM + \vec{u}$ for a fixed vector $\vec{u}$ whose length is sufficiently small so that 
$(\cvarM \mapsto \frac1{\parallel\cvarM + \vec{u}\parallel}(\cvarM + \vec{u}))$ is a diffeomorphism from our sphere to the standard sphere, thanks to Lemma~\ref{lemconsC} and to the definition of $g(.,.)$ that follows it.
In all these cases, the small drawn subconfiguration is a rigid set of vectors determined in a unique way by the constraints.
In addition, $m_x$ must be at the intersection of three surfaces close and parallel to some $\Sigma^S(a^{(\ell)}_{d(\ell,k)})$.
\eop

\begin{lemma}
\label{lemconttiny}
The contribution of an independent tiny subconfiguration
$$
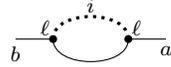
\begin{tikzpicture}\useasboundingbox (-.3,-.6) rectangle (2,.6);
\draw (-.25,0) -- (.25,0) (1.25,0) -- (1.75,0) (.25,0) .. controls (.25,-.4) and (1.25,-.4) .. (1.25,0) (1.15,.15) node[right]{\scriptsize $\ell$} (.35,.15) node[left]{\scriptsize $\ell$} (.75,.2) node[above]{\scriptsize $i$} (-.25,0.05) node[below]{\scriptsize $b$} (1.75,0.05) node[below]{\scriptsize $a$};
\draw [very thick, dotted] (.25,0) .. controls (.25,.4) and (1.25,.4) .. (1.25,0);
\fill (.25,0) circle (1.5pt) (1.25,0) circle (1.5pt);
\end{tikzpicture}$$
where half-edges $a$ and $b$ may be dotted reads
$$\frac12\sum_{d(i,.) \in \CD(i)}\sum_{\sigma \in \mathfrak{A}_3}\CI(d(i,.))
\begin{tikzpicture}\useasboundingbox (-.2,-.1) rectangle (2.2,1.5);
\draw (0,0) -- (2,0) (1,0) -- (1,1) arc (-90:270:.2) (1.15,1.2) node[right]{\scriptsize $B(i,\sigma,d(i,.))(t)$} (.95,.35) node[right]{\scriptsize $\Sigma^{(i)}_{d(i,\sigma(1))}$} (0,.05) node[below]{\scriptsize $b$} (2,.05) node[below]{\scriptsize $a$} (1,0) node[below]{\scriptsize $\ell$} (1,.95) node[above]{\scriptsize $i$};
\draw [->] (1,1) arc (-90:0:.2);
\fill (1,0) circle (1.5pt) (1,1) circle (1.5pt);
\end{tikzpicture}$$
meaning that the vertex colored by $\ell$ is constrained to live in $\Sigma^{(i)}_{d(i,\sigma(1))}=\overline{\Sigma}^S(\check{z}^{(i)}_{d(i,\sigma(1))},r)$ that intersects $A_S^{(\ell)}$
as $$\sum_{k=1}^{g(A^{(\ell)})}lk(z_{d(i,\sigma(1))}^{(i)},z_{k}^{(\ell)})\Sigma^S(a^{(\ell)}_{k})$$
and that the associated normal will be assigned to the vertical half-edge labelled by $\Sigma^{(i)}_{d(i,\sigma(1))}$.
\end{lemma}
\bp Recall $B(i,\sigma,d(i,.))=lk_e(\tilde{z}^{(i)}_{d(i,\sigma(2))},\tilde{z}^{(i)}_{d(i,\sigma(3))})-lk_e(\tilde{z}^{(i)}_{d(i,\sigma(3))},\tilde{z}^{(i)}_{d(i,\sigma(2))})$.
Thanks to Lemma~\ref{lemoredgedir}, we may assume that the orientations of $e(r)$ and $e(s)$ match. Then we compute the contribution
of the tiny independent configuration as in the first part of the proof of Proposition~\ref{propthetaiellone}. Thanks to Lemma~\ref{lemdeccoledge},
this contribution reads
$$\sum_{d(i,.) \in \CD(i)}\sum_{\sigma \in \mathfrak{A}_3}\CI(d(i,.))\begin{tikzpicture}\useasboundingbox (0,-.1) rectangle (4,1.2);
\draw [-<] (.2,0) -- (2,0);
\draw (1.95,0) -- (3.8,0) (.6,.05) node[below]{\scriptsize $b$} (3.4,.05) node[below]{\scriptsize $a$} (1.3,0.05) node[below]{\scriptsize $v_2$} (2.7,0) node[below]{\scriptsize $v_1$}  (2,.85) node[above]{\scriptsize $-B(i,\sigma,d(i,.))(t)$} (2,.95) node[below]{\scriptsize $i$} (1.15,.3) node[left]{\scriptsize $N(S)$} (2.85,.3) node[right]{\scriptsize $\Sigma^{(i)}_{d(i,\sigma(1))}$};
\draw [very thick, dotted, ->] (1,0) .. controls (1,.9) and (1.7,.9) .. (2,.9) (3,0) .. controls (3,.9) and (2.3,.9) .. (2,.9);
\fill (1,0) circle (1.5pt) (3,0) circle (1.5pt);
\end{tikzpicture}$$
where the minus sign in front of $B(i,\sigma,d(i,.))(t)$ comes from the fact that $i \in S$.

Now,
$$\begin{tikzpicture}\useasboundingbox (0,-.1) rectangle (4.85,1.2);
\draw [-<] (.2,0) -- (2,0);
\draw (1.95,0) -- (3.8,0) (.6,.05) node[below]{\scriptsize $b$} (3.4,.05) node[below]{\scriptsize $a$} (1.3,0.05) node[below]{\scriptsize $v_2$} (2.7,0) node[below]{\scriptsize $v_1$}  (2,.85) node[above]{\scriptsize $-B(i,\sigma,d(i,.))(t)$} (2,.95) node[below]{\scriptsize $i$} (1.15,.3) node[left]{\scriptsize $N(S)$} (2.85,.3) node[right]{\scriptsize $\Sigma^{(i)}_{d(i,\sigma(1))}$};
\draw [very thick, dotted, ->] (1,0) .. controls (1,.9) and (1.7,.9) .. (2,.9) (3,0) .. controls (3,.9) and (2.3,.9) .. (2,.9);
\fill (1,0) circle (1.5pt) (3,0) circle (1.5pt);
\end{tikzpicture} =\begin{tikzpicture}\useasboundingbox (0,-.1) rectangle (4.85,1.2);
\draw (.2,0) -- (3.8,0) (.6,.05) node[below]{\scriptsize $v_1$} (3.4,.05) node[below]{\scriptsize $a$} (1.3,0.05) node[below]{\scriptsize $v_2$} (1.2,-0.05) node[above]{\scriptsize $+$} (2.7,0) node[below]{\scriptsize $b$}  (2,.85) node[above]{\scriptsize $B(i,\sigma,d(i,.))(t)$} (2,.95) node[below]{\scriptsize $i$} (1.15,.3) node[left]{\scriptsize $N(S)$} (2.85,.3) node[right]{\scriptsize $\Sigma^{(i)}_{d(i,\sigma(1))}$};
\draw [very thick, dotted, ->] (1,0) .. controls (1,.9) and (1.7,.9) .. (2,.9) (3,0) .. controls (3,.9) and (2.3,.9) .. (2,.9);
\fill (1,0) circle (1.5pt) (3,0) circle (1.5pt);
\end{tikzpicture} =\frac12
\begin{tikzpicture}\useasboundingbox (0,-.1) rectangle (3.2,1.5);
\draw (.2,0) -- (1.8,0) (1,0) -- (1,1) arc (-90:270:.2) (1.15,1.2) node[right]{\scriptsize $B(i,\sigma,d(i,.))(t)$} (.95,.35) node[right]{\scriptsize $\Sigma^{(i)}_{d(i,\sigma(1))}$} (.6,.05) node[below]{\scriptsize $b$} (1.4,.05) node[below]{\scriptsize $a$} (1,0) node[below]{\scriptsize $\ell$} (1,.95) node[above]{\scriptsize $i$};
\draw [->] (1,1) arc (-90:0:.2);
\fill (1,0) circle (1.5pt) (1,1) circle (1.5pt);
\end{tikzpicture}.$$
\eop

Now, according to Lemma~\ref{lembeadvan}, Proposition~\ref{propcruc} follows from Lemmas~\ref{lemtinyconf} and \ref{lemconttiny}.\eop

\noindent{\sc Proof of Proposition~\ref{propkeythmflag}:}

We can perform the following operation $\phi_{i,s}$ on a $\theta$-graph whose vertices are colored by $\ell$, with an edge $e(r)$ from  $v(x)$ to $v(y)$ colored by $i$, with another edge
$e(s)$ such that $r<s$, and a third edge $e(t)$. 
Under $\phi_{i,s}$, the edge $e(r)$ becomes looped at the vertex $v_{\phi}(y)$ that is colored by $i$, 
the edge $e(t)$ becomes looped at the vertex $v_{\phi}(x)$ that is colored by $\ell$.
The edge $e(s)$ goes from $v_{\phi}(x)$ to $v_{\phi}(y)$ if $e(s)$ goes from
$v(x)$ to $v(y)$, and from $v_{\phi}(y)$ to $v_{\phi}(x)$, otherwise. By convention, we dash the edge $e(r)$. 

$$\phi_{i,s}\left(
\begin{tikzpicture}\useasboundingbox (-.5,-.1) rectangle (2.5,.85);
\draw (.2,.1) node[left]{\scriptsize $\ell$} (1.8,.1) node[right]{\scriptsize $\ell$} (.2,-.2) node[left]{\scriptsize $v(y)$} (1.75,-.15) node[right]{\scriptsize $v(x)$} (.85,.1) node[above]{\scriptsize $i$} (1,.05) node[below]{\scriptsize $e(s)$} (1,.7) node{\scriptsize $e(r)$} (1,-.7) node{\scriptsize $e(t)$} (1,-.5) .. controls (1.3,-.5) and (1.9,-.5) .. (1.9,0) (1,0) -- (1.9,0);
\draw [very thick, dotted,-<] (.1,0) .. controls (.1,.5) and (.7,.5) .. (1,.5);
\draw [very thick, dotted] (1,.5) .. controls (1.3,.5) and (1.9,.5) .. (1.9,0);
\draw [-<] (.1,0) .. controls (.1,-.5) and (.7,-.5) .. (1,-.5);
\draw [-<] (.1,0) -- (1,0);
\fill (.1,0) circle (1.5pt) (1.9,0) circle (1.5pt);
\end{tikzpicture}
 \right)\;\;=\;\;
\begin{tikzpicture}\useasboundingbox (-.2,-.1) rectangle (3.2,.4);
\draw [->] (.4,0) -- (1.5,0);
\draw [dashed] (2.8,0) circle (.2);
\draw (1.5,0) -- (2.6,0) (.2,0) circle (.2) (.2,.1) node[above]{\scriptsize $e(t)$} (2.8,.1) node[above]{\scriptsize $e(r)$} (1.5,-.05) node[above]{\scriptsize $e(s)$} (2.7,.15) node[left]{\scriptsize $i$} (.3,.15) node[right]{\scriptsize $\ell$}  (2.75,-.2) node[left]{\scriptsize $v_{\phi}(y)$} (.25,-.15) node[right]{\scriptsize $v_{\phi}(x)$} ;
\fill (.4,0) circle (1.5pt) (2.6,0) circle (1.5pt);
\end{tikzpicture}
\;\; =\;\; \phi_{i,s}\left(
\begin{tikzpicture}\useasboundingbox (-.5,-.1) rectangle (2.5,.85);
\draw (.2,.1) node[left]{\scriptsize $\ell$} (1.8,.1) node[right]{\scriptsize $\ell$} (.2,-.2) node[left]{\scriptsize $v(y)$} (1.75,-.15) node[right]{\scriptsize $v(x)$} (.85,.1) node[above]{\scriptsize $i$} (1,.05) node[below]{\scriptsize $e(s)$} (1,.7) node{\scriptsize $e(r)$} (1,-.7) node{\scriptsize $e(t)$} (1,-.5) .. controls (1.3,-.5) and (1.9,-.5) .. (1.9,0) (1,0) -- (1.9,0);
\draw [very thick, dotted,-<] (.1,0) .. controls (.1,.5) and (.7,.5) .. (1,.5);
\draw [very thick, dotted] (1,.5) .. controls (1.3,.5) and (1.9,.5) .. (1.9,0);
\draw [->] (.1,0) .. controls (.1,-.5) and (.7,-.5) .. (1,-.5);
\draw [-<] (.1,0) -- (1,0);
\fill (.1,0) circle (1.5pt) (1.9,0) circle (1.5pt);
\end{tikzpicture}
 \right).$$

$$\phi_{i,s}\left(
\begin{tikzpicture}\useasboundingbox (-.5,-.1) rectangle (2.5,.85);
\draw (.2,.1) node[left]{\scriptsize $\ell$} (1.8,.1) node[right]{\scriptsize $\ell$} (.2,-.2) node[left]{\scriptsize $v(y)$} (1.75,-.15) node[right]{\scriptsize $v(x)$} (.85,.1) node[above]{\scriptsize $i$} (1,.05) node[below]{\scriptsize $e(s)$} (1,.7) node{\scriptsize $e(r)$} (1,-.7) node{\scriptsize $e(t)$} (1,-.5) .. controls (1.3,-.5) and (1.9,-.5) .. (1.9,0) (1,0) -- (1.9,0);
\draw [very thick, dotted,-<] (.1,0) .. controls (.1,.5) and (.7,.5) .. (1,.5);
\draw [very thick, dotted] (1,.5) .. controls (1.3,.5) and (1.9,.5) .. (1.9,0);
\draw [-<] (.1,0) .. controls (.1,-.5) and (.7,-.5) .. (1,-.5);
\draw [->] (.1,0) -- (1,0);
\fill (.1,0) circle (1.5pt) (1.9,0) circle (1.5pt);
\end{tikzpicture}
 \right)\;\;=\;\;
\begin{tikzpicture}\useasboundingbox (-.2,-.1) rectangle (3.2,.4);
\draw [-<] (.4,0) -- (1.5,0);
\draw [dashed] (2.8,0) circle (.2);
\draw (1.5,0) -- (2.6,0) (.2,0) circle (.2) (.2,.1) node[above]{\scriptsize $e(t)$} (2.8,.1) node[above]{\scriptsize $e(r)$} (1.5,-.05) node[above]{\scriptsize $e(s)$} (2.7,.15) node[left]{\scriptsize $i$} (.3,.15) node[right]{\scriptsize $\ell$}  (2.75,-.2) node[left]{\scriptsize $v_{\phi}(y)$} (.25,-.15) node[right]{\scriptsize $v_{\phi}(x)$} ;
\fill (.4,0) circle (1.5pt) (2.6,0) circle (1.5pt);
\end{tikzpicture}
\;\; =\;\; \phi_{i,s}\left(
\begin{tikzpicture}\useasboundingbox (-.5,-.1) rectangle (2.5,.85);
\draw (.2,.1) node[left]{\scriptsize $\ell$} (1.8,.1) node[right]{\scriptsize $\ell$} (.2,-.2) node[left]{\scriptsize $v(y)$} (1.75,-.15) node[right]{\scriptsize $v(x)$} (.85,.1) node[above]{\scriptsize $i$} (1,.05) node[below]{\scriptsize $e(s)$} (1,.7) node{\scriptsize $e(r)$} (1,-.7) node{\scriptsize $e(t)$} (1,-.5) .. controls (1.3,-.5) and (1.9,-.5) .. (1.9,0) (1,0) -- (1.9,0);
\draw [very thick, dotted,-<] (.1,0) .. controls (.1,.5) and (.7,.5) .. (1,.5);
\draw [very thick, dotted] (1,.5) .. controls (1.3,.5) and (1.9,.5) .. (1.9,0);
\draw [->] (.1,0) .. controls (.1,-.5) and (.7,-.5) .. (1,-.5);
\draw [->] (.1,0) -- (1,0);
\fill (.1,0) circle (1.5pt) (1.9,0) circle (1.5pt);
\end{tikzpicture} \right).$$

Let $\Gamma \in \CS_n^u$ be a graph equipped with an admissible coloring $\CaC$.
Let $\CS(\Gamma,\CaC)$ be the set of choices of one edge $e(s)$ with $r<s$ for each $\theta$-case.
For each element $\bf{s}$ in $\CS(\Gamma,\CaC)$, we get a graph $\phi_{\bf{s}}(\Gamma,\CaC) \in \CS^{\ell}_n$ with a simple coloring
by composing the operations $\phi_{i,s}$ associated with $\bf{s}$ and operations $\phi_i$ of Proposition~\ref{propcruc}. In each tadpole-component \tadpoleplain $\,$ of the graph $\phi_{\bf{s}}(\Gamma,\CaC) \in \CS^{\ell}_n$, one looped edge is dashed.
Note that the operation $\phi_i$ must be performed before the operation $\phi_j$ for Configuration $D2$.
For Configuration $D3$, $\phi_i$ and $\phi_k$ must be performed before $\phi_j$, and they commute so that
$\phi_{\bf{s}}(\Gamma,\CaC)$ is well-defined. Let $n(\Gamma,\CaC)$ be the number of colored edges in $(\Gamma,\CaC)$, then $\phi_{\bf{s}}$ is a composition of $n(\Gamma,\CaC)$ operations.

Then according to Propositions~\ref{propthetaiellone} and~\ref{propcruc},
$$ I(\Gamma,\CaC)=\frac{1}{2^{n(\Gamma,\CaC)}}\sum_{\bf{s} \in \CS(\Gamma,\CaC)}I(\phi_{\bf{s}}(\Gamma,\CaC)).$$

Say that a colored graph $(\Gamma^{\ell},\CaC^{\ell})$ of $\CS^{\ell}_n$ is $\theta$-admissible, if $\CaC^{\ell}$ is a simple coloring, if $(\Gamma^{\ell},\CaC^{\ell})$ has one dashed looped edge in each of its tadpole-components, and if, for any such looped edge $e(r)$, the label $s$ of the adjacent non-looped edge $e(s)$ is bigger than $r$.

Let us prove
$$\sum_{\Gamma \in \CS^u_n}\sum_{\CaC\; \mbox{\scriptsize \textrm admissible coloring of}\;\Gamma}I(\Gamma,\CaC)=\sum_{\Gamma^{\ell} \in \CS^{\ell}_n\; ;(\Gamma^{\ell},\CaC^{\ell}) \theta\mbox{\scriptsize \textrm -admissible,} }I(\Gamma^{\ell},\CaC^{\ell}).$$

Starting with a $\theta$-admissible graph $(\Gamma^{\ell},\CaC^{\ell})$ of $\CS^{\ell}_n$, let us determine the colored graphs such that $\phi_{\bf{s}}(\Gamma,\CaC)=(\Gamma^{\ell},\CaC^{\ell})$, and the associated $\bf{s}$.
If $(\Gamma^{\ell},\CaC^{\ell})$ has only one loop, whose vertex is colored by $i$, then there are two ways of writing $(\Gamma^{\ell},\CaC^{\ell})$ as $\phi_i(\Gamma,\CaC)$, they are obtained from each other by reversing the double edge.
$$
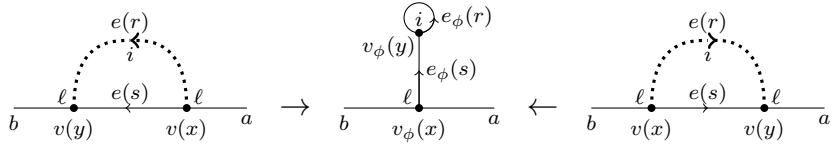
\begin{tikzpicture}\useasboundingbox (0,-.1) rectangle (3.5,.8);
\draw [-<] (.2,0) -- (1.75,0);
\draw (1.7,0) -- (3.3,0) (.2,.05) node[below]{\scriptsize $b$} (3.3,.05) node[below]{\scriptsize $a$} (1,0) node[below]{\scriptsize $v(y)$} (2.5,0) node[below]{\scriptsize $v(x)$} (1.75,-.1) node[above]{\scriptsize $e(s)$}  (1.75,.85) node[above]{\scriptsize $e(r)$} (1.75,.95) node[below]{\scriptsize $i$} (1.05,.15) node[left]{\scriptsize $\ell$} (2.45,.15) node[right]{\scriptsize $\ell$};
\draw [very thick, dotted, ->] (1,0) .. controls (1,.9) and (1.5,.9) .. (1.75,.9)(2.5,0) .. controls (2.5,.9) and (2,.9) .. (1.75,.9);
\fill (1,0) circle (1.5pt) (2.5,0) circle (1.5pt);
\end{tikzpicture}\; \rightarrow \;\begin{tikzpicture}\useasboundingbox (-.2,-.1) rectangle (2.2,1.5);
\draw (0,0) -- (2,0) (1,0) -- (1,1) arc (-90:270:.2) (1.15,1.2) node[right]{\scriptsize $e_{\phi}(r)$} (.95,.5) node[right]{\scriptsize $e_{\phi}(s)$} (0,.05) node[below]{\scriptsize $b$} (2,.05) node[below]{\scriptsize $a$} (1,0) node[below]{\scriptsize $v_{\phi}(x)$} (1,.95) node[above]{\scriptsize $i$} (1.05,.15) node[left]{\scriptsize $\ell$} (1.1,.8) node[left]{\scriptsize $v_{\phi}(y)$};
\draw [->] (1,1) arc (-90:0:.2);
\draw [->] (1,0) -- (1,.5);
\fill (1,0) circle (1.5pt) (1,1) circle (1.5pt);
\end{tikzpicture}\; \leftarrow \;\begin{tikzpicture}\useasboundingbox (0,-.1) rectangle (3.5,.8);
\draw [->] (.2,0) -- (1.75,0);
\draw (1.7,0) -- (3.3,0) (.2,.05) node[below]{\scriptsize $b$} (3.3,.05) node[below]{\scriptsize $a$} (1,0) node[below]{\scriptsize $v(x)$} (2.5,0) node[below]{\scriptsize $v(y)$} (1.75,-.1) node[above]{\scriptsize $e(s)$}  (1.75,.85) node[above]{\scriptsize $e(r)$} (1.75,.95) node[below]{\scriptsize $i$} (1.05,.15) node[left]{\scriptsize $\ell$} (2.45,.15) node[right]{\scriptsize $\ell$};
\draw [very thick, dotted, -<] (1,0) .. controls (1,.9) and (1.5,.9) .. (1.75,.9)(2.5,0) .. controls (2.5,.9) and (2,.9) .. (1.75,.9);
\fill (1,0) circle (1.5pt) (2.5,0) circle (1.5pt);
\end{tikzpicture}$$

Similarly, if $(\Gamma^{\ell},\CaC^{\ell})$ has no meeting lollipops (see below) and if it has $m(\Gamma^{\ell},\CaC^{\ell})$ loops, there are $2^{m(\Gamma^{\ell},\CaC^{\ell})}$ colored graphs of $\CS_n^u$ such that $\phi_{\bf{s}}(\Gamma,\CaC)=(\Gamma^{\ell},\CaC^{\ell})$, and for all these graphs $m(\Gamma^{\ell},\CaC^{\ell})=n(\Gamma,\CaC)$.
If there are several meeting lollipops 

$$
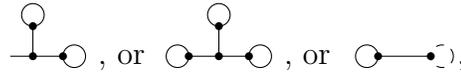
\begin{tikzpicture} \useasboundingbox (.4,-.1) rectangle (1.5,.8);
\draw (.75,.55) circle (.15) (1.3,0) circle (.15) (.45,0) -- (1.15,0) (.75,0) -- (.75,.4);
\fill (.75,0) circle (1.5pt) (1.15,0) circle (1.5pt) (.75,.4) circle (1.5pt);
\end{tikzpicture}\;\mbox{, or}\;\;\begin{tikzpicture} \useasboundingbox (0,-.1) rectangle (1.5,.8);
\draw (.75,.55) circle (.15) (1.3,0) circle (.15) (.2,0) circle (.15) (.35,0) -- (1.15,0) (.75,0) -- (.75,.4);
\fill (.75,0) circle (1.5pt) (1.15,0) circle (1.5pt) (.75,.4) circle (1.5pt) (.35,0) circle (1.5pt);
\end{tikzpicture}\;\mbox{, or}\;\; \tadpole ,$$
first remove all the dashed loops of the tadpoles, and next the loop $e(r)$ with the smallest $r$, thus we get a loop elimination algorithm with $2$ choices at each loop removal that produces $2^k$ graphs $(\Gamma,\CaC)$ of $\CS_n^u$ equipped with admissible colorings and with elements $\bf{s}$ of $\CS(\Gamma,\CaC)$, where $k$ is the number of colored edges of the obtained $(\Gamma,\CaC)$. One easily checks that the $2^k$ obtained $(\Gamma,\CaC,\bf{s})$
are all the $(\Gamma,\CaC,\bf{s})$ such that $\phi_{\bf{s}}(\Gamma,\CaC)=(\Gamma^{\ell},\CaC^{\ell})$.

Now, the proof of Proposition~\ref{propkeythmflag} is reduced to the proof of the following equality. $$\sum_{\Gamma^{\ell} \in \CS^{\ell}_n\; ;(\Gamma^{\ell},\CaC^{\ell})\; \theta\mbox{\scriptsize \textrm -admissible,}}I(\Gamma^{\ell},\CaC^{\ell})=\sum_{\Gamma \in \CS^{\ell}_n}\sum_{\CaC\; \mbox{\scriptsize \textrm simple coloring of}\;\Gamma}I(\Gamma,\CaC). $$
Note that there is a canonical obvious one-to-one correspondence between graphs without tadpoles of the left-hand side and graphs without tadpoles of the right-hand side.
Say that two simply colored graphs $(\Gamma,\CaC)$ and $(\Gamma^{\prime},\CaC^{\prime})$ are equivalent if $(\Gamma^{\prime},\CaC^{\prime})$ is obtained from $(\Gamma,\CaC)$ by permuting some edge labels inside tadpoles. Let $t(\Gamma)$ be the number of tadpoles in 
$\Gamma$.
Since the edges are oriented, there are exactly $6^{t(\Gamma)}$ simply colored graphs  $(\Gamma^{\prime},\CaC^{\prime})$ in the equivalence class of $(\Gamma,\CaC)$. For each of them,
$I(\Gamma^{\prime},\CaC^{\prime})=I(\Gamma,\CaC)$.
Similarly, there are exactly $6^{t(\Gamma)}$ $\theta$-admissible graphs  $(\Gamma^{\prime},\CaC^{\prime})$ in the class of $(\Gamma,\CaC)$. Indeed,
\begin{itemize}
\item for the two permutations of 
the edge labels in a tadpole that put the smallest label in the middle edge, no looped edge can be dashed in an admissible way,
\item for the two permutations of 
the edge labels in a tadpole that put the biggest label in the middle edge, each looped edge can be dashed in an admissible way,
\item for the two remaining permutations, one edge can be dashed. \end{itemize}
\eop

\newpage
\section{Relative Pontrjagin numbers}
\setcounter{equation}{0}
\label{secrelPont}

Recall that any compact oriented $3$-manifold is parallelizable.

\subsection{Definition of relative Pontrjagin numbers}
\label{subpont}
Let $C_0$ and $C_1$ be two compact connected oriented $3$-manifolds whose boundaries  have collars that are identified by a diffeomorphism. Let $\tau_0\colon C_0 \times \RR^3 \rightarrow TC_0$ and $\tau_1\colon C_1 \times \RR^3 \rightarrow TC_1$ be two parallelizations (that respect the orientations) that agree on the collar neighborhoods of $\partial C_0=\partial C_1$. 
Then the {\em relative Pontrjagin number $p_1(\tau_0,\tau_1)$\/} is the Pontrjagin obstruction to extending the trivialization of $TW \otimes \CC$ induced by $\tau_0$ and $\tau_1$ to the interior of a signature $0$ cobordism $W$ from $C_0$ to $C_1$. Details follow.

Let $C$ be a compact connected oriented $3$-manifold. A {\em special complex trivialization\/} of $TC$ is a trivialization of $TC \otimes \CC$ that is obtained from a parallelization $\tau_C \colon C \times \RR^3 \rightarrow TC$ by composing $(\tau^{\CC}_C=\tau_C \otimes_{\RR} \CC ) \colon  C \times \CC^3\rightarrow TC \otimes \CC$ by $$\begin{array}{llll} 
\psi(G): &C \times \CC^3 &\longrightarrow  &C \times \CC^3\\
&(x,y) & \mapsto &(x,G(x)(y))\end{array}$$
for a map $G\colon (C,\partial C) \rightarrow (SL(3,\CC),1)$. The definition and properties of relative Pontrjagin numbers that are given with more details below are valid for pairs of special complex trivializations.

Recall that the \emph{signature} of a $4$-manifold is the signature of the intersection form on its $H_2(.;\RR)$ (number of positive entries minus number of negative entries in a diagonalised version of this form). Also recall that any closed oriented three-manifold bounds a compact oriented $4$-dimensional manifold whose signature may be arbitrarily changed by connected
sums with copies of $\CC P^2$ or $-\CC P^2$.
A {\em cobordism from $C_0$ to $C_1$\/} is a compact oriented $4$-dimensional manifold $W$ with corners such that
$$\partial W=-C_0\cup_{\partial C_0 \sim 0 \times \partial C_0}  (-[0,1] \times \partial C_0)  \cup_{\partial C_1 \sim 1 \times \partial C_0} C_1,$$ 
and $W$ is identified with an open subspace of one of the products $[0,1[ \times C_0$ or $]0,1] \times C_1$ near $\partial W$, as the following picture suggests.

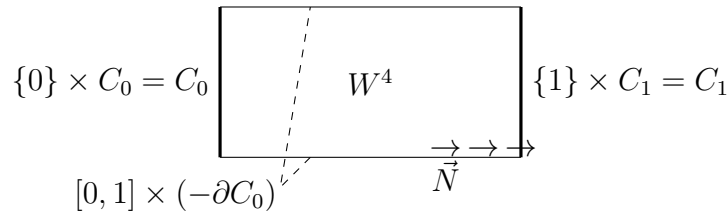
\begin{center}
\begin{tikzpicture}
\useasboundingbox (-3,-.7) rectangle (5.5,2.5);
\draw (0,0) -- (4,0) (4,2) -- (0,2) (2,1) node{$W^4$} (0,1) node[left]{$\{0\} \times C_0=C_0 $} (4,1) node[right]{$\{1\} \times C_1=C_1$};
\draw [very thick] (0,0) -- (0,2) (4,0) -- (4,2);
\draw [dashed] (1.2,0) -- (.8,-.4) -- (1.2,2);
\draw (.9,-.5) node[left]{$[0,1] \times (-\partial C_0)$};
\draw (4,.1) node{$ \rightarrow $} (3.5,.1) node{$ \rightarrow $} (3,.1) node{$ \rightarrow $} node[below](3,-.1) {$ \vec{N}$};
\end{tikzpicture}
\end{center}

Let $W=W^4$ be such a connected cobordism from $C_0$ to $C_1$, with signature $0$.
Consider the complex $4$-bundle $TW \otimes \CC$ over $W$.
Let $\vec{N}$ be the tangent vector to $[0,1] \times \{\mbox{pt}\}$ (under the identifications above), and let 
$\tau(\tau_0,\tau_1)$ denote the trivialization of $TW \otimes \CC$ over $\partial W$ that is obtained by stabilizing either $\tau_0$ or $\tau_1$ into $\vec{N} \oplus \tau_0$ or  $\vec{N} \oplus \tau_1$. Then the obstruction to extending this trivialization to $W$ is the relative first {\em Pontrjagin class\/} $$p_1(W;\tau(\tau_0,\tau_1))=p_1(\tau_0,\tau_1)[W,\partial W] \in H^4(W,\partial W;\ZZ=\pi_3(SU(4)))=\ZZ[W,\partial W]$$ of the trivialization. 

Now, we specify our sign conventions for this Pontrjagin class. They are the same as in \cite{milnorsta}. 
In particular, $p_1$ is the opposite of the second Chern class $c_2$ of the complexified tangent bundle. See \cite[p. 174]{milnorsta}.
More precisely, equip $C_0$ and $C_1$ with Riemannian metrics that coincide near $\partial C_0$, and equip $W$ with a Riemannian metric that coincides with the orthogonal product metric of one of the products $[0,1] \times C_0$ or $[0,1] \times C_1$ near $\partial W$. 
Equip $TW \otimes \CC$ with the associated hermitian structure. The determinant bundle of $TW$
is trivial because $W$ is oriented and $\mbox{det}(TW \otimes \CC)$ is also trivial.
Our parallelization $\tau(\tau_0,\tau_1)$ over $\partial W$ is special with respect to
the trivialization of $\mbox{det}(TW \otimes \CC)$.
Up to homotopy, assume that  $\tau(\tau_0,\tau_1)$
is unitary with respect to the hermitian structure of $TW \otimes \CC$ and the standard hermitian form of $\CC^4$.
Since $\pi_i(SU(4))=\{0\}$ when $i<3$, the trivialization $\tau(\tau_0,\tau_1)$
extends to a special unitary trivialization $\tau$ outside the interior of a $4$-ball $B^4$
and defines 
$$\tau\colon S^3 \times \CC^4 \longrightarrow (TW \otimes \CC)_{|S^3}$$
over the boundary $S^3=\partial B^4$ of this $4$-ball $B^4$.
Over this $4$-ball $B^4$, the bundle $TW \otimes \CC$ admits a trivialization
$$\tau_B\colon B^4 \times \CC^4 \longrightarrow (TW \otimes \CC)_{|B^4}.$$
Then $\tau_B^{-1} \circ \tau(v \in S^3, w \in \CC^4)=(v, \phi(v)(w))$ for a map $\phi: S^3 \longrightarrow  SU(4)$.

Let $i^2(m^{\CC}_r)$ 
be the following map
$$\begin{array}{llll} i^2(m^{\CC}_r): &(S^3 \subset \CC^2) & \longrightarrow & SU(4)\\
& (z_1,z_2) & \mapsto & \left[\begin{array}{cccc} 1&0&0&0\\0&1&0&0\\0&0&z_1&-\overline{z}_2\\0&0&z_2& \overline{z}_1\end{array} \right] \end{array}.$$
When $(e_1,e_2,e_3,e_4)$ is the standard basis of $\CC^4$, the columns of the matrix contain the coordinates of the images of the $e_i$ with respect to $(e_1,e_2,e_3,e_4)$.
Then the homotopy class $[i^2(m^{\CC}_r)]$ of $i^2(m^{\CC}_r)$ generates $\pi_3(SU(4))=\ZZ[i^2(m^{\CC}_r)]$. Define $p_1(\tau_0,\tau_1)$ from the homotopy class $[\phi]$ of $\phi$ so that
$$[\phi]=-p_1(\tau_0,\tau_1)[i^2(m^{\CC}_r)] \in \pi_3(SU(4)).$$

\begin{proposition}
\label{proppontdef}
Let $C_0$ and $C_1$ be two compact connected oriented $3$-manifolds whose boundaries have collars that are identified by a diffeomorphism. Let $\tau_0 \colon C_0 \times \CC^3 \rightarrow TC_0 \otimes \CC  $ and $\tau_1\colon C_1 \times \CC^3 \rightarrow TC_1 \otimes \CC$ be two special complex trivializations (that respect the orientations) that agree on the collar neighborhoods of $\partial C_0=\partial C_1$. 

The first {\em Pontrjagin number\/} $p_1(\tau_0,\tau_1)$ is well-defined by the above 
conditions.
\end{proposition}
\bp 
According to the Nokivov additivity theorem, if a closed (compact, without boundary) $4$-manifold $Y$ reads $Y= Y^+ \cup_{X} Y^-$
where $Y^+$ and $Y^-$ are two $4$-manifolds with boundary, embedded in $Y$ that intersects along a closed $3$--manifold $X$ (their common boundary, up to orientation) then $$\mbox{signature}(Y)=\mbox{signature}(Y^+)+\mbox{signature}(Y^-).$$ According to a Rohlin theorem (see \cite{Roh} or \cite[p. 18]{gm2}),
when $Y$ is a compact oriented $4$--manifold without boundary, then $p_1(Y)=3\;\mbox{signature}(Y)$.

We only need to prove that $p_1(\tau_0,\tau_1)$ is independent of the signature $0$ cobordism $W$. Let $W_E$
be a $4$-manifold of signature $0$ bounded by $(-\partial W)$. Then $W \cup_{\partial W} W_E$ is a $4$-dimensional manifold without boundary whose signature is $\left(\mbox{signature}(W_E) + \mbox{signature}(W)=0\right)$ by the Novikov additivity theorem.
According to the Rohlin theorem, the first Pontrjagin class of $W \cup_{\partial W} W_E$ is also zero. On the other hand,
this first Pontrjagin class is the sum of the relative first Pontrjagin classes of $W$ and $W_E$ with respect to $\tau(\tau_0,\tau_1)$. These two relative Pontrjagin classes are opposite and therefore the relative first Pontrjagin class of $W$ with respect to $\tau(\tau_0,\tau_1)$ does not depend on $W$.
\eop

\begin{remark}
 When $\partial C_1=\emptyset$ and when $C_0=\emptyset$, the map $p_1$ coincides with the map $h$ that is studied by Hirzebruch in \cite[\S 3.1]{hirzebruchEM}, and by Kirby and Melvin in \cite{km} under the name of {\em Hirzebruch defect.\/} 
\end{remark}

\subsection{Properties of relative Pontrjagin numbers}

\begin{proposition}
\label{proppont}
Let $C_0$ and $C_1$ be two compact connected oriented $3$-manifolds whose boundaries  have collars that are identified by a diffeomorphism. Let $\tau_0 \colon C_0 \times \CC^3 \rightarrow TC_0 \otimes \CC  $ and $\tau_1 \colon C_1 \times \CC^3 \rightarrow TC_1 \otimes \CC  $ be two special complex trivializations (that respect the orientations) that agree on the collar neighborhoods of $\partial C_0=\partial C_1$. 

The first Pontrjagin number $p_1(\tau_0,\tau_1)$  satisfies the following properties.
\begin{enumerate}
\item Let $C_2$ be a compact $3$-manifold whose boundary has a collar neighborhood identified with a collar neighborhood of $\partial C_0$.
Let $\tau_2$ be a special complex trivialization of $TC_2$ that agrees with $\tau_0$ near $\partial C_2$.
If two of the Lagrangians of $C_0$, $C_1$ and $C_2$ coincide in $H_1(\partial C_0;\QQ)$, then  $$p_1(\tau_0,\tau_2)=p_1(\tau_0,\tau_1)+p_1(\tau_1,\tau_2).$$
In particular, since $p_1(\tau_0,\tau_0)=0$, $p_1(\tau_1,\tau_0)=-p_1(\tau_0,\tau_1)$.
 \item 
Let $D$ be a connected compact $3$-manifold that contains $C_0$ in its interior, and let $\tau_D$ be a special complex trivialization of $TD$ that restricts as the special complex trivialization $\tau_0$ on $TC_0$, let $D_1$ be obtained from $D_0$ by replacing $C_0$ by $C_1$, and let $\tau_{D_1}$ be the trivialization
of $TD_1$ that agrees with $\tau_1$ on $TC_1$ and with $\tau_D$ on $T(D\setminus C_0)$. If the Lagrangians of $C_0$ and $C_1$ coincide,
then $$p_1(\tau_D,\tau_{D_1})=p_1(\tau_0,\tau_1).$$
\item
For any 
$$g: (C_1, \partial C_1) \longrightarrow (SU(3),1),$$
define
$$\begin{array}{llll} 
\psi(g): &C_1 \times \CC^3 &\longrightarrow  &C_1 \times \CC^3\\
&(x,y) & \mapsto &(x,g(x)(y))\end{array}$$
and 
$$p^{\prime}_1(g)= p_1(\tau^{\CC}_0,\tau^{\CC}_1 \circ \psi(g))-p_1(\tau^{\CC}_0,\tau^{\CC}_1)=p_1(\tau_1^{\CC},\tau^{\CC}_1 \circ \psi(g)).$$
Then $p^{\prime}_1(g)$ does not depend on the special complex trivialization $\tau_1^{\CC}$ of $TC_1$. The group structure of $SU(3)$ induces a group structure on the set $[(C_1, \partial C_1) , (SU(3),1)]$ of homotopy classes   of maps from $C_1$ to $SU(3)$ that map $\partial C_1$ to $1$,
and $p^{\prime}_1$ defines an isomorphism from $[(C_1, \partial C_1) , (SU(3),1)]$  to $\ZZ$. 
\item
For any 
$$g: (C_1, \partial C_1) \longrightarrow (SO(3),1),$$
let $\mbox{deg}(g)$ denote the degree of $g$ and let 
$$\begin{array}{llll} 
\psi_{\RR}(g): &C_1 \times \RR^3 &\longrightarrow  &C_1 \times \RR^3\\
&(x,y) & \mapsto &(x,g(x)(y))\end{array}$$
then 
$$p_1(\tau_0,\tau_1 \circ \psi_{\RR}(g))-p_1(\tau_0,\tau_1)=p_1(\tau_1,\tau_1 \circ \psi_{\RR}(g))=-p_1(\tau_1 \circ \psi_{\RR}(g),\tau_1)=2\mbox{deg}(g).$$ 

\item

If $C_0$ is the unit ball of $\RR^3$, if $\tau_0$ is the standard parallelization of $\RR^3$, and
if $C_1$ has the same integral homology as a point, then $p_1(\tau_0,.)$ defines a bijection from the set of homotopy classes of parallelizations of $C_1$ that are standard near $\partial C_1=S^2$ to $4 \ZZ$.

\end{enumerate}
\end{proposition}

The proof will use a weak form of the Wall Non-Additivity theorem. We quote the weak form we need.

\begin{theorem}[\cite{wall}]
\label{thmwall}
Let $Y$ be a compact oriented $4$--manifold (with possible boundary), and let $X$ be a three manifold properly embedded in $Y$ that separates $Y$ and that induces the splitting  $Y = Y^+ \cup_{X} Y^-$, for two $4$-manifolds $Y^+$ and $Y^-$ in $Y$, whose intersection is $X$, as in the following figure of $Y$:
\begin{center}
\begin{tikzpicture}
\useasboundingbox (-3,-.4) rectangle (3,.4);
\draw (0,-.4) -- (0,.4);
\draw (0,-.4) .. controls (1.2,-.4)  .. (1.2,0)  .. controls (1.2,.4) .. (0,.4);
\draw (0,-.4) .. controls (-1.2,-.4)  .. (-1.2,0)  .. controls (-1.2,.4) .. (0,.4);
\draw (-1.1,-.1) node[left]{\footnotesize $X^-$} (1.1,-.1) node[right]{\footnotesize $X^+$} (-.1,-.1) node[right]{\footnotesize $X$} (.7,0) node{\footnotesize $Y^+$} (-.5,0) node{\footnotesize $Y^-$};
\end{tikzpicture}
\end{center}

Set
$$X^+ =\overline {\partial Y^+ \setminus (-X)}\;\;\;\mbox{and}\;\;\;X^- =-\overline {\partial Y^- \setminus X}.$$ 
Let $\CL$, $\CL^-$ and $\CL^+$ respectively denote the Lagrangians of $X$, $X^-$ and $X^+$, they are the Lagrangian subspaces of $H_1(\partial X, \QQ)$.
Then $$\left(\mbox{signature}(Y)-\mbox{signature}(Y^+)-\mbox{signature}(Y^-)\right)$$ is the signature of an explicit quadratic form on 
$$\frac{\CL \cap \left(\CL^- + \CL^+ \right)}{(\CL \cap \CL^-) + (\CL \cap \CL^+)}.$$
Furthermore, this space is isomorphic to $\frac{\CL^+ \cap \left(\CL + \CL^- \right)}{(\CL^+ \cap \CL) + (\CL^+ \cap \CL^-)}$ and $\frac{\CL^- \cap \left(\CL + \CL^+ \right)}{(\CL^- \cap \CL) + (\CL^- \cap \CL^+)}$.
\end{theorem}
We use this theorem in cases when the space above is trivial. This is why we do not make the involved quadratic form explicit.

\noindent {\sc Proof of Proposition~\ref{proppont}:} 
Let us prove the first property. Let $Y^-=W$ be a signature $0$ cobordism from $X^-=C_0$ to $X=C_1$, and let $Y^+$ be a signature $0$ cobordism from $C_1$ to $X^+=C_2$. Then it is enough to prove that the signature of $Y=Y^+\cup_X Y^-$ is zero. With the notation of Theorem~\ref{thmwall}, under our assumptions, the space $\frac{\CL \cap \left(\CL^- + \CL^+ \right)}{(\CL \cap \CL^-) + (\CL \cap \CL^+)}$
is trivial, therefore, according to the Wall theorem, the signature of $Y$ is zero. The first property follows.

We now prove that under the assumptions of the second property, $p_1(\tau_D,\tau_{D_1})=p_1(\tau_0,\tau_1)$.
Let $Y^+=\left([0,1] \times (D\setminus \mathring{C}_0)\right)$, let  $Y^-=W$ be a signature $0$ cobordism from $C_0$ to $C_1$, and
let $X=-[0,1]\times \partial C_0$. Note that the signature of $Y^+$ is zero.
In order to prove the wanted equality, it is enough to prove that the signature of $Y=Y^+\cup_X Y^-$ is zero.
Here, $H_1(\partial X;\QQ)=H_1(\partial C_0) \oplus H_1(\partial C_0)$.
Let $j \colon H_1(\partial C_0) \rightarrow H_1(D\setminus \mathring{C}_0)$ and let $j_{\partial D} \colon H_1(\partial D) \rightarrow H_1(D\setminus \mathring{C}_0)$ be the maps induced by inclusions. With the notation of Theorem~\ref{thmwall},

$$\begin{array}{ll}
\CL&=\{(x,-x);x \in H_1(\partial C_0)\}\\ 
\CL^-&=\{(x,y);x \in \CL_{C_0}, y \in \CL_{C_1}\}\\ 
\CL^+&= \{(x,y);(j(x),j(y))=(j_{\partial D}(z \in H_1(\partial D)),-j_{\partial D}(z))\}\\
&=\{(y,-y);j(y) \in \mbox{Im}(j_{\partial D})\} \oplus \{(x,0);j(x) =0\}\\
\CL \cap \CL^-  &= \{(x,-x);x \in \CL_{C_0} \cap \CL_{C_1}=\CL_{C_0}\}\\
\CL \cap \CL^+  &= \{(x,-x);j(x) \in \mbox{Im}(j_{\partial D})\}
\end{array}$$
Let us prove that $\CL \cap (\CL^- + \CL^+)=(\CL \cap \CL^-) + (\CL \cap \CL^+)$.
For a subspace $K$ of $H_1(\partial C_0;\QQ)$, set $j_{MV}(K)=\{(x,-x);x \in K\} \subset H_1(\partial X;\QQ)$. Then $\CL=j_{MV}(H_1(\partial C_0))$ and $\CL \cap \CL^+=j_{MV}\left(j^{-1}\left(\mbox{Im}(j_{\partial D})\right)\right)$.
$$\CL \cap (\CL^- + \CL^+) = \CL \cap \CL^+  + j_{MV}\left(\CL_{C_1} \cap (\CL_{C_0} + \mbox{Ker}(j))\right).$$
Since $\CL_{C_0}=\CL_{C_1}$, $\left(\CL_{C_1} \cap (\CL_{C_0} + \mbox{Ker}(j))\right)= \CL_{C_0}$ so that $\CL \cap (\CL^- + \CL^+)=(\CL \cap \CL^+) +j_{MV}(\CL_{C_0}).$
Then the second property is proved, thanks to the Wall theorem that guarantees the additivity of the signature in this case.

The equality $p_1(\tau^{\CC}_0,\tau^{\CC}_1 \circ \psi(g))-p_1(\tau^{\CC}_0,\tau^{\CC}_1)=p_1(\tau_1^{\CC},\tau^{\CC}_1 \circ \psi(g))$ comes from the first property. Then it is easy to observe that $p_1(\tau_1^{\CC},\tau^{\CC}_1 \circ \psi(g))$ does not depend on $\tau^{\CC}_1$. Once it is observed, the fact that $p^{\prime}_1$ is a group homomorphism comes from the first property again.
Now, since $\pi_i(SU(3))$ is trivial for $i<3$ and since $\pi_3(SU(3))=\ZZ$, the group $[(C_1, \partial C_1) , (SU(3),1)]$ is isomorphic to $\ZZ$ and it is generated by the class of a map that maps the complement of a $3$-ball $B$ to $1$ and that factors through a map that generates $\pi_3(SU(3))$. By definition of the Pontrjagin classes, such a generator is sent to $\pm 1$. 

In order to prove the fourth assertion, for an oriented connected $3$-manifold $C$ with possible boundary, and equipped with a parallelization, similarly define the map
$$\begin{array}{llll} p^{\prime}_1 \colon &[(C,\partial C),(SO(3),1)] &\rightarrow & \ZZ \\
   &[g]&\mapsto & p_1(\tau,\tau \circ \psi_{\RR}(g))
  \end{array}$$from the group of homotopy classes
$[(C,\partial C),(SO(3),1)]$ to $\ZZ$.
This map is independent of $\tau$. Furthermore, it is a group homomorphism. Therefore, according to the following independent lemma~\ref{lempreptrivun},
it is proportional to the degree and it is enough to compute it when $g$ is the map $G_C(\rho)$ described in Lemma~\ref{lempreptrivun}, and this is done in \cite[Proposition 1.8]{lesconst}. (See also Lemmas 2.36 and 2.39 in \cite{lesconst}, but be warned that what is denoted by $\tau$ there would be denoted by $\tau^{-1}$ here.)

The fifth assertion is proved in \cite[Lemma 2.40]{lesconst}.
\eop

Thus we are left with the proof of the following version of Lemma 2.32 of \cite{lesconst}.

\begin{lemma}
\label{lempreptrivun}
See $S^3$ as $B^3/\partial B^3$ and see $B^3$ as $([0,2\pi] \times S^2)/(0 \sim \{0\} \times S^2)$.
Let $\chi_{\pi}\colon [0,2\pi] \to [0,2\pi]$ be an increasing smooth bijection whose derivatives
vanish at $0$ and $2\pi$.
Let $\rho \colon B^3 \rightarrow SO(3)$ map $(\theta \in[0,2\pi], v \in S^2)$ to the rotation with axis directed by $v$ and with angle $\chi_{\pi}(\theta)$. This map induces the double covering $\tilde{\rho} \colon S^3 \rightarrow SO(3)$.
Let $C$ be an (oriented, compact) connected $3$-manifold with possible boundary. Let $G_C(\rho):C \longrightarrow SO(3)$ 
be a map that sends the complement of a ball $B^3$ of the interior of $C$ to the identity, and that coincides with $\rho$ on this 3-ball. Note that all such maps induce the same element
$[G_C(\rho)]$ in the group of relative homotopy classes $[(C,\partial C),(SO(3),1)]$. Also note that the elements of $[(C,\partial C),(SO(3),1)]$ have a well-defined degree.
\begin{enumerate}
\item Any homotopy class of a map $G$ from $(C,\partial C)$ to $(SO(3),1)$, such that 
$$\pi_1(G): \pi_1(C) \longrightarrow \pi_1(SO(3))\cong \ZZ/2\ZZ$$
is trivial, belongs to the subgroup $<[G_C(\rho)]>$ of 
$[(C,\partial C),(SO(3),1)]$
generated by $[G_C(\rho)]$. 
\item For any $[G] \in [(C,\partial C),(SO(3),1)]$, 
$$[G]^2 \in <[G_C(\rho)]>.$$
\item The group $[(C,\partial C),(SO(3),1)]$ is abelian.
\item The degree is a group homomorphism from $[(C,\partial C),(SO(3),1)]$ to $\ZZ$.
\item The morphism $$\begin{array}{llll}\frac{\mbox{deg}}{2}:&[(C,\partial C),(SO(3),1)]\otimes_{\ZZ} \QQ &\longrightarrow
&\QQ[G_C(\rho)]\\
&[g] \otimes 1 &\mapsto &\frac{\mbox{deg}(g)}{2}[G_C(\rho)]\end{array}$$
is an isomorphism.
\end{enumerate}
\end{lemma}
\bp
For a connected topological group $\Gamma$, the $\pi_i$-group structure of $\pi_i(\Gamma,1)$ coincides with its group structure induced by the multiplication of maps using the group structure of $\Gamma$.
Also recall that $\mbox{deg} \colon \pi_3(SO(3)) \rightarrow 2\ZZ$ is an isomorphism.

 Assume that $\pi_1(G)$ is trivial. Choose a cell decomposition of $C$ with respect to its boundary, with only one three-cell, no zero-cell if $\partial C \neq \emptyset$, one zero-cell if $\partial C = \emptyset$,  one-cells, and two-cells. Then after a homotopy relative to $\partial C$, we may assume that $G$ maps the one-skeleton of $C$ to $1$. Next, since $\pi_2(SO(3)) = 0$, we may assume that $G$ maps the two-skeleton of $C$ to $1$, and therefore that $G$ maps the exterior of some $3$-ball to $1$. Now $G$ becomes a map from $B^3/\partial B^3=S^3$ to $SO(3)$, and its homotopy class is $k[\tilde{\rho}]$ in $\pi_3(SO(3))=\ZZ[\tilde{\rho}]$, where $(2k)$ is the degree of the map $G$ from $S^3$ to $SO(3)$. Therefore $G$ is homotopic to $G_C(\rho)^k$. This proves the first assertion.

Since $\pi_1(G^2)=2\pi_1(G)$ is trivial, the second assertion follows.

For the third assertion, first note that $[G_C(\rho)]$ belongs to the center of
$[(C,\partial C),(SO(3),1)]$ because it can be supported in a small ball disjoint 
from the support (preimage of $SO(3) \setminus \{1\}$) of a representative of any other element. Therefore, according to the second assertion any square will be in the center. Furthermore, since any commutator induces the trivial map on $\pi_1(M)$, any commutator is in $<[G_C(\rho)]>$.
In particular, if $f$ and $g$ are elements of $[(C,\partial C),(SO(3),1)]$, 
$$(gf)^2=(fg)^2=(f^{-1}f^2g^2f)(f^{-1}g^{-1}fg)$$ 
where the first factor equals $f^2g^2=g^2f^2$. Exchanging $f$ and $g$ yields
$f^{-1}g^{-1}fg=g^{-1}f^{-1}gf$. Then the commutator that is a power of $[G_C(\rho)]$ has a vanishing square, and thus a vanishing degree. Then it must be trivial.

For the fourth assertion, it is easy to see that $\mbox{deg}(fg)=\mbox{deg}(f)+\mbox{deg}(g)$ when $f$ or $g$
is a power of $[G_C(\rho)]$.
Let us prove that $\mbox{deg}(f^k)=k\mbox{deg}(f)$ for any $f$.
In general, there is an unoriented embedded surface $S_f$ of the interior of $C$ such that the restriction of $f$ to a loop $\gamma$ is homotopic to a constant if and only the intersection mod 2 of $S_f$ and $\gamma$ is $1$. Choose a trivialization $\tau_C$ of $TC$ and consider a tubular neighborhood of $S_f$. Locally, the tubular neighborhood reads $[-1,1] \times D$ for a disk $D$ of $S_f$. This defines a local coorientation for $D$. Then define the map from $[-1,1] \times D$ to $SO(3)$ that maps $(t,x)$ to the rotation of angle $(\pi (t+1))$ and with axis directed by the positive normal to $D$. Since changing the coorientation of $D$ reverses both $t$ and the direction of the normal, this process defines a map $g(\tau_C,S_f)$ that is homotopic to $f$ on the two-skeleton of $C$. Therefore $f$ is homotopic to $g(\tau_C,S_f)$ composed with some $G_C(\rho)^r$, and we are left with the proof that the degree of $g^k$ is $k\mbox{deg}(g)$
for $g=g(\tau_C,S_f)$. This can easily be done by noticing that $g^k$ reads $g(\tau_C,S^{(k)}_f)$ where $S^{(k)}_f$ is the $k$-covering of 
$S_f$ embedded in the tubular neighborhood of $S_f$ and intersecting each $[-1,1] \times D$ as $k$ parallel copies of $D$. In general,
$\mbox{deg}(fg)=\frac{1}2\mbox{deg}((fg)^2)=\frac{1}2\mbox{deg}(f^2g^2)=\frac{1}2\left( \mbox{deg}(f^2)+\mbox{deg}(g^2)\right)$, and the fourth assertion is proved.

In particular, since $$\frac{\mbox{deg}}{2}\colon [(C,\partial C),(SO(3),1)]\otimes_{\ZZ} \QQ \rightarrow
\QQ[G_C(\rho)]$$
is a homomorphism that maps $[G_C(\rho)]$ to $[G_C(\rho)]$, it
is a surjective morphism, and since $\mbox{deg}(g)=0 \Rightarrow \mbox{deg}(g^2)=0 \Rightarrow [g^2]=0$, the last assertion follows, too.
\eop

\newpage
\section{Alternative definition of $\tilde{Z}$ with differential forms}
\label{secaltdefzform}
\setcounter{equation}{0}

In this section, we present an alternative definition of $\tilde{Z}$ with differential forms. This section has its own interest. It will be used in the proof of Theorem~\ref{thmdefinvpseudo} in order to refer to some computations that were performed in this context in \cite{lessumgen}.
It will also be used to study the behaviour of $\tilde{Z}$ under connected sum in Section~\ref{secconnsum}.

We use the notation of Subsection~\ref{subcyc}.
For a parallelization $\tau\colon \EXT \times \RR^3 \rightarrow T\EXT$ that coincides with the standard parallelization $\tau_{\pi}$ on $\EXT_{]1,4]} \times \RR^3$,
define a map $$\pi(\tau) \colon \left(\partial \tilde{C}_2(\EXT)\cup (\tilde{C}_2(\EXT) \setminus \tilde{C}_2(\EXT_{[0,2[}))\right) \rightarrow S^2$$
\begin{itemize}
\item that coincides with the map $\pi$ on $(\tilde{C}_2(\EXT) \setminus \tilde{C}_2(\EXT_{[0,2[}))$,
\item that is the natural projection of $U\EXT_{[0,2[}$ to $S^2$ induced by $\tau^{-1}$ on $U\EXT_{[0,2[}=\{0\} \times U\EXT_{[0,2[}$, 
\item that maps $(\NN \setminus \{0\}) \times U\EXT_{[0,2[}$ to the upward vertical vector, and $(-\NN \setminus \{0\}) \times U\EXT_{[0,2[}$ to the downward vertical vector.
\end{itemize}

Let $\omega_i$ be a $2$--form on $S^2$, supported in $S^2_H$, such that $\int_{S^2}\omega_i=1$, and let $\pi(\tau)^{\ast}(\omega_i)$ denote its pull-back under $\pi(\tau)$ on $\left(\partial \tilde{C}_2(\EXT)\cup (\tilde{C}_2(\EXT) \setminus \tilde{C}_2(\EXT_{[0,2[}))\right)$.

Let $\sK_i \in D_{\{1\}}$, let $D(\sK_i) \times S^1$ be a tubular neighborhood of $\sK_i \times S^1$.
Define 
$$p(\sK_i)\colon\partial \tilde{C}_2(\EXT)\cup (\tilde{C}_2(\EXT) \setminus \tilde{C}_2(\EXT_{[0,2[})) \rightarrow S^2$$
so that $p(\sK_i)$ maps the complement of $\{0\} \times U(D(\sK_i) \times S^1)$ to the upward vertical vector $N$, and $p(\sK_i)$ factors through the natural projection from $\{0\} \times U(D(\sK_i) \times S^1)$ to $D(\sK_i)$ on $\{0\} \times U(D(\sK_i) \times S^1)$.
Let $\omega_{i,D}$ be a $2$-form with compact support on $S^2 \setminus \{N\}$ such that $\int_{S^2}\omega_{i,D}=1$.

\begin{lemma}
\label{lemOmegai}
Recall $J_{\Delta(K)}(t)=\frac{t\Delta(K)^{\prime}(t)}{\Delta(K)(t)}$.
 There exists a closed $2$-form $\Omega_i$ with compact support on $\tilde{C}_2(\EXT)$ that restricts on $\left(\partial \tilde{C}_2(\EXT_{[0,2[})\cup (\tilde{C}_2(\EXT) \setminus \tilde{C}_2(\EXT_{[0,2[})\right)$
as $$\delta(K)(\thetap^{\ast})\pi(\tau)^{\ast}(\omega_i) + (\delta(K)J_{\Delta(K)})(\thetap^{\ast})(p(\sK_i)^{\ast}(\omega_{i,D})).$$
\end{lemma}
\bp
The obstruction to such an existence would be cohomological. Since Proposition~\ref{propdefbord} guarantees the existence of a form that is Poincar\'e dual to a natural extension of $C(\cvarM_i,\sK_i,\tau)$ and that satisfies the above boundary conditions, this obstruction vanishes.
\eop

Now, fix a graph $\Gamma$ of $\CS_n^u$ and consider the following associated covering $\tilde{C}_{2n}(\EXT,\Gamma)$ of $C_{2n}(\EXT)$.
The fiber over a configuration $(m_1,\dots,m_{2n})$ of $2n$ distinct points of $\EXT$ is the space of $3n$-tuples $([\gamma_1],\dots,[\gamma_{3n}])$ of rational homology classes of paths $\gamma_i$ associated with edges such that if the edge $e(i)$ goes from $v(j)$ to $v(k)$, then $\gamma_i\colon [0,1] \rightarrow \EXT$ is a path of $\EXT$ such that $\gamma_i(0)=m_j$ and $\gamma_i(1)=m_k$,
and $[\gamma_i]$ is the rational homology class of $\gamma_i$ that contains all paths $\gamma^{\prime}_i$ from $m_j$ to $m_k$ such that 
$\gamma^{\prime}_i\overline{\gamma}_i$ vanishes in $H_1(\EXT;\QQ)$.

Then there are natural projections $\tilde{p}(\Gamma,i) \colon \tilde{C}_{2n}(\EXT,\Gamma) \rightarrow \tilde{C}_2(\EXT)$ associated with the edges. These projections equip
$ \tilde{C}_{2n}(\EXT,\Gamma)$ with the $2n$-form $$\Omega_{\Gamma}((\Omega_i)_i)=\bigwedge_{i=1}^{3n}\tilde{p}(\Gamma,i)^{\ast}(\Omega_i).$$

\begin{lemma}
The support of $\Omega_{\Gamma}((\Omega_i)_i)$ is compact.
\end{lemma}
\bp Since the support of $\Omega_j$ is compact, the map $(x,y) \mapsto |\functE(y)-\functE(x)|$ is bounded by some integer $M_j$ for every $j$, on the support of $\Omega_j$. Let $\mbox{Supp}(\Omega_{\Gamma}((\Omega_i)_i))$ denote the support of $\Omega_{\Gamma}((\Omega_i)_i)$. Consider the map $$\begin{array}{llll}p_j \colon &\mbox{Supp}(\Omega_{\Gamma}((\Omega_i)_i)) &\rightarrow &[-M_j,M_j]\\ 
&([\gamma_1],\dots,[\gamma_{3n}])&\mapsto&\functE(\tilde{\gamma}_j(1))-\functE(\tilde{\gamma}_j(0))                                \end{array}
$$
The union of $(\ZZ +[-1/4,1/4])$ and $(\ZZ +[1/4,3/4])$ covers $\RR$.
Let $$C_{j,0}=p_j^{-1}\left([-M_j,M_j]\cap( \ZZ +[-1/4,1/4])\right)$$ and $$C_{j,1/2}=p_j^{-1}\left([-M_j,M_j]\cap( \ZZ +[1/4,3/4])\right)$$
Then $$\mbox{Supp}(\Omega_{\Gamma}((\Omega_i)_i))=C_{j,0} \cup C_{j,1/2}=\bigcup_{L\colon \{1,2,\dots,3n\} \rightarrow \{0,1/2\}} \cap_{i=1}^{3n} C_{i,L(i)}$$
where each $\cap_{i=1}^{3n} C_{i,L(i)}$ is a union of finitely many copies of a closed subset of the compact ${C}_{2n}(\EXT)$. 
\eop

Recall that $H_1(\EXT;\ZZ)/\mbox{Torsion}$ is identified with $\ZZ$ where $1$ is the image of the oriented meridian of $K$ that generates $H_1(\EXT;\ZZ)/\mbox{Torsion}$.
Then an element of $\tilde{C}_{2n}(\EXT,\Gamma)$ induces the map from $H_1(\Gamma;\ZZ)$ to $\ZZ$ that maps a loop of $\Gamma$ to the homology class of the composition of the paths associated with the edges of the loop (these paths are reversed for edges with coefficient $(-1)$).
In particular, $\tilde{C}_{2n}(\EXT,\Gamma)$ splits as a disjoint union over the elements $\rho$ of $H^1(\Gamma;\ZZ)$ of the subsets $\tilde{C}_{2n}(\EXT,\Gamma,\rho)$ of 
$\tilde{C}_{2n}(\EXT,\Gamma)$ made of the elements inducing $\rho$.

Let $\Gamma_{\rho,\delta(K)}$ be obtained from the element $\Gamma_{\rho}$ of Lemma~\ref{lemgenhol} by dividing all the beads by ${\delta(K)(t)}$.

Now, set $I_{\Gamma,\rho}((\Omega_i)_i)=\int_{\tilde{C}_{2n}(\EXT,\Gamma,\rho)}\Omega_{\Gamma}((\Omega_i)_i)$ and define the following elements of $\CA^h_n(\delta) \otimes_{\QQ} \RR$:
$$I_{\Gamma}((\Omega_i)_i)=\sum_{\rho \in H^1(\Gamma;\ZZ)} I_{\Gamma,\rho}((\Omega_i)_i)[\Gamma_{\rho,\delta(K)}],$$
$$I_n((\Omega_i)_i)=\sum_{\Gamma \in \CS^u_n}\frac{I_{\Gamma}((\Omega_i)_i)}{2^{3n}(3n)!(2n)!}\;\;\;\;\;\;\; \mbox{and}\;\;\;\;\;\;\;
I^c_n((\Omega_i)_i)=\sum_{\Gamma \in \CS_n}\frac{I_{\Gamma}((\Omega_i)_i)}{2^{3n}(3n)!(2n)!}$$
where $\CS^u_n$ and $\CS_n$ are defined in the beginning of Subsection~\ref{subeqint}.

\begin{lemma}
Let $\{\deleqprop_i(\tau)=\delta(K)(t)\eqprop_i(\tau)\}_i$ be a collection of cycles of $(\tilde{C}_{2}(\EXT), \partial \tilde{C}_{2}(\EXT))$ as in Theorem~\ref{thminvwithtau}. For $\rho \in H^1(\Gamma;\ZZ)$, let $I_{\Gamma,\rho}(\{\eqprop_i(\tau)\})$ denote the algebraic intersection of the $\tilde{p}(\Gamma,i)^{-1}(\deleqprop_i(\tau))$ in the component $\tilde{C}_{2n}(\EXT,\Gamma,\rho)$, then 
$$I_{\Gamma}(\{\eqprop_i(\tau)\})=\sum_{\rho \in H^1(\Gamma;\ZZ)} I_{\Gamma,\rho}(\{\eqprop_i(\tau)\})[\Gamma_{\rho,\delta(K)}].$$
\end{lemma}
\bp The intersection points coincide, and their coefficients also do.
\eop

\begin{lemma}
\label{lemhfourtc}
$H_4(\tilde{C}_2(\EXT);\RR)=0$.
\end{lemma}
\bp
Let $\tilde{\theta}=\theta_{\EXT} \times \theta_{\EXT}$ act on $\tilde{\EXT}^2$, and let 
$\widetilde{\EXT^2}=\frac{\tilde{\EXT}^2}{x\sim \tilde{\theta}(x)}$ denote the quotient of $\tilde{\EXT}^2$ by the induced $\ZZ$-action. 
Recall that $\mbox{Int}(\tilde{C}_2(\EXT))$ is the quotient of $\tilde{\EXT}^2\setminus (p_{\EXT} \times p_{\EXT})^{-1}(\mbox{diag}(\EXT^2))$  by the same action of $\ZZ$.
According to \cite[Proposition 12.11]{lesbetaone},
$H_4(\widetilde{\EXT^2};\RR)=0$. Therefore any $4$-cycle of $\mbox{Int}(\tilde{C}_2(\EXT))$ bounds a $5$--chain in the interior of $\widetilde{\EXT^2}$. Such a $5$-chain generically intersects the quotient of $(p_{\EXT} \times p_{\EXT})^{-1}(\mbox{diag}(\EXT^2))$ (that is $\ZZ$ copies of $\mbox{diag}(\EXT^2)$) as a $2$-cycle. Since $H_2(\EXT;\QQ)=0$, this $2$-cycle may be assumed to be empty. Therefore $H_4(\mbox{Int}(\tilde{C}_2(\EXT));\RR)=0$ and since $\mbox{Int}(\tilde{C}_2(\EXT))$ is homotopy equivalent to $\tilde{C}_2(\EXT)$, the lemma follows.
\eop

\begin{theorem}
\label{thmdefzform}
For any $(\Omega_i)_i$ that satisfies the hypotheses of Lemma~\ref{lemOmegai}, $\tilde{Z}(K,\tau)=(I_n((\Omega_i)_i))_{n \in \NN}$ and $\tilde{z}(K,\tau)=(I^c_n((\Omega_i)_i))_{n \in \NN}$. 
\end{theorem}
\bp
In order to prove this theorem, we prove that $I_n((\Omega_i)_i)$ is independent of the specific $(\Omega_i)_i$, and then we choose $\Omega_i$ and $\deleqprop_i$ Poincar\'e dual to each other such that $I_{\Gamma}((\Omega_i)_i)=I_{\Gamma}(\{\eqprop_i\})$.

In order to prove that $I_n((\Omega_i)_i)$ is independent of the specific $(\Omega_i)_i$, it is enough to prove that it does not vary when
$\Omega_1$ is changed to some similar $\Omega^{\prime}_1$
associated with $\omega^{\prime}_1=\omega_1 +d\eta_1$ for a one-form $\eta_1$ supported in $S^2_H$,
and $\omega^{\prime}_{1,D}=\omega_{1,D}+d\eta_{1,D}$ for a one-form $\eta_{1,D}$ compactly supported in $S^2 \setminus \{N\}$.
Extend $\pi(\tau)$ and $p(\sK_1)$ in a neighborhood of $\left(\partial \tilde{C}_2(\EXT_{[0,2[})\cup (\tilde{C}_2(\EXT) \setminus \tilde{C}_2(\EXT_{[0,2[})\right)$, and consider a map $\chi$ valued in $[0,1]$ on this neighborhood that maps a smaller neighborhood to $1$ and the complement of an intermediate neighborhood to $0$.
Set $$\Omega_{1,1}=\Omega_1+d\left(\delta(K)(\thetap^{\ast})\pi(\tau)^{\ast}(\chi\eta_1)\right) +d\left( (\delta(K)J_{\Delta(K)})(\thetap^{\ast})p(\sK_1)^{\ast}(\chi\eta_{1,D})\right)$$
Then $(\Omega^{\prime}_1-\Omega_{1,1})$ is a closed $2$-form on $\tilde{C}_2(\EXT)$ with compact support that vanishes on $$\left(\partial \tilde{C}_2(\EXT_{[0,2[})\cup (\tilde{C}_2(\EXT) \setminus \tilde{C}_2(\EXT_{[0,2[})\right).$$ 
Thus its cohomology class belongs to $H^2_c(\tilde{C}_2(\EXT_{[0,2 -\varepsilon[}));\RR)$ that is isomorphic to $H_4(\tilde{C}_2(\EXT_{[0,2 -\varepsilon[});\RR)$ that is trivial according to Lemma~\ref{lemhfourtc}.
Then $(\Omega^{\prime}_1-\Omega_{1,1})$ has a primitive $\zeta_{1,1}$ with compact support in $\tilde{C}_2(\EXT_{[0,2 -\varepsilon[})$
and $$\zeta_1=\zeta_{1,1}+\delta(K)(\thetap^{\ast})\pi(\tau)^{\ast}(\chi\eta_1) + (\delta(K)J_{\Delta(K)})(\thetap^{\ast})p(\sK_1)^{\ast}(\chi\eta_{1,D})$$
is a primitive of $(\Omega^{\prime}_1-\Omega_1)$.
Let $$\Omega_{\Gamma}(\zeta_1,(\Omega_i)_i)=\tilde{p}(\Gamma,1)^{\ast}(\zeta_1) \wedge \bigwedge_{i=2}^{3n}\tilde{p}(\Gamma,i)^{\ast}(\Omega_i).$$
Then, by the Stokes theorem, the variation of $I_{\Gamma,\rho}((\Omega_i)_i)=\int_{\tilde{C}_{2n}(\EXT,\Gamma,\rho)}\Omega_{\Gamma}((\Omega_i)_i) \in \RR$ when $\Omega_1$ is changed to $\Omega^{\prime}_1$
is $$\int_{\partial \tilde{C}_{2n}(\EXT,\Gamma,\rho)}\Omega_{\Gamma}(\zeta_1,(\Omega_i)_i).$$
This integral is a sum of integrals over the codimension $1$ faces of $\tilde{C}_{2n}(\EXT,\Gamma,\rho)$ of the $(2n-1)$-form $\Omega_{\Gamma}(\zeta_1,(\Omega_i)_i)$.
\begin{lemma}
\label{lempfbordext}
The $(2n-1)$-form $\Omega_{\Gamma}(\zeta_1,(\Omega_i)_i)$ vanishes on the faces where one point is on $\partial \EXT$.
\end{lemma}
\bp The proof is dual to the proof of \cite[Lemma 2.11]{lesbonn}, but we write it in more details below. Recall that $\varepsilon \in ]0,\frac{1}{2n}[$ is involved in the definition of $\pi(\tau)$ in the beginning of Section~\ref{secaltdefzform} with notation of Subsection~\ref{subcyc}. This definition is valid for a fixed $n$. 
Let $v(\ell)$ be the vertex that is mapped to $m_{\ell} \in \partial \EXT$.
With a configuration $c$ of the face that maps $v(i)$ to $m_i \in  \EXT$, associate a subgraph $G$ of $\Gamma$ with edge-orientations
independent from the initial edge-orientations of $\Gamma$, so that $G$ satisfies the following properties:
\begin{itemize}
\item $G$ contains $v(\ell)$, 
\item there is an increasing (with respect to the new edge-orientation) path from $v(\ell)$ to every other vertex of $G$
\item if an edge $e(a)$ of $G$ goes from a vertex $v(i)$ of $G$ to another such $v(j)$ then $r(m_j) >r(m_i)-\frac{4n}{2n-1}\varepsilon$,
\item $G$ contains all the edges of $\Gamma$ between two elements of $G$,
\item if an edge $e(a)$ of $\Gamma$ contains a vertex $v(i)$ of $G$ and a vertex $v(j)$ of $\Gamma \setminus G$, then $r(m_j) <r(m_i)-2\varepsilon$.
\end{itemize}
Such a graph $G$ can be constructed inductively, starting with $v(\ell)$: For each constructed vertex $v(i)$ with less than $3$ incoming edges, successively add each other edge of $\Gamma$ adjacent to $v(i)$ whose other end $v(j)$ satisfies $r(m_j) >r(m_i)-\frac{4n}{2n-1}\varepsilon$, together with the vertex $v(j)$, and orient such an edge from $v(i)$ to $v(j)$. Add the edges between two elements of $G$. The graph $G$ is not uniquely defined from $c$. But any such graph $G$ defines an open subset $F(v(\ell);G)$ of our face $F(v(\ell))$ that is covered by these open sets. It is enough to prove that the restriction of $\Omega_{\Gamma}(\zeta_1,(\Omega_i)_i)$ to
$F(v(\ell);G)$ vanishes. 

Let $E(G)$ denote the set of edges of $G$, let $E_a(G)$ denote the set of edges of $\Gamma$ adjacent to $G$ (i.e. that have one end in $G$), and let $V(G)$ be the set of vertices of $G$.
The form $\Omega_{\Gamma}(\zeta_1,(\Omega_i)_i)$ reads $\Omega_G \wedge \Omega_{o(G)}$ where $\Omega_G$ is the exterior product of the forms associated with the edges of $E(G) \cup E_a(G)$ that read $\pi(\tau)^{\ast}(\omega_i)$ or $\pi(\tau)^{\ast}(\eta_1)$, and $\Omega_{o(G)}$ is the exterior product of the forms associated with the other edges of $\Gamma$.
Write $m_i \in D_{[2,4]} \times S^1$ as $(p_D(m_i) \in D_{[2,4]},\funcE(m_i))$.
For an edge $a$ adjacent to $G$ between $v(i) \in G$ and $v(j) \notin G$, $r(m_j) \leq r(m_i)-2\varepsilon$ and $\pi(\tau)$ factors through
$$(m_i,m_j) \rightarrow (p_D(m_i), \functE(m_i)-\functE(m_j)).$$
(Locally, we identify the $m_i$ with their lifts.)
Write $p_{S,a}(c)=m_i$ and $q_{S,a}(c)=\functE(m_j)$.
The degree $\mbox{deg}(\Omega_G)$ of $\Omega_G$ is $2(\sharp E(G) + \sharp E_a(G))$ or $(2(\sharp E(G) + \sharp E_a(G))-1)$.
On the other hand, because of the above factorisation, if $\Omega_G \neq 0$,
$$\mbox{deg}(\Omega_G) \leq \mbox{dim}\left(\mbox{Im}\left(\prod_{a \in E(G) \cup E_a(G) }\pi(\tau) \circ \tilde{p}(\Gamma,a)\right)\right)\leq  \mbox{dim}\left(\mbox{Im}(P_G)\right)$$
where $P_G=\prod_{a \in E_a(G)}q_{S,a} \times \prod_{a \in E_a(G)}p_{S,a} \times \prod_{a \in E(G)} \pi(\tau) \circ \tilde{p}(\Gamma,a)$.
Now, $P_G$ only depends on the positions of the vertices of $G$ that make a $(3\sharp V(G)-1)$-manifold (since $m_{\ell} \in \partial \EXT$), and on the ''vertical coordinates $q_{S,a}$``. Since a global vertical
translation of these positions does not affect $P_G$, $P_G$ factors through a $(3\sharp V(G)-2 +\sharp E_a(G))$-manifold, and if $\Omega_G \neq 0$,
$$\mbox{deg}(\Omega_G) \leq \sharp E_a(G) + 3\sharp V(G)-2.$$
A count of half-edges shows that $3\sharp V(G)=\sharp E_a(G) +2\sharp E(G)$, and that $\mbox{deg}(\Omega_G) \leq 2(\sharp E(G) + \sharp E_a(G))-2$. This shows that $\Omega_G=0$. Therefore $\Omega_{\Gamma}(\zeta_1,(\Omega_i)_i)$ vanishes on $F(v(\ell);G)$.
\eop

The other codimension one faces correspond to simple collapses of 
subgraphs $\Gamma_I$ of $\Gamma$ made of the vertices numbered in a set $I$ with cardinality at least $2$ and the edges between two of them.
The projection in $C_{2n}(M)$ of a configuration of such a face $\tilde{C}_{2n}(\EXT,\Gamma,\rho;I)$ maps $\Gamma_I$ to a point $m(\Gamma_I)$ and injects the set $(V(\Gamma) \setminus V(\Gamma_I))$ of other vertices of $\Gamma$ to $\EXT \setminus \{m(\Gamma_I)\}$. The face $\tilde{C}_{2n}(\EXT,\Gamma,\rho;I)$ fibers over the complement of all diagonals in $\left( \mathring{\EXT} \times \mathring{\EXT}^{(V(\Gamma) \setminus V(\Gamma_I))} \right)$, and the fiber over 
$(m(\Gamma_I), (m_j)_j)$ is a covering of the quotient of $(T_{m(\Gamma_I)}\EXT)^{V(\Gamma_I)}$ by the translations and the dilations (homotheties with positive ratio).

When the subgraph $\Gamma_I$ is disconnected, the product of the $\tilde{p}(\Gamma,i)$ factors through the effect of translating a connected component of $\Gamma_I$ without moving the other ones; therefore the $(2n-1)$-form $\Omega_{\Gamma}(\zeta_1,(\Omega_i)_i)$ pulls back through a space of dimension strictly less than $(2n-1)$ on $\tilde{C}_{2n}(\EXT,\Gamma,\rho;I)$ and it vanishes there.

If $\Gamma_I$ is connected, if $\Gamma_I$ is not an edge, and if $\Gamma_I$ has a univalent vertex, the product of the $\tilde{p}(\Gamma,i)$ factors through the effect of translating this univalent vertex in the direction of its edge on $\Gamma_i$, and the form $\Omega_{\Gamma}(\zeta_1,(\Omega_i)_i)$ vanishes on $\tilde{C}_{2n}(\EXT,\Gamma,\rho;I)$.

The faces corresponding to edge collapses are treated by the so-called IHX or Jacobi identification that shows that they do not yield any variation, thanks to the IHX relation, see \cite[Lemma 2.21]{lesconst}. Here, we furthermore note that on an edge collapse face, the form associated with the edge is $\delta(K)(\thetap^{\ast}) \pi(\tau)^{\ast}(\omega_i)$ since $p(\sK_i)^{\ast}(\omega_{i,D})$ does not see the direction of the infinitesimal edge. Then the edge will be first beaded by $\delta(K)(t)$ and next by $1$ after the division. Therefore, we can directly pretend that the edge is beaded by $1$.

The faces where $\Gamma_I$ has at least one vertex of valency $2$ is treated by the Bott and Taubes parallelogramm identification (that works as in \cite[Section 2.6]{lesbonn}, see also \cite[Lemma 2.20]{lesconst}).

The faces where $\Gamma_I$ is a component of $\Gamma$ won't contribute because, via the parallelization $\tau$, the product of the $\pi(\tau) \circ \tilde{p}(\Gamma,j)$ for the edges of $\Gamma_I$ does not depend on the position of $m(\Gamma_I)$, while $p(\sK_i) \circ
\tilde{p}(\Gamma,i)$ only depends on a projection to $S^2$ of the image of $m(\Gamma_I)$.

Now, we have proved that $I_n((\Omega_i)_i)$ is independent of the $\Omega_i$ satisfying the above conditions.
We can fix $\cvarM_i$, let the supports of the $\omega_i$ be small enough disks around $\cvarM_i$, let the supports of $\omega_{i,D}$
also be small enough, and choose the $\Omega_i$ Poincar\'e dual to $\deleqprop_i$ such that their supports are small enough neighborhoods of the $\deleqprop_i$ so that the support of $\Omega_{\Gamma}((\Omega_i)_i)$ is concentrated around the intersection points of the $\tilde{p}(\Gamma,i)^{-1}(\deleqprop_i)$. Then $I_n((\Omega_i)_i)=I_n(\{\eqprop_i(\tau)\})=\tilde{Z}_n(K,\tau)$ and we are done.
\eop

\begin{remark}
Let $\iota$ be the involution of $\tilde{C}_2(\EXT)$ that lifts the exchange of two points in a pair, and that preserves $\{0\} \times U\EXT$, $\iota \circ \thetap = \thetap^{-1} \circ \iota$. Assume $\delta(t)=\delta(t^{-1})$. Set $\Omega =\frac12(\Omega_1-\iota^{\ast}(\Omega_1))$. Then $\Omega$ satisfies the same properties as the $\Omega_i$, and it satisfies the additional property that $\Omega=-\iota^{\ast}(\Omega)$. Choose $\Omega_i=\Omega$ for all $i$, then for any $\Gamma \in \CS_n^u$, $I_{\Gamma}((\Omega_i)_i)$ is independent of the numberings of the edges and vertices and of the orientation of the edges. Therefore, $I_n((\Omega_i)_i)$ reads as a sum over the isomorphism classes of unoriented graphs $\Gamma$ of $\frac{I_{\Gamma}((\Omega_i)_i)}{\sharp \mbox{Aut}(\Gamma)}$ where $\sharp \mbox{Aut}(\Gamma)$ is the number of automorphisms of $\Gamma$ (permutations of the half-edges that respect the graph structure).
\end{remark}

\section{Pseudo-parallelizations and associated notions}
\setcounter{equation}{0}
\label{secpseudotriv}

Sections~\ref{secdefpseudointform} and \ref{secdefpseudoint} will present the more flexible definition of $\tilde{Z}$ that is needed to prove Theorem~\ref{thmflag} in general.
In order to give this more flexible definition, we first recall the notion of {\em pseudo-parallelization\/} (or pseudo-trivialization) from \cite[Section 4.3 and 4.2]{lessumgen}. See also \cite[Section~10]{lesbetaone}.

\begin{definition}
\label{defpseudotriv}
A {\em pseudo-parallelization\/} $\tilde{\tau}=(N(c);\tau_e,\tau_b)$ of a compact $3$-manifold $\Aman$ with possible boundary consists of 
\begin{itemize}
\item a framed link $c$ of the interior of $\Aman$ equipped with a neighborhood $N(c)=[a,b] \times c \times [-1,1]$,
\item a parallelization $\tau_e$ of 
$\Aman$ outside $N(c)$,
\item a parallelization $\tau_b \colon N(c) \times \RR^3 \rightarrow TN(c)$ of $N(c)$ such that
$$\tau_b=\left\{\begin{array}{ll} \tau_e & \mbox{over}\; \partial([a,b] \times c \times [-1,1])\setminus (\{a\} \times c \times [-1,1]) \\
 \tau_e \circ  \CT_c& \mbox{over}\; \{a\} \times c \times [-1,1].
\end{array} \right.$$
\end{itemize}
where 
$$\CT_c(t,\gamma \in c,u \in [-1,1];\cvarM\in \RR^3)=(t,\gamma,u,\rho_{\alpha(u)}(\cvarM))$$
where 
$\rho_{\alpha(u)}=\rho(\alpha(u),(0,0,1))$ denotes the rotation of $\RR^3$ 
with axis directed by $(0,0,1)$ and with angle $\alpha(u)$, and $\alpha$ is a smooth map from $[-1,1]$ to $[0,2 \pi]$ 
that maps $[-1,-1+\varepsilon]$ to $0$,
that increases from $0$ to $2\pi$ on $[-1+\varepsilon, 1-\varepsilon]$, and such that $\alpha(-u)+\alpha(u)=2\pi$ for any $u \in [-1,1]$.
\end{definition}

\begin{definition}
\label{defpseudotrivpone}[Trivialisation $\tilde{\tau}_{\CC}$ of $T\Aman \otimes_{\RR} \CC$]
Let $F_U$ be a smooth map such that
$$\begin{array}{llll}F_U:&[a,b] \times [-1,1] &\longrightarrow &SU(3)\\
& (t,u) & \mapsto & \left\{\begin{array}{ll}\mbox{Identity}&
 \mbox{if}\; |u|>1-\varepsilon\\
\rho_{\alpha(u)} &  \mbox{if}\; t<a+\varepsilon\\
\mbox{Identity} &  \mbox{if}\; t>b-\varepsilon.\end{array}\right.\end{array}$$
$F_U$ extends to $[a,b] \times [-1,1]$  because $\pi_1(SU(3))$ is trivial.
Define the trivialization $\tilde{\tau}_{\CC}$ of $T\Aman \otimes_{\RR} \CC$ as follows.
\begin{itemize}
\item On $(\Aman\setminus N(c)) \times \CC^3$, $\tilde{\tau}_{\CC} =\tau_e \otimes 1_{\CC}$,
\item Over $[a,b] \times c \times [-1,1]$, 
$\tilde{\tau}_{\CC}(t,\gamma,u;\cvarM) =\tau_b(t,\gamma,u;F_U(t,u)^{-1}(\cvarM))$.
\end{itemize}
Since $\pi_2(SU(3))$ is trivial, the homotopy class of $\tilde{\tau}_{\CC}$
is well-defined.
\end{definition}

\begin{definition}[Pseudo-sections $s_{\tilde{\tau}}(.;\cvarM)$]
\label{defpseudosec}
 Let $\varepsilon>0$ be a small positive number, recall the map $\alpha$ from Definition~\ref{defpseudotriv}, and define a smooth map 
$$\begin{array}{llll}F:&[a,b] \times [-1,1] &\longrightarrow &SO(3)\\
& (t,u) & \mapsto & \left\{\begin{array}{ll}\mbox{Identity}&
 \mbox{if}\; |u|>1-\varepsilon\\
\rho_{\alpha(u)} &  \mbox{if}\; t<a+\varepsilon\\
\rho_{-\alpha(u)} &  \mbox{if}\; t>b-\varepsilon .\end{array}\right.\end{array}$$

The map $F$ extends to $[a,b] \times [-1,1]$  because its restriction to the boundary
is trivial in $\pi_1(SO(3))$.

Let ${F}(c,\tau_b)$ be defined on $UN(c) \stackrel{\tau_b}{=}   [a,b] \times   c \times [-1,1] \times S^2$ as follows
$$\begin{array}{llll} {F}(c,\tau_b): &  [a,b] \times   c \times [-1,1] \times S^2 & \longrightarrow & S^2\\
&(t,\gamma,u; \svarM) & \mapsto &  F(t,u)(\svarM).\end{array}$$

Let $\cvarM \in S^2$ and let $S^1(\cvarM)$ be the circle (or point) in $S^2$ that lies in a plane orthogonal to the axis generated by $(0,0,1)$ and that contains $\cvarM$.
There is a $2$-dimensional chain $C_2(\cvarM)$ in $[-1,1] \times S^1(\cvarM)$ whose
boundary is 
$\{\left(u,\rho_{-\alpha(u)}(\cvarM)\right), u \in[-1,1]\}
+\{\left(u,\rho_{\alpha(u)}(\cvarM)\right), u \in[-1,1]\} - 2[-1,1]\times \{\cvarM\}.$
Then $s_{\tilde{\tau}}(\Aman;\cvarM)$ is the following $3$--cycle of $(U\Aman, U\Aman_{|\partial \Aman})$
$$s_{\tilde{\tau}}(\Aman;\cvarM)=s_{\tau_e}(\Aman\setminus \mathring{N}(c);\cvarM) + 
\frac{s_{\tau_b \circ \CT_c^{-1}}(N(c);\cvarM) + F(c,\tau_b)^{-1}(\cvarM) + \{b\} \times c \times C_2(\cvarM)}{2}
$$
where $\tau_b$ and $\tau_e$ identify $U\Aman_{|\{b\} \times c \times [-1,1]}$ with $\{b\} \times c \times [-1,1] \times S^2$ in the same way.
When $\Sigma$ is a $2$--chain that intersects $N(c)$ along sections $N_{\gamma}(c)=[a,b]\times\{\gamma\}\times [-1,1],$ $$s_{\tilde{\tau}}(\Sigma;\cvarM)=s_{\tilde{\tau}}(\Aman;\cvarM)\cap U\Aman_{|\Sigma}$$ so that $$s_{\tilde{\tau}}(N_{\gamma}(c);\cvarM)= 
\frac{s_{\tau_b \circ \CT_c^{-1}}(N_{\gamma}(c);\cvarM) + F(c,\tau_b)^{-1}(\cvarM)\cap U\Aman_{|N_{\gamma}(c)} - \{b\} \times \{\gamma\} \times C_2(\cvarM)}{2}.
$$
We also introduce small deformations of these sections associated with $\varepsilon_i$ such that $0<\varepsilon_i < \varepsilon$ as follows.
Let $N(c,\varepsilon_i)=[a,b-\varepsilon_i] \times c \times [-1,1] \times S^2$, and let ${F}(c,\tau_b,\varepsilon_i)$ be the restriction of ${F}(c,\tau_b)$ to $N(c,\varepsilon_i)$. Then $s_{\tilde{\tau}}(\Aman;\cvarM,\varepsilon_i)$ is the following $3$--cycle of $(U\Aman, U\Aman_{|\partial \Aman})$
$$s_{\tilde{\tau}}(\Aman;\cvarM,\varepsilon_i)=$$ $$s_{\tau_e}(\Aman\setminus \mathring{N}(c,\varepsilon_i);\cvarM) + 
\frac{s_{\tau_b \circ \CT_c^{-1}}(N(c,\varepsilon_i);\cvarM) + F(c,\tau_b,\varepsilon_i)^{-1}(\cvarM) + \{b-\varepsilon_i\} \times c \times C_2(\cvarM)}{2}.
$$
\end{definition}

\begin{lemma}
\label{lemdiffz}
Let $A$ be a $\QQ$-handlebody equipped with two pseudo-parallelizations $\tilde{\tau}_0$ and $\tilde{\tau}_1$ that coincide with the same genuine parallelization near the boundary of $A$. Let $\cvarM \in S^2$. Let $\eta \in ]0,\varepsilon[$. 
There exists a rational $4$-chain $H(\tilde{\tau}_0,\tilde{\tau}_1,\cvarM,\eta)$ in $UA$ such that
$$\partial H(\tilde{\tau}_0,\tilde{\tau}_1,\cvarM,\eta)=s_{\tilde{\tau}_1}(A;\cvarM,\eta) - s_{\tilde{\tau}_0}(A;\cvarM,\eta).$$
\end{lemma}
\bp
Let us show that $C=s_{\tilde{\tau}_1}(A;\cvarM,\eta) - s_{\tilde{\tau}_0}(A;\cvarM,\eta)$ vanishes in 
$$H_3(UA;\QQ)=H_1(A;\QQ) \otimes H_2(S^2;\QQ).$$
In order to do so, we prove that the algebraic intersections of $C$ with $s_{\tilde{\tau}_0}(S;\svarM)$ vanish for $2$--cycles $S$ of $(A,\partial A)$ that generate $H_2(A,\partial A)$, for $\svarM \in S^2 \setminus S^1(\cvarM)$.
Of course, $s_{\tilde{\tau}_0}(\partial A;\cvarM)$ does not intersect $s_{\tilde{\tau}_0}(S;\svarM)$. According to Lemma~10.5 in \cite{lesbetaone}, $s_{\tilde{\tau}_0}(S;\cvarM)$ does not intersect $s_{\tilde{\tau}_0}(S;\svarM)$ algebraically in $UA_{|S}$, so that $s_{\tilde{\tau}_0}(A;\cvarM)$  does not intersect $s_{\tilde{\tau}_0}(S;\svarM)$, either. Since Lemma~10.7 in \cite{lesbetaone} guarantees that $(s_{\tilde{\tau}_0}(S;\cvarM)-s_{\tilde{\tau}_1}(S;\cvarM))$ bounds in $UA_{|S}$, $s_{\tilde{\tau}_1}(S;\cvarM)$ does not intersect $s_{\tilde{\tau}_0}(S;\svarM)$ algebraically in $UA_{|S}$, so that $s_{\tilde{\tau}_1}(A;\cvarM)$ does not intersect $s_{\tilde{\tau}_0}(S;\svarM)$, algebraically. This proves the existence of the $4$--chain $H(\tilde{\tau}_0,\tilde{\tau}_1,\cvarM,\eta)$.
\eop

\begin{definition}[$2$--forms associated with pseudo-parallelizations on $U\Aman$ with $S^2_H$-support]
\label{defpseudoform}
Recall the part $S^2_H$ of $S^2$ introduced in Proposition~\ref{propdefbord}. The axis $(0,0,1)$ of $\rho_{\alpha(.)}$ is vertical.
Let $\omega_s$ be a $2$--form compactly supported in $S^2_H$, invariant under the rotations around the vertical axis, such that $\int_{S^2}\omega_s=1$.
Let $\omega_i$ be a $2$--form compactly supported in $S^2_H$, such that $\int_{S^2}\omega_i=1$.
Let $\eta_{i,s,1}$ be a $1$-form with compact support in $S^2_H$ such that $\omega_i = \omega_s +d\eta_{i,s,1}.$ Let $\varepsilon_i \in ]0,\frac{\varepsilon}2[$ and let $k$ be a large integer. Let $\pi(\tau_b)$ denote the projection from $UN(c)$ to $S^2$ induced by $\tau_b$. $$\pi(\tau_b)(\tau_b(t,\gamma,u;\cvarM \in S^2))=\cvarM .$$
Let $\eta_{i,s}$ be a $1$-form with compact support on $[b-\varepsilon_i-2\varepsilon_i^k,b-\varepsilon_i+2\varepsilon_i^k]\times c\times [-1,1] \times S^2_H$ that pulls back through $p_{[a,b]} \times p_{[-1,1]} \times \pi(\tau_b)$ and such that
$$\eta_{i,s}=\left\{\begin{array}{ll}
\frac{\pi(\tau_b \circ \CT_c^{-1})^{\ast}(\eta_{i,s,1}) +\pi(\tau_b \circ \CT_c)^{\ast}(\eta_{i,s,1})}{2} & \mbox{on}\;[b-\varepsilon_i-2\varepsilon_i^k,b-\varepsilon_i-\varepsilon_i^k]\times c\times [-1,1] \times S^2_H\\
\pi(\tau_b)^{\ast}(\eta_{i,s,1})& \mbox{on}\;[b-\varepsilon_i+\varepsilon_i^k,b-\varepsilon_i+2\varepsilon_i^k]\times c\times [-1,1] \times S^2_H
\\ &\mbox{and on}\;[b-\varepsilon_i-2\varepsilon_i^k,b-\varepsilon_i+2\varepsilon_i^k]\times c\times \partial [-1,1] \times S^2_H.
\end{array}\right.$$ 

Define
$$\omega(\tilde{\tau},\omega_i,k,\varepsilon_i,\eta_{i,s})=\left\{\begin{array}{ll}
\pi(\tau_e)^{\ast}(\omega_i) & \mbox{on}\;U\left(\Aman \setminus \left([a,b-\varepsilon_i+\varepsilon_i^k] \times c \times [-1,1]\right)\right)\\
\frac{\pi(\tau_b \circ \CT_c^{-1})^{\ast}(\omega_i) + F(c,\tau_b)^{\ast}(\omega_i)}{2}& \mbox{on}\;U([a,b-\varepsilon_i-\varepsilon_i^k] \times c \times [-1,1])
\\
\pi(\tau_b)^{\ast}(\omega_s) + d \eta_{i,s}& \mbox{on}\;U([b-\varepsilon_i-\varepsilon_i^k,b-\varepsilon_i+\varepsilon_i^k]  \times c \times ]-1,1[)
\end{array}\right.$$with the notation of Definition~\ref{defpseudosec}.
\end{definition}

\begin{lemma}
Definition~\ref{defpseudoform} of  $\omega(\tilde{\tau},\omega_i,k,\varepsilon_i,\eta_{i,s})$ is consistent.
\end{lemma}
\bp
There exists a function $g_s \colon ]-1/50,1/50[ \rightarrow \RR$ with compact support such that $\omega_s = (g_s\circ p_v) d\theta \wedge d p_v$ where $p_v$ is the vertical coordinate in $\CC \times \RR$ and $\theta$ is the angular coordinate around the vertical axis (the argument in the horizontal $\CC$).
Recall $\CT_c (t,\gamma \in c,u \in [-1,1];\cvarM\in S^2)=(t,\gamma,u,\rho_{\alpha(u)}(\cvarM))$. Thus
$\tau_b \circ \CT^{\kappa}_c(t,\gamma,u;\cvarM)=\tau_b(t,\gamma,u;\rho_{\kappa\alpha(u)}(\cvarM))$ and
$$\pi(\tau_b \circ \CT^{\kappa}_c)(\tau_b(t,\gamma,u;\svarM))=\rho_{-\kappa\alpha(u)}(\svarM).$$
so that $F(c,\tau_b)^{\ast}(\omega_i)=\pi(\tau_b \circ \CT_c)^{\ast}(\omega_i)$ on $[b-\varepsilon_i-2\varepsilon_i^k,b-\varepsilon_i-\varepsilon_i^k] \times c \times [-1,1] \times S^2$ and
$$\pi(\tau_b \circ \CT^{\kappa}_c)^{\ast}(\omega_s)=\pi(\tau_b)^{\ast}(\omega_s)-\kappa(g_s\circ p_v)\alpha^{\prime}(u) du \wedge d p_v.$$
Furthermore
$\left(\pi(\tau_b \circ \CT_c^{-1})^{\ast}(\omega_s) + F(c,\tau_b)^{\ast}(\omega_s)\right)$ reads $2\pi(\tau_b)^{\ast}(\omega_s)$ on $[b-\varepsilon_i-2\varepsilon_i^k,b-\varepsilon_i-\varepsilon_i^k] \times c \times [-1,1] \times S^2$ and extends as $2\pi(\tau_b)^{\ast}(\omega_s)$ on $U([b-\varepsilon_i-\varepsilon_i^k,b-\varepsilon_i+\varepsilon_i^k]  \times c \times ]-1,1[)$.
This shows that the definition is consistent when $\omega_i=\omega_s$, $\eta_{i,s,1}=0$ and $\eta_{i,s}=0$.
The general case follows easily.
\eop

\begin{definition}[$2$--forms associated with pseudo-parallelizations on $U\EXT$]
\label{defpseudoformgen}
Definition~\ref{defpseudoform} can be generalized starting with forms $\omega_i$ with arbitrary support in $S^2$ such that $\int_{S^2}\omega_i=1$ by removing all the constraints involving $S^2_H$ from Definition~\ref{defpseudoform}.
In particular, $\omega_i$ can be the homogeneous volume $1$ form of $S^2$ that was already treated in \cite[Section 4.3]{lessumgen}, and the definition here coincides with the definition of \cite[Section 4.3]{lessumgen}.
\end{definition}

Finally, recall Lemma~10.2 in \cite{lesbetaone}.
\begin{lemma}
\label{lem102lesbetaone}
 Let $A$ be a $\QQ$-handlebody and let $\tau$ be a parallelization of $A$ defined on a collar $[-1,0]\times \partial A$ of $\partial A$.
Then there is a pseudo-parallelization of $A$ that extends the restriction of $\tau$ to $[-1,0]\times \partial A$.
\end{lemma}

\section{On the anomaly $\ansothree$}
\setcounter{equation}{0}
\label{secanomaly}

\begin{lemma}
\label{lemdiffzbis}
Let $A$ be a $\QQ$-handlebody equipped with two pseudo-parallelizations $\tilde{\tau}$ and $\tilde{\tau}_1$ that coincide near the boundary of $A$. 
For $i=1, \dots, 3n$, let $\omega_i$ be a $2$--form on $S^2$ such that $\int_{S^2}\omega_i=1$, let $\omega_i(\tilde{\tau})=\omega(\tilde{\tau},\omega_i,k,\varepsilon_i,\eta_{i,s})$ be a form on $UA$ associated with $\omega_i$ and $\tilde{\tau}$ by Definition~\ref{defpseudoformgen}, and let $\omega_i(\tilde{\tau}_1)$ be a form similarly associated with $\omega_i$ and $\tilde{\tau}_1$. 
For $i=1, \dots, 3n$, there exists a $1$--form $\eta_i=\eta_i(\omega_i,\tilde{\tau},\tilde{\tau}_1)$ on $UA$ with compact support in the interior of $A$ such that $d\eta_i=\omega_i(\tilde{\tau}_1)-\omega_i(\tilde{\tau})$.
\end{lemma}
\bp First note that the cohomology class of $(\omega_i(\tilde{\tau}_1)-\omega_i(\tilde{\tau}))$ is independent of $(\omega_i,k,\eta_{i,s})$. In particular, $\omega_i$ can be assumed to be supported in a very small neighborhood of some $\cvarM \in S^2_H$ so that $\omega_i(\tilde{\tau})$ is supported in a neighborhood of $s_{\tilde{\tau}}(\Aman;\cvarM,\varepsilon_i)$ and locally dual to this chain outside $U([b-\varepsilon_i-\varepsilon_i^k,b-\varepsilon_i+\varepsilon_i^k]  \times c \times ]-1,1[)$. Furthermore, its support intersects $U([b-\varepsilon_i-\varepsilon_i^k,b-\varepsilon_i+\varepsilon_i^k]  \times c \times ]-1,1[)$ in $[b-\varepsilon_i-\varepsilon_i^k,b-\varepsilon_i+\varepsilon_i^k]  \times c \times ]-1,1[ \times S^2_H$.
In order to prove Lemma~\ref{lemdiffzbis} as Lemma~\ref{lemdiffz}, it suffices to prove that $$\int_{s_{\tilde{\tau}}(S;\svarM)}\omega_i(\tilde{\tau}_1)-\omega_i(\tilde{\tau}) =0$$ for any $2$--cycle $S$ of $(A,\partial A)$, for $\svarM \in S^2 \setminus S^2_H$. Since this integral is the algebraic intersection with $s_{\tilde{\tau}_1}(A;\cvarM,\eta) - s_{\tilde{\tau}}(A;\cvarM,\eta)$, it vanishes according to Lemma~\ref{lemdiffz}.
\eop

The preimage under the blow-down map $p_b \colon C_{2n}(\Mbet) \rightarrow \Mbet^{2n}$ of $\Delta_{2n}(\Mbet)=\{(m,m,\dots,m);m \in \Mbet\}$ is a fibered space $p_S \colon S_{2n}(T\Mbet)\rightarrow \Mbet$, its fiber $p_S^{-1}(m)$ is a compactification
$S_{2n}(T_m\Mbet)$ of the space of injections from $\{1,2,\dots, 2n\}$ to $T_m\Mbet$ up to translations and dilations. The space $S_{2n}(T\Mbet)$ is oriented as a codimension $0$ part of
$\partial C_{2n}(\Mbet)$.
When $\Aman$ is a submanifold of $\Mbet$, $S_{2n}(T\Aman)=p_S^{-1}(\Aman) \subset S_{2n}(T\Mbet)$. Note that $S_{2}(T\Aman)$ is the unit tangent bundle $U\Aman$ to $\Aman$.
When a graph $\Gamma$ of $\CS_n$ is given, the maps $p(\Gamma,i)$ introduced in Subsection~\ref{subeqint} restrict to natural maps $p(\Gamma,i) \colon S_{2n}(T\Aman) \rightarrow  U\Aman$.

\begin{proposition}
\label{propxinzform}
Under the hypotheses of Lemma~\ref{lemdiffzbis}, 
let $p_{[0,1]}$ (resp. $p_{UA}$) denote the natural projection from $[0,1] \times UA$
to $[0,1]$ (resp. to $UA$).
For $i=1, \dots, 3n$, there exists a closed $2$--form $\zeta_i$ on $[0,1] \times UA$ that restricts as $p_{UA}^{\ast}(\omega_i(\tilde{\tau}))$ on 
$(\{0\} \times UA) \cup ([0,1] \times UA_{|\partial A})$ and as $p_{UA}^{\ast}(\omega_i(\tilde{\tau}_1))$ on $(\{1\} \times UA)$.
Then $$\frac{1} {2^{3n}(3n)!(2n)!}\sum_{\Gamma \in \CS_n}\int_{[0,1] \times S_{2n}(TA)}\left(\bigwedge_{j=1}^{3n}p(\Gamma,j)^{\ast}(\zeta_j)\right)[\Gamma] \in \CA_n^h(1)
$$ only depends on $A$, $\tilde{\tau}$ and $\tilde{\tau}_1$. It will be denoted by $\ansothree_n(A,\tilde{\tau},\tilde{\tau}_1)$.
\end{proposition}
\bp 
The existence of $\zeta_i$ comes from Lemma~\ref{lemdiffzbis}: If $\chi_i\colon [0,1] \rightarrow [0,1]$ is a map that maps
a neighborhood of $0$ to $0$ and a neighborhood of $1$ to $1$, we can choose $$\zeta_i=d(p_{[0,1]}^{\ast}(\chi_i) p_{UA}^{\ast}(\eta_i))+p_{UA}^{\ast}(\omega_i(\tilde{\tau})).$$ The independence of the choices of the $\zeta_i$ when the $\omega_i$ are fixed can easily be proved like Theorem~\ref{thmdefzform} since $H^2([0,1]\times UA,\partial ([0,1]\times UA);\QQ)=0$. When the $\omega_i$ are not fixed, first use partial primitives of $(\omega^{\prime}_i(\tilde{\tau}_1)-\omega_i(\tilde{\tau}_1))$ and $(\omega^{\prime}_i(\tilde{\tau})-\omega_i(\tilde{\tau}))$
that factor through $(p_{[a,b]} \times p_{[-1,1]} \times \pi(\tau))$ on $UN(c)$ and through $\pi(\tau)$ elsewhere on $UA$.
\eop

\begin{proposition}
\label{propxinzformbis} Under the hypotheses of Lemma~\ref{lemdiffzbis} and with the notation of Proposition~\ref{propxinzform},
$$\ansothree_n(A,\tilde{\tau},\tilde{\tau}_1)=$$
$$\frac{1} {2^{3n}(3n)!(2n)!}\sum_{i=1}^{3n}\sum_{\Gamma \in \CS_n}\int_{S_{2n}(TA)}p(\Gamma,i)^{\ast}(\eta_i) \bigwedge \wedge_{j=1}^{i-1}p(\Gamma,j)^{\ast}(\omega_j(\tilde{\tau}_1)) \bigwedge \wedge_{j=i+1}^{3n}p(\Gamma,j)^{\ast}(\omega_j(\tilde{\tau}))[\Gamma].
$$
\end{proposition}
\bp
Let $\chi_i\colon [0,1] \rightarrow [0,1]$ be a map that maps $[\frac{i}{3n}-\frac{1}{9n} ,1]$ to $1$ and $[0,\frac{i}{3n}-\frac{2}{9n}]$ to $0$, choose $\zeta_i=d(p_{[0,1]}^{\ast}(\chi_i) p_{UA}^{\ast}(\eta_i))+p_{UA}^{\ast}(\omega_i(\tilde{\tau}))$ in Proposition~\ref{propxinzform}. Then a simple calculation allows us to conclude. 
\eop

\begin{definition}
\label{defxin}
 Let $\rho \colon (B^3,\partial B^3) \rightarrow (SO(3),1)$ be the map defined in Lemma~\ref{lempreptrivun} that induces the double covering of $SO(3)$.
Consider the standard parallelization $\tau_B \colon B^3 \times \RR^3 \rightarrow TB^3$ of $B^3$, and the parallelization $\tau_B\circ \psi_{\RR}(\rho)$ where $\psi_{\RR}(\rho)(x,v)=(x,\rho(x)(v))$. With the notation of Proposition~\ref{propxinzform}, set
$$\ansothree_n=\ansothree_n(B^3,\tau_B,\tau_B\circ \psi_{\RR}(\rho)).$$
See \cite[Subsection 1.6]{lesconst} for an alternative definition of $\ansothree_n=-\xi_n$. According to \cite[Proposition 1.10]{lesconst}, $\ansothree_{2k}=0$ for any integer $k$.
\end{definition}

\begin{theorem}
\label{thmxin}
Let $A$ be a $\QQ$-handlebody equipped with two pseudo-parallelizations $\tilde{\tau}$ and $\tilde{\tau}_1$ that coincide near the boundary of $A$. Then $\ansothree_n(A,\tilde{\tau},\tilde{\tau}_1)=\frac{p_1(\tilde{\tau},\tilde{\tau}_1)}{4}\ansothree_n$
\end{theorem}
\bp
If $A$ is a rational homology ball, this is Proposition~4.8 in \cite{lessumgen}, because the variation that is expressed there is
$\ansothree_n(A,\tilde{\tau},\tilde{\tau}_1)(2^{3n}(3n)!(2n)!)$. See \cite[Subsection 2.3]{lesconst}.
Then any $\QQ$-handlebody $A$ can be embedded in such a ball $B$ and the pseudo-parallelizations can be extended in the same way to $B\setminus A$, so that $\ansothree_n(A,\tilde{\tau},\tilde{\tau}_1)=\ansothree_n(B,\tilde{\tau},\tilde{\tau}_1)=\frac{p_1(\tilde{\tau}_{|B},\tilde{\tau}_{1|B})}{4}\ansothree_n=\frac{p_1(\tilde{\tau}_{|A},\tilde{\tau}_{1|A})}{4}\ansothree_n$, according to Proposition~\ref{proppont} (2).
\eop

\section{A more flexible definition of $\tilde{Z}$ with pseudo-parallelizations and differential forms}
\setcounter{equation}{0}
\label{secdefpseudointform}

The definition of $\tilde{z}(K,\tau)$ with differential forms extends to the case when $\tau$ is a pseudo-parallelization that coincides with $\tau_{\pi}$ on $\EXT_{[1,4]}$.

\begin{theorem}
\label{thmdefzformpseudo}
With the notation of Definition~\ref{defpseudoform}, and of the beginning of Section~\ref{secaltdefzform}, let $\tilde{\tau}$ be a pseudo-parallelization that coincides with $\tau_{\pi}$ on $\EXT_{]1/2,4]}$,
there exists a closed $2$-form $\Omega_i$ with compact support on $\tilde{C}_2(\EXT)$ that restricts
as $$\delta(K)(\thetap^{\ast})\pi(\tau_{\pi})^{\ast}(\omega_i) + (\delta(K)J_{\Delta(K)})(\thetap^{\ast})(p(\sK_i)^{\ast}(\omega_{i,D}))$$ on $\left(\partial \tilde{C}_2(\EXT_{]1/2,2[})\cup (\tilde{C}_2(\EXT) \setminus \tilde{C}_2(\EXT_{[0,2[}))\right)$,
and as $\delta(K)(\thetap^{\ast})\omega(\tilde{\tau},\omega_i,k,\varepsilon_i,\eta_{i,s})$ on $U\EXT_{[0,1/2]}$.
Recall the definition of $(I^c_n((\Omega_i)_i))_{n \in \NN}$ in Section~\ref{secaltdefzform}.
Then $(I^c_n((\Omega_i)_i))_{n \in \NN}$ is independent of the chosen $(\Omega_i)_i$, it is an invariant of $(K,\tilde{\tau})$ that is denoted by $\tilde{z}_n(K,\tilde{\tau})$.
Furthermore, if $\tilde{\tau}_1$ is another pseudo-parallelization that coincides with $\tau_{\pi}$ on $\EXT_{[1/2,4]}$, then
$$\tilde{z}_n(K,\tilde{\tau}_1)-\tilde{z}_n(K,\tilde{\tau})=\frac{p_1(\tilde{\tau},\tilde{\tau}_1)}{4}\ansothree_n.$$
\end{theorem}
\bp Let $\tau$ be a parallelization of $\EXT$ that coincides with $\tau_{\pi}$ on $\EXT_{[1/2,4]}$. Let $\Omega_i(\tau)$ be a form that satisfies the hypotheses of Theorem~\ref{thmdefzform} with respect to $\tau$.
 Let $\chi\colon\tilde{C}_2(\EXT) \rightarrow [0,1]$ be a map supported in a neighborhood of $U\EXT_{[0,1/2]}$ that maps a smaller neighborhood of $U\EXT_{[0,1/2]}$ to $1$. Extending the form $\eta_i=\eta_i(\omega_i,\tilde{\tau},\tau)$ of Lemma~\ref{lemdiffzbis} to $\tilde{C}_2(\EXT)$ so that $\Omega_i=\Omega_i(\tau) -\delta(K)(\thetap^{\ast})(d\chi\eta_i)$ satisfies the conditions of the statement shows the existence of $\Omega_i$. Conversely, for any $\Omega_i$ that satisfies the hypotheses of the statement, the form
$\Omega_i(\tau)=\Omega_i+\delta(K)(\thetap^{\ast})(d\chi\eta_i)$ satisfies the hypotheses of
 Theorem~\ref{thmdefzform} with respect to $\tau$.

Let $$\Omega_{\Gamma}(\eta_i,(\Omega_j)_j)=\tilde{p}(\Gamma,i)^{\ast}(\delta(K)(\thetap^{\ast})(\chi\eta_i)) \wedge \bigwedge_{j=1}^{i-1}\tilde{p}(\Gamma,j)^{\ast}(\Omega_j(\tau))\wedge \bigwedge_{j=i+1}^{3n}\tilde{p}(\Gamma,j)^{\ast}(\Omega_j).$$
The variation of $I_{\Gamma,\rho}((\Omega_i)_i)=\int_{\tilde{C}_{2n}(\EXT,\Gamma,\rho)}\Omega_{\Gamma}((\Omega_i)_i) \in \RR$ when $\Omega_1$ is changed to $\Omega_1(\tau)$
is $$\int_{\partial \tilde{C}_{2n}(\EXT,\Gamma,\rho)}\Omega_{\Gamma}(\eta_1,(\Omega_j)_j).$$
This integral is a sum of integrals over the codimension $1$ faces of $\tilde{C}_{2n}(\EXT,\Gamma,\rho)$ of the $(2n-1)$-form $\Omega_{\Gamma}(\eta_1,(\Omega_j)_j)$, where the faces that do not correspond to a total collapse of the connected $\Gamma$ cancel as in the proof of Theorem~\ref{thmdefzform}.

The $\tilde{p}(\Gamma,j)$ project the faces that correspond to a total collapse to $\partial \tilde{C}_{2}(\EXT)$. 
Since  $\eta_1$ is supported near $U\EXT_{[0,1/2]}$ where $\Omega_j$ reads $\delta(K)(\thetap^{\ast})\omega(\tilde{\tau},\omega_i,k,\varepsilon_i,\eta_{i,s})$, and since the beads will be divided by $\delta(t)$, we change $\Omega_j$ to $\Omega_j^s=\omega(\tilde{\tau},\omega_i,k,\varepsilon_i,\eta_{i,s})$ in the following discussion. Then the $\tilde{p}(\Gamma,j)$ project the support of $$\Omega_{\Gamma}^s(\eta_1,(\Omega^s_j)_j)=\tilde{p}(\Gamma,i)^{\ast}(\chi\eta_1) \wedge \bigwedge_{j=2}^{3n}\tilde{p}(\Gamma,j)^{\ast}(\Omega^s_j)$$ on these faces to $U\EXT_{[0,1/2]} \times \ZZ$, and more precisely to $U\EXT_{[0,1/2]} \times \{0\}$ because of the known restriction of $\Omega^s_j$ on the boundary.
Finally, the variation of $I_{\Gamma,\rho}((\Omega_i)_i)$ when $\Omega_1$ is changed to $\Omega_1(\tau)$ is zero if $\rho \neq 0$, and it is 
$$\int_{{S}_{2n}(T\EXT)}\Omega_{\Gamma}^s(\eta_1,(\Omega^s_j)_j),$$ otherwise.
Similarly computing the variations when $\Omega_i$ is changed to $\Omega_i(\tau)$, successively, for $i=2, \dots, 3n,$ shows that
$$(I^c_n((\Omega_i(\tau))_i))_{n \in \NN}-(I^c_n((\Omega_i)_i))_{n \in \NN}=\ansothree_n(\EXT,\tilde{\tau},\tau)$$
with the expression of Proposition~\ref{propxinzformbis}.
This proves that $(I^c_n((\Omega_i)_i))_{n \in \NN}$ is an invariant $z(K,\tilde{\tau})$ of $(K,\tilde{\tau})$ and that 
$z_n(K,\tau)-z_n(K,\tilde{\tau})=\ansothree_n(\EXT,\tilde{\tau},\tau)$. Now, we may apply the same technique to prove $$\tilde{z}_n(K,\tilde{\tau}_1)-\tilde{z}_n(K,\tilde{\tau})=\ansothree_n(\EXT,\tilde{\tau},\tilde{\tau}_1)$$ in general, or use that 
$\ansothree_n(\EXT,\tilde{\tau},\tilde{\tau}_1)=\ansothree_n(\EXT,\tilde{\tau},\tau)-\ansothree_n(\EXT,\tilde{\tau}_1,\tau)$, that is clear by Proposition~\ref{propxinzform}. Conclude with Theorem~\ref{thmxin}.
\eop

\section{The definition of $\tilde{Z}_n(K,.)$ for pseudo-parallelizations with equivariant intersections.}
\setcounter{equation}{0}
\label{secdefpseudoint}

The definition of $\tilde{Z}_n(K,.)$ of Section~\ref{secrecdef} naturally extends to pseudo-parallelizations as follows.
\begin{theorem}
\label{thmdefinvpseudo}
There exist
\begin{itemize}
\item $\varepsilon \in ]0,\frac{1}{2n}[$, $\funcE$ (involved in the definition of the map $\pi$ in Subsection~\ref{subcyc}), $k\neq 0 \in \NN$,
\item $3n$ distinct elements $(\sK_1, \dots, \sK_{3n})$ of $S^1$,
\item $3n$ distinct $\varepsilon_i$ in $]0,\varepsilon[$,
\item a {\em regular\/} $3n$-tuple $(\cvarM_1, \dots,\cvarM_{3n})$  of $(S^2_H)^{3n}$,
\item $3n$ integral chains $C(\cvarM_i,\sK_i,\tilde{\tau})$ of $\tilde{C}_2(\EXT_{[0,2]})$  in {\em general $3n$-position\/} whose boundaries are $$k\delta(K)(\thetap)
\left(\pi_{|\partial \tilde{C}_2(\EXT_{[0,2]}) \setminus \partial \tilde{C}_2(\EXT_{[0,2[})}^{-1}(\cvarM_i) \cup s_{\tilde{\tau}}(\EXT_{[0,2]};\cvarM_i,\varepsilon_i) \cup J_{\Delta}(\thetap) U\EXT_{|{\sK_i \times S^1}}\right),$$ respectively. (See Definition~\ref{defpseudosec} for $s_{\tilde{\tau}}$.)
\end{itemize}
 Under these assumptions, set $$\eqprop_i(\tilde{\tau})=\frac{1}{k\delta(K)(t)}C(\cvarM_i,\sK_i,\tilde{\tau}) + \overline{\pi^{-1}(\cvarM_i)} \subset C_2(\EXT).$$
Then $I^c_n(\{\eqprop_i(\tilde{\tau})\})=\tilde{z}_n(K,\tilde{\tau})$, where $\tilde{z}_n(K,\tilde{\tau})$ has been defined in Theorem~\ref{thmdefzformpseudo}, and $ I_n(\{\eqprop_i(\tilde{\tau})\})=\tilde{Z}_n(K,\tilde{\tau})$
where $$(\tilde{Z}_n(K,\tilde{\tau}))_{n\in \NN}=\exp\left(\sum_{j>0}\tilde{z}_j(K,\tilde{\tau})\right)=\tilde{Z}(K,\tilde{\tau})=\tilde{Z}(K)\exp\left(\frac{p_1(\tau_\pi,\tilde{\tau})}{4}\ansothree\right).$$
\end{theorem}

The proof of this theorem is very similar to the proofs of Theorems~\ref{thminvwithtau} and \ref{thmhighloop}. We only detail the additional different arguments below.

\subsection{Regular points of $(S^2_H)^{3n}$}
\label{subuplegen}
Extend $\pi$ to $\left(\tilde{C}_2(\EXT_{[0,4]}) \setminus \tilde{C}_2(\EXT_{[0,1.5[})\right)$.
The restriction of $p_C \colon \tilde{C}_2(\EXT)\rightarrow C_2(\EXT)$ to $$\left(\tilde{C}_2(\EXT_{[0,4]}) \setminus \tilde{C}_2(\EXT_{[0,1.5[})\right) \cap \pi^{-1}(\overline{S^2_H})$$ is a diffeomorphism $p_H$ to its image $P_H$ in $\left({C}_2(\EXT_{[0,4]}) \setminus {C}_2(\EXT_{[0,1.5[})\right)$.
Define $\overline{\pi}$ as $\pi \circ p_H^{-1}$ on $P_H$, for the map $\pi$ of Subsection~\ref{subcyc}.

Consider a graph $\Gamma \in \CS_n$, its set of vertices $V(\Gamma)$, a coloring $L\colon V(\Gamma) \rightarrow \{\EXT_{[0,1.5]},\EXT_{[1.5,4]}\}$ of its vertices, and a set $B$ of edges of $\Gamma$ that contains all the edges between two vertices of color
$\EXT_{[0,1.5]}$. Set $A=\{1,\dots,3n\} \setminus B$ and 
$${C}(\Gamma,L,B)=p_b^{-1}\left(\EXT_{[0,1.5]}^{L^{-1}(\{\EXT_{[0,1.5]}\})} \times \EXT_{[1.5,4]}^{L^{-1}(\{\EXT_{[1.5,4]}\})}\right) \cap \cap_{a\in A}p(\Gamma, a)^{-1}(P_H) \subset C_{2n}(\EXT)$$
where $p_b\colon C_{2n}(\EXT) \rightarrow \EXT^{2n}$ denotes the blow-down map.
Consider the product $$g(\Gamma,L,B)= \prod_{a\in A}\overline{\pi} \circ p(\Gamma, a)\times \mbox{Identity}(\overline{S^2_H})^B \colon {C}(\Gamma,L,B)\times (\overline{S^2_H})^B \rightarrow (\overline{S^2_H})^A \times (\overline{S^2_H})^B.$$

In general, a {\em regular value\/} of a map $g\colon C \rightarrow Y$ is a point $y$ of $Y$ such that for any point $x$ of $g^{-1}(y)$, the tangent map to $g$ at $x$ is surjective.
Here, our spaces $C$ have stratifications, and a {\em regular value\/} of a map $g\colon C \rightarrow Y$ is a point $y$ of $Y$ that is a regular value for all the restrictions of $g$ to the faces of $C$.

In order to be regular, a point of $(\overline{S^2_H})^{3n}$ will have to be regular (in our sense) for all the maps $g(\Gamma,L,B)$, but it will also have to be regular (in the same sense) for other maps associated with the definition of pseudo-sections (Definition~\ref{defpseudosec}).
We describe all these maps below.

Let $p_v$ denote the height in $\tilde{\EXT}_{[1.5,1.7]}$ that is the natural projection to the vertical factor $\RR$.
Let $U$ be the disjoint union of a line and the trivial infinite cyclic covering of a circle. Embed an equivariant tubular neighborhood $[a,b]\times U \times [-1,1]$ of $U$ in $\tilde{\EXT}_{[1.5,1.7]}$
so that $p_v$ is constant on $[a,b]\times (k,t)$ for any $(k,t) \in U \times [-1,1]$, and the restriction of $p_v$ to $\{\frac{a+b}2\}\times U \times \{0\}$ is regular on the line and it has only two critical points in every circle: a non degenerate maximum and a degenerate minimum with respective heights $n$ and $(n+1)$.
(This $U$ and the normalization of $\funcE$ in $N(c)$ in Lemma~\ref{lemtransboun} below are only here in order to have a  notion
of regularity in $(S^2_H)^{3n}$ independent of $N(c)$.)

Consider a partition $\CP$ of the subset of vertices $L^{-1}(\EXT_{[1.5,4]})$ into subsets $P_0$, $P_1$,
\dots $P_k$, where every $P_i$ has at least $2$ elements for $i\geq 1$, and $P_0$ may be empty; also consider the configuration subspace ${C}(\Gamma,L,B;\CP)$ of ${C}(\Gamma,L,B)$ made of the configurations that map $P_0$ to $\EXT_{[2,4]}$, and every $P_i$ to some point $m_i$ of $[a,b]\times U \times [-1,1]$.
Assume that $B$ contains all the edges whose vertices lie in different $P_i$ with $i>0$.

Say that an edge $e$ is a $P_i$-edge if both ends of $e$ are in $P_i$, and say that an edge is a $P_{>0}$-edge if $e$ is a $P_i$-edge for some $i>0$. Let $A_{>0}$ be the set of the $P_{>0}$-edges of $A$.

Let $e\in A_{>0}$, then $p(\Gamma,e)$ projects ${C}(\Gamma,L,B;\CP)$ to $U(\EXT_{|[a,b]\times U \times [-1,1]})$, this projection can be naturally composed
to yield two projections $q_1(\Gamma,e)\colon {C}(\Gamma,L,B;\CP) \rightarrow [a,b]$ and $$q_2(\Gamma,e)\colon {C}(\Gamma,L,B;\CP) \rightarrow [-1,1].$$
Then the projection $p(\Gamma,e)$ can be twisted in various ways with the help of the functions of Definition~\ref{defpseudosec}.
Namely, let $\pi(\tau_{\pi}) \colon U\EXT_{|[a,b]\times U \times [-1,1]} \rightarrow S^2$ be the projection induced by the parallelization $\tau_{\pi}$, set
$$p_F(\Gamma,e)(\gamma)=F(q_1(\Gamma,e)(\gamma),q_2(\Gamma,e)(\gamma))(\pi(\tau_{\pi}) \circ p(\Gamma,e)(\gamma)),$$
and, for $j=1,2,3$,
$$p_j(\Gamma,e)(\gamma)=\rho_{(j-2)\alpha(q_2(\Gamma,e)(\gamma))}(\pi(\tau_{\pi}) \circ p(\Gamma,e)(\gamma)).$$

Let $A_1$, $A_2$, $A_C$ and $A_F$ be four disjoint subsets of $A_{>0}$ such that, for any $i$, $A_C$ contains at most one $P_i$-edge, and if $A_C$ contains a $P_i$-edge,
$A_F$ does not. Let $A_3$ denote the complement of the union $A_1\cup A_2 \cup A_C\cup A_F$ in $A_{>0}$. Define
$$g(\Gamma,L,B;\CP,A_1,A_2,A_C,A_F) \colon {C}(\Gamma,L,B;\CP) \times (\overline{S^2_H})^B \rightarrow (\overline{S^2_H})^{\{1,2,\dots,3n\} \setminus A_C}\times [-1/50,1/50]^{A_C}$$ by
$$\begin{array}{ll}g(\Gamma,L,B;\CP,A_1,A_2,A_C,A_F)=&\prod_{e \in (A\setminus A_{>0})} \overline{\pi} \circ p(\Gamma,e)\\
& \times \prod_{e \in A_1} p_1(\Gamma,e) \times \prod_{e \in A_2}p_2(\Gamma,e) \times \prod_{e \in A_3}p_3(\Gamma,e)
\\&\times \prod_{e \in A_F} p_F(\Gamma,e)\times \prod_{e \in A_C} p_v\circ \overline{\pi} \circ p(\Gamma,e) \times \mbox{Identity}(\overline{S^2_H})^B.\end{array}$$
where $p_v \colon \overline{S^2_H} \rightarrow [-1/50,1/50]$ is the first coordinate.

\begin{definition}
 \label{defuplegen}
Let $P_V(X_i)$ be the vertical plane of $\RR^3$ that contains $0$ and $X_i$.
An element $(X_i)_i$ of $(S^2_H)^{3n}$ is {\em regular\/} if 
\begin{itemize}
\item it is a regular value for all the maps $g(\Gamma,L,B)$,
\item its image under $\mbox{Identity}(\overline{S^2_H})^{\{1,2,\dots,3n\} \setminus A_C}\times p^{A_C}_v$ is a regular value of $g(\Gamma,L,B;\CP,A_1,A_2,A_C,A_F)$ for all these maps,
\item for any pair $\{i,j\}$ of $\{1,2,\dots,3n\}$, $P_V(X_i) \cap P_V(X_j)$ is the vertical line through the origin. 
\end{itemize}
\end{definition}

\begin{lemma}
The set of regular $(X_i)_i$ is dense and open in $(\overline{S^2_H})^{3n}$.
\end{lemma}
\bp
The set of regular values of $g(\Gamma,L,B)_{|g(\Gamma,L,B)^{-1}((S_H^2)^{3n})}$ is open and dense in $(S_H^2)^{3n}$.
The density is a direct corollary of the Morse-Sard theorem \cite[Chapter 3, Section 1]{hirsch}. Since
the set where the tangent maps are not all surjective is closed in the compact source of $g(\Gamma,L,B)$, its image is compact, and the set of regular values is open.
Similarly, the sets of regular values of the $g(\Gamma,L,B;\CP,A_1,A_2,A_C,A_F)$ are open and dense.
Then the preimage of such a set under the open map $\mbox{Identity}(\overline{S^2_H})^{\{1,2,\dots,3n\} \setminus A_C}\times p^{A_C}_v$ is also open and dense.
Thus the set of regular $(X_i)_i$ is a finite intersection of dense open subsets of $(\overline{S^2_H})^{3n}$.
\eop

\subsection{Proof of Theorem~\ref{thmdefinvpseudo}}

\begin{lemma}
\label{lemtransboun} Assume that $\funcE$ satisfies the following assumptions $(\ast \ast)$
\begin{itemize}
\item the critical points of the restriction of $\funcE$ to the framed link $c$ of Definition~\ref{defpseudotriv} are non-degenerate, they are minima and maxima $x$ mapped to $\ZZ$ by $\functE$ respectively that have neighborhoods  $[a,b] \times [x-\varepsilon,x+\varepsilon] \times [-1,1]$ that may be naturally identified with neighborhoods of the maxima or the minima of $U$ by maps
$(\alpha,y,t) \mapsto (\alpha,g(y),t)$ so that $\functE(\alpha,y,t)=\functE(\alpha,g(y),t)$,
 \item the restriction of $\functE$ to $[b-\varepsilon,b] \times\{(\gamma,t)\}$ is constant, for any $(\gamma,t) \in  c \times [-1,1]$
\item the map from $[a,b] \times c \times [-1,1]$ that is the product of $\funcE$ and the projections
$p_{[a,b]}$ to $[a,b]$ and $p_{[-1,1]}$ to $[-1,1]$ is a local diffeomorphism except for the neighborhoods of the extrema normalized above.
\end{itemize}
Then when $(X_i)_i$ is regular, the intersection of the $p(\Gamma,i)^{-1}(s_{\tilde{\tau}}(\EXT_{[0,2]};\cvarM_i,\varepsilon_i) + U\EXT_{|{\sK_i \times S^1}} +\overline{\pi^{-1}(\cvarM_i)})$ over any non empty subset $I$ of $\{1,2,\dots,2n\}$ is transverse in $$C_{2n}(\EXT) \cap \cap_{i=1}^n p(\Gamma,i)^{-1}\left( U \EXT_{[0,2[} \cup C_2(\EXT)\setminus C_2(\EXT_{[0,2[})\right).$$
\end{lemma}
\bp 
Consider an intersection point $m$ in such a partial intersection.
For $i \in I$, $p(\Gamma,i)(m)$ could be in $U\EXT_{[0,2[}$ or in $\overline{\pi^{-1}(\cvarM_i)}$.
In both cases, $\left(s_{\tilde{\tau}}(\EXT_{[0,2]};\cvarM_i,\varepsilon_i) + U\EXT_{|{\sK_i \times S^1}} +\overline{\pi^{-1}(\cvarM_i)}\right)$ is locally the preimage of a point under a map $\pi_i$ valued in $S^2$, 
except in the cases listed below.
In order to show transversality, we prove that at the intersection, the tangent map to the product of the constraint maps is onto.
Note that when $p(\Gamma,i)(m)$ is in $U\EXT_{|\sK_i \times S^1}$, this constrains the position of the projection to $\EXT$ and yields a constraint valued in the horizontal $\RR^2$ that is naturally independent from the other ones (moving horizontally a point only changes the value of this constraint map). We won't consider this part of the boundary anymore. Let us now list the bad cases using notation from Definition~\ref{defpseudosec}.
\begin{itemize}
\item When $p(\Gamma,i)(m)$ is in $U\EXT_{|({[a,b] \setminus \{b-\varepsilon_i\}) \times c \times [-1,1]}}$, then $s_{\tilde{\tau}}(\EXT_{[0,2]};\cvarM_i,\varepsilon_i)$ can be a combination of preimages $\pi_{i,j}^{-1}(\cvarM_i)$ of such $\pi_{i,j}$, and 
all preimages involved in these combinations must be transverse to each other.
\item When $p(\Gamma,i)(m)$ in $U\EXT_{|{\{b-\varepsilon_i\} \times c \times [-1,1]}}$, then it must be in $\{b-\varepsilon_i\} \times \{c\} \times C_2(\cvarM_i)$.
\end{itemize}

Note that the only impact of the points that are in $\EXT_{[0,1.8[}$ on the constraint associated
with their edges 
that have another end in $\EXT_{[2,4]}$ factors through $\functE$. In particular, the influence of the points of $[a,b] \times c \times [-1,1]$ factors through $(\functE,p_{[a,b]},p_{[-1,1]})$. Then according to our assumptions, and because of the invariance of $\overline{\pi}$ under vertical global translations, it will be enough to pretend that our link $c$ lifts as $U$. Therefore the regularities of our maps $g(\Gamma,L,B;\CP,A_1,A_2,\emptyset,A_F)$ show transversality when no chain $C_2(\cvarM)$ is involved.
Let us now study the impact of these chains.
With respect to the parallelization $\tau_b$, the boundary of the chain $C_2(\cvarM)\subset [-1,1] \times S^1(\cvarM)$ reads 
$$\{\left(u,\rho_{-\alpha(u)}(\cvarM)\right), u \in[-y,x]\}
+\{\left(u,\rho_{\alpha(u)}(\cvarM)\right), u \in[-y,x]\} - 2[-y,x]\times \{\cvarM\}$$
for any arbitrary $x$ and $y$ in $[1-\varepsilon,1]$.
Then $C_2(\cvarM)$ can be chosen as the union of two embedded triangles $T_1$ and $T_2$ in $[-1,1] \times S^1(\cvarM)$
with respective boundaries $\{\left(u,\rho_{-\alpha(u)}(\cvarM)\right), u \in[-y,x]\}
+\{x\} \times S^1(\cvarM) - [-y,x]\times \{\cvarM\}$ and $\{\left(u,\rho_{\alpha(u)}(\cvarM)\right), u \in[-y,x]\}
-\{x\} \times S^1(\cvarM) - [-y,x]\times \{\cvarM\}$ (Note that when $\cvarM \in S^2_H$, $S^1(\cvarM)$ is not reduced to a point). 
In the interior of the triangles, the constraint reads as the product of $p_{[a,b]} \circ p(\Gamma,e)$ that must go to $(b-\varepsilon_i)$ and $p_v \circ \overline{\pi} \circ p(\Gamma,e)$ that must be the height of the horizontal circle $S^1(\cvarM)$.
Therefore, since $p_{[a,b]} \circ p(\Gamma,e)$ is independent from the constraints that factor through $\functE$ and from the constraints of the other $P_i$-edges, transversality holds there.
On the edges $\{\left(u,\rho_{(2-j)\alpha(u)}(\cvarM)\right), u \in[-y,x]\}$, the constraint map $p_v \circ \overline{\pi} \circ p(\Gamma,e)$ is replaced by $p_j(\Gamma,e)$, and on the edges $\{x\} \times S^1(\cvarM)$ the constraint that $p_{[-1,1]} \circ p(\Gamma,e)$ goes to $x$ must be added to the constraint map $p_v \circ \overline{\pi} \circ p(\Gamma,e)$. On the vertices, a constraint that $p_{[-1,1]} \circ p(\Gamma,e)$ goes to $x$ or $-y$ must be added to the constraint map $p_2(\Gamma,e)$.
Thus when $x$ and $y$ are regular values of finitely many well-chosen functions on compact spaces, transversality holds.
\eop

\noindent{\sc Proof of Theorem~\ref{thmdefinvpseudo}:}
According to Lemmas~\ref{lemgenpos} and~\ref{lemdiffz}, there is no homological obstruction to the existence of $C_i=C(\cvarM_i,\sK_i,\tilde{\tau})$. 
The existence of the chains $C_i$ of $\tilde{C}_2(\EXT_{[0,2]})$ with the prescribed boundaries in general $3n$-position is proved as in \cite[Section 2.5]{lesbonn} once we have transversality on the prescribed faces, as we have, thanks to Lemma~\ref{lemtransboun}.
Now, it suffices to find forms $\Omega_i$ Poincar\'e dual to our given $\deleqprop_i(\tilde{\tau})=\delta(K)(t)\eqprop_i(\tilde{\tau})$, such that the supports of the $\Omega_i$ are small enough neighborhoods of the $\deleqprop_i(\tilde{\tau})$ so that $I_{\Gamma}(\{\eqprop_i(\tilde{\tau})\}_i)=I_{\Gamma}((\Omega_i)_i)$ for any $\Gamma$ in  $\CS_n^u$, and the $\Omega_i$ satisfy the hypotheses of Theorem~\ref{thmdefzformpseudo}. Again, we let the supports of the $\omega_i$ be small enough disks around the points $\cvarM_i$ associated with the $\{\eqprop_i(\tilde{\tau})\}_i$, and the only additional thing to see there is that
our $\Omega_i$ can be assumed to have the form prescribed by Theorem~\ref{thmdefzformpseudo} 
(and Definition~\ref{defpseudoform}) around $$U(N_i(k)=[b-\varepsilon_i-\varepsilon_i^k,b-\varepsilon_i+\varepsilon_i^k]  \times c \times ]-1,1[)=_{\tau_b} N_i(k) \times S^2.$$

The chain $\frac{1}{\delta(K)(t)}C_1$ meets $\ZZ \times N_1(k) \times S^2$ in $\{0\} \times N_1(k) \times S^1(X_1)$.

Let $p_v \colon S^2 \rightarrow [-1,1]$ be the height map, and let $N_k(X_1)$ be a thin neighborhood 
$p_v^{-1}( ]p_v(X_1)-\varepsilon_1^k,p_v(X_1)+\varepsilon_1^k[)$ of $S^1(X_1)=p_v^{-1}( \{p_v(X_1)\})$.

By transversality, we may assume that $p(\Gamma,1)$ projects 
$\bigcap_{j\neq 1} p(\Gamma,j)^{-1}(p_C(\deleqprop_j(\tilde{\tau})))$ 
outside $\{b-\varepsilon_1\}\times c \times ]-1,1[ \times S^1(X_1) \subset \partial C_2(E)$.
Then when $k$ is big enough, for any $\Gamma$, $p(\Gamma,1)$ projects   
$\bigcap_{j\neq 1} p(\Gamma,j)^{-1}(p_C(\deleqprop_j(\tilde{\tau})))$ 
outside a neighborhood of $N_1(k) \times N_k(X_1)$ in $C_2(E)$.
Assume that $\omega_1$ is supported in $N_k(X_1)$ and that $\Omega_1$ has the prescribed form
on $\partial \tilde{C}_2(E) \setminus (\{0\} \times N_1(k) \times S^2)$.

\begin{sublemma}
\label{subexppseud}
There exists an explicit closed form of $(N_1(k) \times S^2)$ supported on $N_1(k) \times N_k(X_1)$ that extends $\Omega_1$ to a form that satisfies the hypotheses of Theorem~\ref{thmdefzformpseudo}.
\end{sublemma}
\noindent {\sc Proof of Sublemma~\ref{subexppseud}:}
Let $g_s \colon ]-1/50,1/50[ \rightarrow \RR$ be a function with compact support in $]p_v(X_1)-\varepsilon_1^{k+1},p_v(X_1)+\varepsilon_1^{k+1}[ \cap ]-1/50,1/50[$ such that $\int_{-1/50}^{1/50} g_s =\frac{1}{2\pi}$, let $\theta$ be the angular coordinate around the vertical axis (the argument in the horizontal $\CC$) and let
$\omega_s = (g_s\circ p_v) d\theta \wedge d p_v$. Assume that 
$\omega_1=(g_s\circ p_v)(g_{\theta} \circ \theta) d\theta \wedge d p_v$ for a function 
$g_{\theta}$ on $[0,2\pi]/(0\sim 2\pi)=\RR/(x\sim x+2\pi)$ such that $\int_0^{2\pi}g_{\theta}d\theta=2\pi$.
Set $$G_{\theta}(\theta)=\int_0^{\theta}(g_{\theta}-1)d\theta.$$
Set $\eta_{1,s,1}=(g_s\circ p_v)G_{\theta}d p_v$.
Then $\omega_1=\omega_s + d\eta_{1,s,1}$.
Define a smooth function $\chi \colon [a,b] \rightarrow [0,1]$ that maps $[a,b-\varepsilon_1-\varepsilon_1^k+\varepsilon_1^{k+1}]$ to $1$ and $[b-\varepsilon_1+\varepsilon_1^k-\varepsilon_1^{k+1},b]$ to $0$.
Assume that the map $\alpha$ of Definition~\ref{defpseudotriv} is constant around $0$.
For $(t,\gamma,u ;\cvarM) \in [a,b] \times c \times [-1,1]\times S^2_H$, define $$G_{\theta,\chi}(t,\gamma,u ;\cvarM)=\left\{ \begin{array}{ll}\frac{G_{\theta}(\theta(\cvarM)+\chi(t)\alpha(u)) +G_{\theta}(\theta(\cvarM)-\chi(t)\alpha(u))}{2} & \mbox{if}\;u\leq 0\\
\frac{G_{\theta}(\theta(\cvarM)+\chi(t)(\alpha(u)-2\pi)) +G_{\theta}(\theta(\cvarM)-\chi(t)(\alpha(u)-2\pi))}{2} & \mbox{if}\;u\geq 0 \end{array}\right.$$ and 
$\eta_{1,s}=G_{\theta,\chi}\pi(\tau_b)^{\ast}((g_s\circ p_v)d p_v)$.
Then $\eta_{1,s}$ and $\eta_{1,s,1}$ satisfy the conditions of Definition~\ref{defpseudoform}, and the restriction of the form
$\omega(\omega_1,k,\varepsilon_1,\eta_{1,s})$ constructed as in this definition to $N_1(k) \times S^2$ is supported in $N_1(k) \times N_k(X_1)$.
\eop

Thus we have two closed $2$-forms on $\{0\} \times N_1(k) \times \overline{N_k(X_1)} \subset \partial \tilde{C}_2(\EXT)$ that coincide near the boundary: the former $\tilde{\Omega}_1$ Poincar\'e dual to our given $\deleqprop_1(\tilde{\tau})$ and the explicit extension of Sublemma~\ref{subexppseud}.
Since $H^2(N_1(k) \times \overline{N_k(X_1)}, \partial (N_1(k) \times \overline{N_k(X_1)}))=0$, their difference has a primitive $\eta$ that vanishes on the boundary that can be used to interpolate them in a collar neighborhood
 $[-1,0] \times \{0\} \times  N_1(k) \times \overline{N_k(X_1)}$ in $\tilde{C}_2(\EXT)$ as $\tilde{\Omega}_1 +d(\chi \eta)$.
Thus our $\tilde{\Omega}_1$ may be modified in a neighborhood of $\{0\} \times N_1(k) \times \overline{N_k(X_1)}$ in $\tilde{C}_2(\EXT)$ that does not meet $p_C^{-1}\left(p(\Gamma,1)\left(\bigcap_{j\neq 1} p(\Gamma,j)^{-1}(p_C(\deleqprop_j(\tilde{\tau})))\right)\right)$ so that it satisfies the hypotheses of 
Theorem~\ref{thmdefzformpseudo} without changing the $I_{\Gamma}((\Omega_i)_i)$.
Applying the same process for each $i$ for pairwise disjoint $[b-\varepsilon_i-2\varepsilon_i^{k_i},b-\varepsilon_i+2\varepsilon_i^{k_i}]$ yields the result.

\eop

\subsection{Other expressions of $\ansothree(A,\tilde{\tau}_0,\tilde{\tau}_1)$ }
\label{subxiint}
Let $A$ be a $\QQ$-handlebody.

When $3n$ $4$-cycles $H_i$ of $([0,1]\times UA, \partial ([0,1]\times UA))$ are given, we can define their connected $\CA_n$-intersection
$$I^c_{n,[0,1]\times S_{2n}(TA)}(\{H_i\})=\sum_{\Gamma \in \CS_n}\frac{\langle p(\Gamma,1)^{-1}(H_1), p(\Gamma,2)^{-1}(H_2), \dots, p(\Gamma,3n)^{-1}(H_{3n})
\rangle_{[0,1]\times S_{2n}(TA)}}{2^{3n}(3n)!(2n)!}[\Gamma] \in \CA_n$$ as soon as they are in general $3n$-position.
Similarly, when a $4$--chain $H_k$ and $(3n-1)$ $3$--cycles $S_i$ of $(UA,UA_{|\partial A})$, $i \in \{1,2,\dots,3n\}\setminus \{k\}$ are given, their connected $\CA_n$-intersection is
$$I^c_{n,S_{2n}(TA)}(H_k,\{S_i\}_{i \neq k})=\sum_{\Gamma \in \CS_n}\frac{\langle p(\Gamma,i)^{-1}(S_i)_{i \in \{1,2,\dots,3n\}\setminus \{k\}}, p(\Gamma,k)^{-1}(H_k)
\rangle_{S_{2n}(TA)}}{2^{3n}(3n)!(2n)!}[\Gamma] \in \CA_n$$ as soon as they are in general $3n$-position (i.e. as soon as this intersection is transverse for all $\Gamma$). 

Recall that $$\tilde{z}_n(K,\tilde{\tau}_1)-\tilde{z}_n(K,\tilde{\tau})=\ansothree_n(\EXT,\tilde{\tau},\tilde{\tau}_1)=\frac{p_1(\tilde{\tau},\tilde{\tau}_1)}{4}\ansothree_n$$ according to Theorems~\ref{thmdefzformpseudo} and \ref{thmxin}. The following proposition gives a definition of $p_1(\tilde{\tau},\tilde{\tau}_1)\ansothree_n$ based on intersections.

\begin{proposition}
\label{propdiffz}
Let $A$ be a $\QQ$-handlebody equipped with two pseudo-parallelizations $\tilde{\tau}$ and $\tilde{\tau}_1$ that coincide near the boundary of $A$. Let $\cvarM \in S^2$.
Pick $3n$ distinct $\varepsilon_i$ in $]0,\varepsilon[$, and let $(\cvarM_1, \dots,\cvarM_{3n})$ be a regular $3n$-tuple of $(S^2_H)^{3n}$. There exist $3n$ $4$--chains $\eqprop_i=\eqprop(\tilde{\tau},\tilde{\tau}_1,\cvarM_i,\varepsilon_i)$ in general $3n$-position in $[0,1] \times UA$ with boundary 
$$\partial(\tilde{\tau},\tilde{\tau}_1,\cvarM_i,\varepsilon_i)= -\left(\{0\} \times s_{\tilde{\tau}}(A;\cvarM_i,\varepsilon_i)\right) - \left([0,1] \times s_{\tilde{\tau}}(\partial A;\cvarM_i,\varepsilon_i)\right) +\left(\{1\} \times s_{\tilde{\tau}_1}(A;\cvarM_i,\varepsilon_i)\right),$$
and
$$I_n^c(\{\eqprop_i\}_{i \in \{1,2,\dots,3n\}})_{[0,1] \times UA}=\frac{p_1(\tilde{\tau},\tilde{\tau}_1)}{4}\ansothree_n.$$
\end{proposition}
\bp
Lemma~\ref{lemdiffz} guarantees that there is no homological obstruction to the existence of $\eqprop_i$. The $4$--chain $\eqprop_i$ can be assumed to be transverse to the boundary. Then this proposition is dual to  Proposition~\ref{propxinzform} and its proof is similar to the proof of Theorem~\ref{thmdefinvpseudo}. 
\eop

We now give a last expression for $p_1(\tilde{\tau},\tilde{\tau}_1)\ansothree_n$ similar to the expression in Proposition~\ref{propxinzformbis}.

\begin{proposition}
\label{propvarpseudo}
Under the assumptions of Proposition~\ref{propdiffz}, there exist a regular $(\cvarM_i)_i$ of $(S^2_H)^{3n}$
and $3n$ $4$--chains $H_i=H(\tilde{\tau},\tilde{\tau}_1,\cvarM_i,\varepsilon_i)$ with compact support in $U(\mathring{A})$, such that $$\partial H_i=s_{\tilde{\tau}_1}(A;\cvarM_i,\varepsilon_i) - s_{\tilde{\tau}}(A;\cvarM_i,\varepsilon_i)$$
and, for any $k$, 
$$\left\{\begin{array}{l}p(\Gamma,k)^{-1}(H_k),\\
\mbox{the}\;p(\Gamma,j)^{-1}(s_{\tilde{\tau}_1}(A;\cvarM_j,\varepsilon_j))\;\mbox{for}\;1\leq j<k\;\mbox{, and}\\
\mbox{the}\;p(\Gamma,j)^{-1}(s_{\tilde{\tau}}(A;\cvarM_j,\varepsilon_j)) \;\mbox{for}\;k< j\leq 3n\end{array}\right.$$
are in general $3n$-position. 
Then $$\frac{p_1(\tilde{\tau},\tilde{\tau}_1)}{4}\ansothree_n=\sum_{k=1}^{3n}I^c_{n,S_{2n}(TA)}(H_k,\{s_{\tilde{\tau}_1}(A;\cvarM_j,\varepsilon_j)\}_{j < k},\{s_{\tilde{\tau}}(A;\cvarM_j,\varepsilon_j)\}_{j > k}).$$
\end{proposition}
\eop

\newpage
\section{Additivity of $\tilde{z}$ under connected sum}
\label{secconnsum}
\setcounter{equation}{0}

Let $\Kone$ and $\Ktwo$ be two null-homologous knots in two $\QQ$-spheres $\Rone$ and $\Rtwo$.
Remove an open $3$-ball met by the knot along a diameter from $(\Rone, \Kone)$ and from $(\Rtwo, \Ktwo)$. Glue the obtained pairs along their spherical boundaries so that the end points of the cut knots match, and the cut knots together form a knot $\Kone \sharp \Ktwo$ whose orientation agrees with the orientations of $\Kone$ and $\Ktwo$.
The knot $\Kone \sharp \Ktwo$ is a null-homologous knot in the $\QQ$-sphere $\Rone \sharp \Rtwo$.
It is called the \emph{connected sum} of $\Kone$ and $\Ktwo$.

This section is devoted to the proof of the following theorem.

\begin{theorem}
 \label{thmconnsum}
Under the above assumptions, with the notation of Theorem~\ref{thmhighloop},
$$\tilde{z}(\Kone \sharp \Ktwo) =\tilde{z}(\Kone) + \tilde{z}(\Ktwo).$$
\end{theorem}

\subsection{Sketch of the proof}
\label{subskconn}
Consider the disks $D(-\frac14,\frac18)$ and $D(\frac14,\frac18)$ in $D_{[0,4]} \subset \CC$ with common radius $\frac18$ and with respective centers $-\frac14$ and $\frac14$.
Identify the boundary of the exterior $\EXTone$ of $\Kone$ (resp. of the exterior $\EXTtwo$ of $\Ktwo$) with $\partial D(-\frac14,\frac18) \times S^1$ (resp. $\partial D(\frac14,\frac18) \times S^1$) via an identification that has the same properties as the identification of $\partial \EXT$ with $\partial D_{[0,4]} \times S^1$ in Subsection~\ref{subcyc}, up to a translation and a dilation of the horizontal plane.
Precisely, when $z \in S^1$, the curve $\partial D(\pm 1/4,1/8) \times \{z\}$ is a parallel of $K^{\pm}$ that bounds in $E^{\pm}$, and for $x \in \partial D(\pm 1/4,1/8)$, $\{x\} \times S^1$ is a meridian of $K^{\pm}$ in $R^{\pm}$.

Let $\BOF=D(-\frac14,\frac18) \sqcup D(\frac14,\frac18)$.
Then the exterior $\EXTonestwo$ of $\Kone \sharp \Ktwo$ reads 
$$\EXTonestwo=\left( (D_{[0,4]} \setminus \mathring{D}^{+, -} ) \times S^1 \right) \cup_{\partial \BOF \times S^1} (\EXTone \sqcup \EXTtwo).$$

(It is homeomorphic to the union of $\EXTone$ and $\EXTtwo$ glued along an annulus of their boundaries, where the annulus is both a regular neighborhood of a meridian of $\Kone$ in $\partial \EXTone$ and of $\Ktwo$ in $\partial \EXTtwo$.)

The exteriors $\EXTone$ and $\EXTtwo$ are equipped with parallelizations $\tauone$ and $\tautwo$ that coincide with the standard $\tau_{\pi}$ near their boundaries so that $\tauone$ and $\tautwo$ induce a parallelization $\tau$ of $\EXTonestwo$ that coincides with $\tau_{\pi}$ on $(D_{[0,4]} \setminus \BOF ) \times S^1$.

With the definitions of $p_1(\tau_{\pi},.)$ in Subsection~\ref{subpont} and with Proposition~\ref{proppont}, it is easy to see that
$p_1(\tau_{\pi},\tau)= p_1(\tau_{\pi},\tauone) + p_1(\tau_{\pi},\tautwo)$.
Thus according to Theorem~\ref{thmhighloop}, in order to prove Theorem~\ref{thmconnsum}, it suffices to prove that 
$$\tilde{z}_n(\Kone \sharp \Ktwo,\tau)=\tilde{z}_n(\Kone,\tauone)+ \tilde{z}_n(\Ktwo,\tautwo) $$
for our arbitrary fixed trivializations, for any $n \in \NN$.

Very roughly, the idea of the proof is the following one.
Varying the radius of the disks $D(-\frac14,\frac18)$ and  $D(\frac14,\frac18)$ to let it approch $0$ concentrates
$\EXTone$ and $\EXTtwo$ near the circles $\Cone=\{-1/4\} \times S^1$ and $\Ctwo=\{1/4\} \times S^1$.  Thus $\EXTonestwo$ may be thought of as the result of an undefined generalized blow-up of $\Cone$ and $\Ctwo$ that are respectively replaced by $\EXTone$ and $\EXTtwo$. 
The $S^1$-coordinate of a limit point on $\Cone$ or $\Ctwo$ is the value of $\funcE$, and it is the only information that is seen from the exterior of the circle. Thus a limit configuration may be thought of as a large scale configuration in $D^2 \times S^1$ equipped with $\Cone$ and $\Ctwo$ with possibly colliding points
on $\Cone$ or $\Ctwo$, and configurations in $\EXTone$ and $\EXTtwo$ of the points on $\Cone$ and $\Ctwo$ whose vertices $\funcE$-values 
are given by the $S^1$-coordinates on $\Cone$ or $\Ctwo$ of the large scale configuration.
Now, the form corresponding to the edges that do not connect vertices inside a circle $\Cone$ or $\Ctwo$ only depends on the large scale configuration and factors through global $S^1$--translation of this large scale configuration.
Thus unless all the edges connect vertices of the same circle, the form factors through a configuration space of smaller dimension. Thus the vertices of a contributing connected graph must all lie either on $\Cone$ (i.e. in $\EXTone$) or on $\Ctwo$ (i.e. in $\EXTtwo$), and these contributions add up to $\tilde{z}_n(\Kone \sharp \Ktwo,\tau)$, for a given $n$.

In the proof below, we follow the idea above without attempting to define the mentioned blow-ups

\subsection{The proof}
\label{subpfconn}

Set $D^-_{[0,\frac{3}{2}]}=\{z \in D_{[0,\frac{3}{2}]}; \mbox{Re}(z) <0\}$,
 $D^+_{[0,\frac{3}{2}]}=\{z \in D_{[0,\frac{3}{2}]}; \mbox{Re}(z) >0\}$ and $\EXTonetwo=\EXTone \sqcup \EXTtwo$.
Change the map $r$ of Subsection~\ref{subcyc} to a homotopic map $r_{\sharp}\colon \EXT \rightarrow [\frac12,4]$ that is unchanged on $r^{-1}([\frac{3}{2},4])$ such that $r^{-1}(]\frac{3}{2},4])=r_{\sharp}^{-1}(]\frac{3}{2},4])$,
$r_{\sharp}^{-1}(\frac12)= \EXTonetwo$
and $r_{\sharp}$ maps a neighborhood of the straight segment from $-\frac{3i}{2}$ to $\frac{3i}{2}$ to $\frac{3}{2}$.
Since the values of $r$ on $[0,\frac{3}{2}]$ do not affect the definition of $\tilde{z}_n$, $\pi$ can be defined with $r_{\sharp}$ instead of $r$ without whanging $\tilde{z}_n$. For $I \subset [0,\frac{3}{2}]$, set $D^-_{I}=D^-_{[0,\frac{3}{2}]} \cap r_{\sharp}^{-1}(I)$ and $D^+_{I}=D^+_{[0,\frac{3}{2}]} \cap r_{\sharp}^{-1}(I)$.

According to Theorem~\ref{thmdefzform} and its proof, the map $p(\sK_i)$ defined before Lemma~\ref{lemOmegai} may be changed on $\{0\}\times U\EXT \subset \partial \tilde{C}_2(\EXT)$ to any homotopic map that maps the complement of the neighborhood of an $S^1$ fiber in $\{0\}\times U((D_{[0,\frac{3}{2}]} \setminus \BOF ) \times S^1)$ to the North Pole, and that factors through the projection to $D_{[0,\frac{3}{2}]} \setminus \BOF$ around this fiber.

For a one form $\eta$ supported in the interior of $\{0\}\times U((D_{[0,\frac{3}{2}]} \setminus \BOF ) \times S^1)$, that factors through the projection to $D_{[0,\frac{3}{2}]} \setminus \BOF$, one may similarly add $P(\thetap^{\ast})d\eta$ to the restriction to the boundary of the form $\Omega_i$
of Lemma~\ref{lemOmegai}, for any polynomial $P$.

Since $\Delta(\Kone \sharp \Ktwo)=\Delta(\Kone)\Delta(\Ktwo)$,
$$\frac{\Delta(K)^{\prime}}{\Delta(K)}=\frac{\Delta(\Kone)^{\prime}}{\Delta(\Kone)}+\frac{\Delta(\Ktwo)^{\prime}}{\Delta(\Ktwo)},$$
and 
$(\delta(K)J_{\Delta(K)})(\thetap^{\ast})(p(\sK_i)^{\ast}(\omega_{i,D}))$ may be changed to

$$(\delta(K)J_{\Delta(\Kone)})(\thetap^{\ast})(p(\sK_i^-)^{\ast}(\omega_{i,D})) + (\delta(K)J_{\Delta(\Ktwo)})(\thetap^{\ast})(p(\sK_i^+)^{\ast}(\omega_{i,D}))$$
where $p(\sK_i^-)$ is supported in $D^-_{[0,1[} \times S^1 \times S^2$, near a fiber $\{x\} \times S^1 \times S^2$ and 
$p(\sK_i^+)$ is supported in $D^+_{[0,1[} \times S^1 \times S^2$, near another fiber.
Furthermore assume that $p(\sK_i^-)$ and $p(\sK_i^+)$ are supported in $r_{\sharp}^{-1}(]\frac{5}{8},\frac34[)$.

Now, extend the map $\pi(\tau)$ defined in the beginning of Section~\ref{secaltdefzform} with notation of Subsection~\ref{subcyc}, constructed with $r_{\sharp}$ instead of $r$ as follows.
Set $$\EXT^-_{[0,\frac{3}{2}]} = \left((D^-_{[0,\frac{3}{2}]}\setminus \mathring{D}(-\frac14,\frac18)) \times S^1\right) \cup \EXTone, \;\;\;\EXT^+_{[0,\frac{3}{2}]} = \left((D^+_{[0,\frac{3}{2}]}\setminus \mathring{D}(\frac14,\frac18)) \times S^1\right) \cup \EXTtwo,$$
and for $I \subset [0,\frac{3}{2}]$, set $\EXT^-_{I}=\EXT^-_{[0,\frac{3}{2}]} \cap r_{\sharp}^{-1}(I)$ and $\EXT^+_{I}=\EXT^+_{[0,\frac{3}{2}]} \cap r_{\sharp}^{-1}(I)$.

Define $\pi(\tau)$ and extend it with the same formula on $$\tilde{C}_2(\EXT) \setminus \left(\tilde{C}_2(r_{\sharp}^{-1}[0,\frac34[) \cup p_C^{-1}\left((\EXT^-_{[0,\frac{3}{2}[} \times \EXT^+_{[0,\frac{3}{2}[}) \cup (\EXT^+_{[0,\frac{3}{2}[} \times \EXT^-_{[0,\frac{3}{2}[})\right) \right).$$

On $p_C^{-1}\left((\EXT^-_{[0,\frac{3}{2}[} \times \EXT^+_{[0,\frac{3}{2}[}) \cup (\EXT^+_{[0,\frac{3}{2}[} \times \EXT^-_{[0,\frac{3}{2}[})\right)$, extend $\pi(\tau)$ by using
$$U^{-,+}(u,v)=(1-\chi(r_{\sharp}(u) -\frac{3}{2}))(0,\functE(u)) + \chi(r_{\sharp}(u)-\frac{3}{2})u$$
$$V^{-,+}(u,v)=(1-\chi(r_{\sharp}(v)-\frac{3}{2}))(0,\functE(v)) + \chi(r_{\sharp}(v)-\frac{3}{2})v$$
for $(u,v) \in (\tilde{\EXT}^-_{[0,\frac{3}{2}[} \times \tilde{\EXT}^+_{[0,\frac{3}{2}[}) \cup (\tilde{\EXT}^+_{[0,\frac{3}{2}[} \times \tilde{\EXT}^-_{[0,\frac{3}{2}[})$,
instead of $U(u,v)$ and $V(u,v)$.

The forms that read 
$$\delta(K)(\thetap^{\ast})\pi(\tau)^{\ast}(\omega_i) + (\delta(K)J_{\Delta(\Kone)})(\thetap^{\ast})(p(\sK_i^-)^{\ast}(\omega_{i,D})) + (\delta(K)J_{\Delta(\Ktwo)})(\thetap^{\ast})(p(\sK_i^+)^{\ast}(\omega_{i,D}))$$
on $\partial \tilde{C}_2(\EXT) \cup \left(\tilde{C}_2(\EXT) \setminus (\tilde{C}_2(\EXT^-_{[0,\frac34]}) \cup \tilde{C}_2(\EXT^+_{[0,\frac34]})) \right)$
extend to $\tilde{C}_2(\EXT^-_{[0,\frac34]})$ and $\tilde{C}_2(\EXT^+_{[0,\frac34]})$ as in Lemma~\ref{lemOmegai}.

Under these assumptions, the configurations that map at least one point in $r_{\sharp}^{-1}(]1,4])$ do not contribute:

Indeed, choose such a point $u$ with maximal $r_{\sharp}(u) \geq 1$. Fix it as a root of an oriented subgraph $G$ and 
add edges like in the proof of Lemma~\ref{lempfbordext} so that the form associated with the added edges of $G$ and the adjacent ones factors through the quotient by global vertical translations of the positions of the vertex of $G$ and the heights of the adjacent vertices, and use the same argument, with $\varepsilon \in ]0,\frac1{16n}[$.

Furthermore assume that the supports of the $\omega_i$ avoid the vertical plane that contains the horizontal real line. Then an edge between $\EXT^+_{[0,\frac{3}{2}-2\varepsilon[}$ and $\EXT^-_{[0,\frac{3}{2}-2\varepsilon[}$ contributes to zero.

Then contributing configurations of connected graphs map all the vertices to $\EXT^+_{[0,1]}$ or to $\EXT^-_{[0,1]}$.

Furthermore, when $\EXT^+_{[0,1]}$ is a standard torus, the above vertical translation argument also applies to configurations of $\EXT^+_{[0,1]}$ (equipped with natural propagator extensions that are invariant under vertical translation) that therefore do not contribute.
This shows, that in this case, the contributions of the configurations of $\EXT^-_{[0,1]}$
add up to $\tilde{z}_n(\Kone,\tauone)$. Similarly, the contributions of the configurations of $\EXT^+_{[0,1]}$
add up to $\tilde{z}_n(\Ktwo,\tautwo)$ so that
$$\tilde{z}_n(\Kone \sharp \Ktwo,\tau)=\tilde{z}_n(\Kone,\tauone)+\tilde{z}_n(\Ktwo,\tautwo).$$

\eop

\def\cprime{$'$}
\providecommand{\bysame}{\leavevmode ---\ }
\providecommand{\og}{``}
\providecommand{\fg}{''}
\providecommand{\smfandname}{\&}
\providecommand{\smfedsname}{\'eds.}
\providecommand{\smfedname}{\'ed.}
\providecommand{\smfmastersthesisname}{M\'emoire}
\providecommand{\smfphdthesisname}{Th\`ese}

\end{document}